\documentclass[11pt,a4paper]{article}
\usepackage[T1]{fontenc}
\usepackage[utf8]{inputenc}
\usepackage{lmodern}
\usepackage{amsfonts,amssymb,amsmath,color}
\usepackage[margin=25mm]{geometry}
\usepackage{microtype}
\usepackage[hidelinks,breaklinks]{hyperref}
\hypersetup{pdftitle={Hadamard Local Well-Posedness for Compressible Liquids with a Free Surface},pdfauthor={Tao Luo and Ruixi Zhang}}
\allowdisplaybreaks[1]
\numberwithin{equation}{section}
\newtheorem{definition}{Definition}[section]
\newtheorem{proposition}[definition]{Proposition}
\newtheorem{theorem}[definition]{Theorem}
\newtheorem{lemma}[definition]{Lemma}
\newtheorem{corollary}[definition]{Corollary}
\DeclareMathOperator{\supp}{supp}

\DeclareMathOperator{\divg}{div}
\DeclareMathOperator{\curl}{curl}

\newcommand{\proofheading}[1]{\par\smallskip\noindent\textit{#1.}\ }
\newcommand{\proofqed}{\unskip\nobreak\hfill\ensuremath{\Box}\par\medskip}
\title{Hadamard Local Well-Posedness for Compressible Liquids with a Free Surface}
\author{Tao Luo\thanks{Department of Mathematics, City University of Hong Kong. Email: \href{mailto:taoluo@cityu.edu.hk}{\texttt{taoluo@cityu.edu.hk}}.}
\and Ruixi Zhang\thanks{Department of Mathematical Sciences, Tsinghua University. Email: \href{mailto:zhangrx24@mails.tsinghua.edu.cn}{\texttt{zhangrx24@mails.tsinghua.edu.cn}}. Corresponding author.}}
\date{}
\begin{document}
\maketitle

\begin{abstract}
We prove local well-posedness in the Hadamard sense for the three-dimensional
compressible Euler equations governing a liquid with a free surface and no
surface tension. Under the Taylor sign condition
and the boundary compatibility conditions, we establish existence,
uniqueness, and strong continuous dependence in a hybrid Sobolev class for
every $s>3$. The velocity and interface have $H^s$ regularity, while the
enthalpy and velocity divergence belong to $H^{s+1/2}$ and $H^{s-1/2}$,
respectively. This separation reflects the coupling of a nondegenerate
acoustic equation with the free-surface and vorticity dynamics. The proof
uses regularization operators that smooth the moving domain and extend the
nonlinear compatibility hierarchy, together with $L^2$ and partial $H^2$
distance estimates for solutions on different domains. Regular solutions
are constructed by a regularized forward Euler scheme. High-order energy
estimates and frequency envelopes then yield strong convergence of smooth
approximations and continuity of the data-to-solution map. Integer Sobolev
indices are treated using a Lions--Magenes endpoint compatibility condition.
\end{abstract}


\section{Introduction}

\subsection{The free-boundary problem}

In this paper we study the three-dimensional compressible Euler equations
with a free boundary and without surface tension. Let $\Omega_t$ be the
bounded domain occupied by the fluid at time $t$ and let
$\Gamma_t=\partial\Omega_t$. The density $\rho$ and velocity $u$ satisfy
\[
\left\{
\begin{aligned}
 &(\partial_t+u\cdot\nabla)\rho+\rho\divg u=0,
     &&\text{in }\Omega_t,\\
 &\rho(\partial_t+u\cdot\nabla)u+\nabla P(\rho)=0,
     &&\text{in }\Omega_t.
\end{aligned}
\right.
\]
We consider the linear pressure law
\[
 P(\rho)=C_v\rho,
\]
where $C_v>0$. This choice is made only for convenience and simplicity
of presentation. The main ideas and techniques developed here also apply
to general sufficiently smooth pressure laws $P=P(\rho)$ with
$P'(\rho)>0$ on the relevant density range, with the corresponding
modifications to the enthalpy formulation and the variable-coefficient
energy and compatibility estimates.

The free boundary moves with the fluid velocity, namely,
\begin{equation}\label{eq:kinbc}
 D_t:=\partial_t+u\cdot\nabla
 \quad\text{is tangential to }\bigcup_t\{t\}\times\Gamma_t.
\end{equation}
Thus its normal velocity is $u\cdot n_{\Gamma_t}$, where $n_{\Gamma_t}$ is
the outward unit normal. We also assume that the pressure is constant on
$\Gamma_t$. Normalize the boundary density to one and take $C_v=1$.
Setting $h=\log\rho$, we can write the equations as
\begin{equation}\label{eq:CE}
\left\{
\begin{aligned}
 &D_th+\divg u=0, &&\text{in }\Omega_t,\\
 &D_tu+\nabla h=0, &&\text{in }\Omega_t,
\end{aligned}
\right.
\end{equation}
with \eqref{eq:kinbc} and
\begin{equation}\label{eq:vacbc}
 h=0\qquad\text{on }\Gamma_t.
\end{equation}
Here $\rho=1$ on $\Gamma_t$, so the density is positive at the boundary.
The associated wave equation is therefore nondegenerate up to the boundary.
This is the liquid case; in the physical-vacuum problem the boundary density
vanishes.

We assume the Taylor sign condition on the initial free boundary:
\begin{equation}\label{eq:TS}
 a_0:=-\nabla h_0\cdot n_{\Gamma_0}\geq c>0
 \qquad\text{on }\Gamma_0.
\end{equation}
The Taylor coefficient $a=-\nabla h\cdot n_{\Gamma_t}$ appears in the
boundary part of the linearized energy, and its positivity is needed for
coercivity. A negative Taylor coefficient leads to Rayleigh--Taylor
instability. For the incompressible free-boundary Euler equations, the
ill-posedness in the absence of the stability condition was proved in
\cite{Ebin}.

The equations contain both acoustic waves and the motion of the free
surface. Applying $D_t$ to the first equation in \eqref{eq:CE}, we obtain
\[
 D_t^2h-\Delta h=\operatorname{tr}\bigl((\nabla u)^2\bigr).
\]
Thus $h$ and $\divg u=-D_th$ satisfy an acoustic equation with a Dirichlet
boundary condition. For the boundary normal we have
\[
 D_tn_{\Gamma_t}=-\nabla^\top u\cdot n_{\Gamma_t}
 \qquad\text{on }\Gamma_t.
\]
The equations for the free surface, tangential derivatives of the velocity,
and vorticity are similar to those in the incompressible problem. We will
refer to $(h,\divg u)$ as the compressible component, and to the free
surface, tangential velocity, and vorticity as the incompressible component.
The two components are coupled through \eqref{eq:CE}.

Their different regularity requirements can also be seen from scaling.
The interior equations are invariant under
\[
 h^\lambda(t,x)=h(\lambda t,\lambda x),\qquad
 u^\lambda(t,x)=u(\lambda t,\lambda x).
\]
The characteristic scaling for the free-surface motion is
\[
\begin{aligned}
 a^\lambda(t,x)&=a(\lambda^{1/2}t,\lambda x),\\
 u^\lambda(t,x)&=\lambda^{-1/2}u(\lambda^{1/2}t,\lambda x),\qquad
 \Gamma_t^\lambda=\lambda^{-1}\Gamma_{\lambda^{1/2}t}.
\end{aligned}
\]
These are the scalings of the two parts of the problem, rather than a
common scaling of the full system. In the spaces used below, the
compressible component has an additional half derivative.

\subsection{Function spaces and the main result}

Let $\Omega_*$ be a bounded, connected domain with smooth boundary
$\Gamma_*$. We consider free surfaces that are normal graphs over
$\Gamma_*$ and write
\begin{equation}\label{eq:graphparameter}
 \Gamma=(\operatorname{Id}+\eta_\Gamma n_{\Gamma_*})\Gamma_*.
\end{equation}
For $\alpha\in[0,1]$ and $\rho>0$, let
$\Lambda(\Gamma_*,\alpha,\rho)$ be the class of such surfaces with
\[
 \|\eta_\Gamma\|_{C^{1,\alpha}(\Gamma_*)}<\rho.
\]
Fix small positive constants $\delta,\delta_0$, put
$\Lambda_*:=\Lambda(\Gamma_*,\frac\delta2,\delta_0)$, and define
\[
 \|\Gamma\|_{H^s}:=\|\eta_\Gamma\|_{H^s(\Gamma_*)}.
\]
We use the following spaces. The estimates will first be proved in
$\mathbf H^s$ and $\mathbf H^s_1$, while the local well-posedness result
is stated in $\mathbf H^s_{1/2}$.

\begin{definition}\label{def:Hs}
For a state $(h,u,\Omega)$ with $\Gamma=\partial\Omega\in\Lambda_*$, define
\[
 (h,u,\Omega)\in\mathbf H^s
 \quad\Longleftrightarrow\quad
 \Gamma\in H^s,\qquad h,u\in H^s(\Omega),
\]
and
\[
 (h,u,\Omega)\in\mathbf H^s_1
 \quad\Longleftrightarrow\quad
 (h,u,\Omega)\in\mathbf H^s,\qquad
 (\divg u,\Delta h)\in H^s(\Omega)\times H^{s-1}(\Omega).
\]
If $s\notin\mathbb N$, define
\[
 (h,u,\Omega)\in\mathbf H^s_{1/2}
 \quad\Longleftrightarrow\quad
 (h,u,\Omega)\in\mathbf H^s,\qquad
 (h,\divg u)\in
 H^{s+\frac12}(\Omega)\times H^{s-\frac12}(\Omega).
\]
The corresponding norms are
\begin{equation}\label{eq:Hs0}
 \|(h,u,\Omega)\|_{\mathbf H^s}
 :=\|\Gamma\|_{H^s}+\|(h,u)\|_{H^s(\Omega)},
\end{equation}
\begin{equation}\label{eq:Hs1}
 \|(h,u,\Omega)\|_{\mathbf H^s_1}
 :=\|\Gamma\|_{H^s}+\|(h,u)\|_{H^s(\Omega)}
   +\|(\divg u,\Delta h)\|_{H^s(\Omega)\times H^{s-1}(\Omega)},
\end{equation}
and
\begin{equation}\label{eq:Hs1/2}
 \|(h,u,\Omega)\|_{\mathbf H^s_{1/2}}
 :=\|\Gamma\|_{H^s}+\|u\|_{H^s(\Omega)}
   +\|(h,\divg u)\|_{H^{s+\frac12}(\Omega)\times
      H^{s-\frac12}(\Omega)}.
\end{equation}
\end{definition}

Since $D_t$ is tangential to the moving boundary, differentiating
\eqref{eq:vacbc} gives $D_t^jh|_{\Gamma_t}=0$ whenever the trace is
defined. This imposes compatibility conditions on the initial data. Set
\begin{equation}\label{eq:compatcond0}
 z_j(h_0,u_0):=(-D_t)^jh\big|_{t=0}.
\end{equation}
By \eqref{eq:CE}, the time derivatives can be written in terms of spatial
derivatives of $(h_0,u_0)$. For a general state $(h,u)$, we write
\begin{equation}\label{eq:zjwj}
 z_j(h,u)=Z_j+N_j(h,u),
\end{equation}
where $Z_j$ is the linear part and $N_j$ is the nonlinear part. In
particular, $Z_0=h$, $Z_1=\divg u$, and $Z_2=\Delta h$. The compatibility
conditions are $z_j(h_0,u_0)|_{\Gamma_0}=0$.

When $s\in\mathbb N$, the highest quantity $z_s$ has only $H^{1/2}$
regularity, so its boundary trace need not be defined. We require instead
that $z_s\in H^{1/2}_{00}(\Omega)$, the Lions--Magenes space in which
extension by zero belongs to $H^{1/2}(\mathbb R^3)$. Accordingly, we make
the following definition for integer $s$.

\begin{definition}\label{def:ints}
If $s\in\mathbb N$, a state $(h,u,\Omega)$ belongs to
$\mathbf H^s_{1/2}$ if and only if
\[
 (h,u,\Omega)\in\mathbf H^s,\qquad
 (h,\divg u)\in
 H^{s+\frac12}(\Omega)\times H^{s-\frac12}(\Omega),\qquad
 z_s\in H^{1/2}_{00}(\Omega).
\]
In this case we set
\[
 \|(h,u,\Omega)\|_{\mathbf H^s_{1/2}}
 :=\|\Gamma\|_{H^s}+\|u\|_{H^s(\Omega)}
 +\|(h,\divg u)\|_{H^{s+\frac12}(\Omega)\times
 H^{s-\frac12}(\Omega)}
 +\|z_s\|_{H^{1/2}_{00}(\Omega)}.
\]
\end{definition}

Since the fluid domains vary with the initial data, we also need to specify
the meaning of convergence. We compare the surfaces as normal graphs over
$\Gamma_*$ and the fluid variables by extending them to $\mathbb R^3$.

\begin{definition}\label{def:cvgHs1/2}
A sequence $(h^m,u^m,\Omega^m)\in\mathbf H^s_{1/2}$ converges to
$(h,u,\Omega)$ if
\[
 \|\eta_{\Gamma^m}-\eta_\Gamma\|_{H^s(\Gamma_*)}\longrightarrow0
\]
and there exist extensions $(\widetilde h^m,\widetilde u^m)$ and
$(\widetilde h,\widetilde u)$ to $\mathbb R^3$ such that
\[
 \|\widetilde u^m-\widetilde u\|_{H^s(\mathbb R^3)}
 +\|(\widetilde h^m-\widetilde h,
       \divg\widetilde u^m-\divg\widetilde u)\|_{
       H^{s+\frac12}(\mathbb R^3)\times H^{s-\frac12}(\mathbb R^3)}
 \longrightarrow0.
\]
If $s\in\mathbb N$, we additionally require
\[
 \|z_s^m\mathbf1_{\Omega^m}-z_s\mathbf1_\Omega\|_
 {H^{1/2}(\mathbb R^3)}\longrightarrow0.
\]
\end{definition}

By the continuously domain-dependent extension operators in \cite{IE},
this definition is equivalent to the formulation using intermediate
functions on enlarged domains. Our main result is the following.

\begin{theorem}[Hadamard local well-posedness]\label{thm:LWP}
Let $s>3$.  Suppose that
$(h_0,u_0,\Omega_0)\in\mathbf H^s_{1/2}$ satisfies the Taylor sign condition
\eqref{eq:TS} and the compatibility conditions
\[
 z_j(h_0,u_0)|_{\Gamma_0}=0,\qquad j<s.
\]
Then there exists a time $T>0$, depending only on
$\|(h_0,u_0,\Omega_0)\|_{\mathbf H^s_{1/2}}$ and
$\inf_{\Gamma_0}a_0$, for which the free surface compressible Euler system has a unique
solution
\[
 (h,u,\Omega_t)\in C([0,T],\mathbf H^s_{1/2})
\]
with initial state $(h_0,u_0,\Omega_0)$.  Moreover, the data-to-solution map is
continuous in this topology: if a sequence of initial states converges to
$(h_0,u_0,\Omega_0)$ in $\mathbf H^s_{1/2}$, then the corresponding solutions,
on a common time interval, converge to $(h,u,\Omega_t)$ in
$C([0,T],\mathbf H^s_{1/2})$.
\end{theorem}

The reference surface $\Gamma_*$ can be any smooth boundary as above.
No irrotationality assumption is made, and the vorticity is controlled in
$H^{s-1}(\Omega_t)$.

The regularity requirements on the acoustic variables are connected to
the evolution of the Taylor coefficient and to the comparison of solutions
on different domains; these issues are discussed in
Section~\ref{sub:technical}. We do not claim that the threshold $s>3$ is
optimal. The incompressible results \cite{WZZZ,IE} reach the lower threshold
$H^{5/2+}$ for the velocity and the free surface in their respective
geometric settings. Reaching that threshold
for compressible liquids would require finer acoustic distance estimates and a
regularization procedure that preserves compatibility at lower regularity.

\paragraph{Comments on general pressure laws.}
The restriction to a linear pressure law serves only to simplify the
presentation. Extending the argument to smooth, nondegenerate liquid
pressure laws requires no new ideas or techniques, but only more involved
calculations with variable coefficients. More precisely,
for a sufficiently smooth pressure law $P=P(\rho)$ with $P'(\rho)>0$
and constant boundary density $\rho_b>0$, introduce the normalized enthalpy
and squared sound speed
\[
 h_P(\rho)=\int_{\rho_b}^{\rho}\frac{P'(r)}{r}\,dr,
 \qquad q(h_P)=P'\bigl(\rho(h_P)\bigr),
\]
where $\rho(h_P)$ denotes the inverse enthalpy map. Then the equations take
the form
\[
 D_th_P+q(h_P)\divg u=0,\qquad D_tu+\nabla h_P=0,
 \qquad h_P|_{\Gamma_t}=0;
\]
see \cite[Section~2]{L2}. If the density ranges in a compact subset of
$(0,\infty)$ and $0<q_*\leq q(h_P)\leq q^*<\infty$, the acoustic
operator remains nondegenerate up to the boundary. The coupling to the
free surface, in particular the need to control $D_ta$ with
$a=-\partial_{n_{\Gamma_t}}h_P>0$, thus motivates the same Sobolev scale
\[
 \Gamma\in H^s,\qquad u\in H^s(\Omega),\qquad
 h_P\in H^{s+\frac12}(\Omega),\qquad
 \divg u\in H^{s-\frac12}(\Omega),\qquad s>3.
\]
The Sobolev exponents are therefore unchanged, although the energy weights
and compatibility expressions depend on $P$. In particular, the linearized
bulk energy has integrand $q(h_P)^{-1}z^2+|w|^2$ for variations $(z,w)$ of
$(h_P,u)$, giving an equivalent $L^2$ norm under the stated bounds.
Similarly, the quantities $z_j^{(P)}=(-D_t)^jh_P$ are computed from the
general-pressure equations and satisfy $z_j^{(P)}|_\Gamma=0$ for $j<s$,
together with $z_s^{(P)}\in H^{1/2}_{00}(\Omega)$ when $s\in\mathbb N$,
as in Definition~\ref{def:ints}. Thus the precise admissible states still
depend on $P$. Derivatives of $q$ introduce additional terms in the energy,
commutator, and compatibility calculations, including those used in
compatibility-preserving regularization; these terms are treated by the
same estimates and constructions developed here. The additional work is
therefore computational and does not require new analytical ingredients.
To keep the presentation concise, we give the detailed proof of
Theorem~\ref{thm:LWP} only for the linear pressure law. This observation
concerns nondegenerate liquids and does not apply to a physical vacuum,
where vanishing density and sound speed require weighted bulk spaces.

\subsection{Related work}

The local existence theory for compressible liquids with a free surface
was developed by Lindblad \cite{L1,L2}, who proved estimates for the
linearized problem and obtained nonlinear local existence by a Nash--Moser
scheme. Trakhinin \cite{Trakhinin} obtained local existence in Eulerian
coordinates using a Nash--Moser scheme and extended the theory to full gas
dynamics and relativistic compressible Euler equations. Wang, Zhang, and
Zhao \cite{WZZ} proved local well-posedness and the incompressible limit
in Eulerian coordinates by relating the Taylor sign condition to the
hyperbolicity of an evolution equation for the free surface. These works
provide the local existence background for the liquid problem with positive
boundary density considered here. Coutand, Hole, and Shkoller \cite{CHS} proved
well-posedness with surface tension and justified the zero-surface-tension
limit. Lindblad and Luo \cite{LindLuo} established a priori energy
estimates for compressible liquids with a free surface, allowing nonzero
vorticity, and used estimates uniform in the sound speed to justify the
incompressible limit. Luo \cite{LuoGravity} obtained a priori estimates
and an incompressible limit for compressible gravity water waves with
vorticity in an unbounded domain. Luo and Zhang \cite{LZh} subsequently
proved local well-posedness for this problem using tangential smoothing.

Huang, Li, and Wang \cite{HLWCorner} proved local well-posedness for
two-dimensional compressible capillary--gravity water waves in bounded
corner domains with acute contact angles. Their geometric formulation
couples the acoustic evolution to the curvature dynamics, and combines
analysis of corner singularities with dissipation at the contact points
to obtain energy estimates and construct solutions. Their setting includes
surface tension and moving contact points; the present paper addresses
three-dimensional compressible liquids with a closed free surface and no
surface tension.

Geometric energy estimates also apply to related compressible free-boundary
models. Ginsberg, Lindblad, and Luo \cite{GLL} proved local well-posedness
for a compressible self-gravitating liquid. Luo, Trivisa, and Zeng
\cite{LTZ} derived  a priori estimates without loss of
derivatives for non-isentropic Euler flows and for Euler--Poisson flows
with a nonlinear electric potential, under the respective boundary sign
conditions. For highly subsonic heat-conducting inviscid flows, Luo and
Zeng \cite{LZHeat} obtained a priori estimates coupling the free-surface
geometry to the interior fluid variables. The additional entropy,
potential, and thermal variables in these models lead to different
couplings, while the control of boundary geometry remains relevant to
the liquid problem studied here.

For ideal compressible magnetohydrodynamics, Trakhinin and Wang
\cite{TW} established local existence and uniqueness for the
nonrelativistic and relativistic free-boundary problems with vanishing
total pressure on the interface. Their argument uses tame estimates in
anisotropic Sobolev spaces and Nash--Moser iteration under a
Rayleigh--Taylor sign condition on the total pressure. This is a related
characteristic free-boundary setting in which the magnetic field enters
both the interior dynamics and the boundary stability condition.

For compressible gases with a physical vacuum, the density and sound
speed vanish at the free boundary, and the local theory uses weighted
regularity adapted to this degeneracy. Coutand and Shkoller \cite{CSVac}
and Jang and Masmoudi \cite{JMVac} established local well-posedness for
three-dimensional polytropic Euler flows in this regime. Luo, Xin, and
Zeng \cite{LXZ} proved uniqueness for general three-dimensional physical
vacuum motions and a local existence theory for spherically symmetric
flows, with or without self-gravitation. Gu and Lei \cite{GLVac}
established local well-posedness for the three-dimensional
Euler--Poisson physical-vacuum problem. These results concern a
degenerate acoustic equation, whereas the liquid equation in the present
paper is nondegenerate up to the boundary.

The physical-vacuum literature also includes global evolution in special
regimes. Had\v{z}i\'c and Jang proved nonlinear stability of expanding
affine Euler motions \cite{HJExpand} and of radially symmetric expanding
stars for the mass-critical Euler--Poisson system \cite{HJStars}.
For Euler equations with damping, Luo and Zeng \cite{LZDamping}
proved global existence and convergence to Barenblatt profiles in one
dimension; Zeng obtained corresponding global results for
three-dimensional spherically symmetric motions \cite{ZengDamping}
and almost global existence near Barenblatt profiles without symmetry
assumptions \cite{ZengAlmost}. More recently, Zeng \cite{ZengGlobal}
proved global existence and asymptotic equivalence to Barenblatt-type
solutions for small perturbations in both two and three dimensions,
treating constant damping and damping $(1+t)^{-\lambda}$ with
$0<\lambda<1$.
For self-gravitating gases, Jang
\cite{JangStars} proved nonlinear instability of Lane--Emden equilibria,
and Guo, Had\v{z}i\'c, and Jang \cite{GHJ} constructed solutions
undergoing continued gravitational collapse. These results address
long-time dynamics in physical-vacuum regimes and are complementary to
the local liquid theory developed here.

Low-regularity local well-posedness for incompressible gravity water waves
with vorticity and without surface tension was established by Wang, Zhang,
Zhao, and Zheng \cite{WZZZ}. In a finite-depth domain with a flat bottom
and a graph free surface, their result gives velocity and free-surface
regularity $H^r$, $r>5/2$ in three dimensions, under the Taylor sign and
positive-depth conditions. Their paradifferential approach uses a good
unknown for the boundary velocity, estimates for the Dirichlet--Neumann
operator, and a modified energy argument that handles the apparent
half-derivative loss in the vorticity terms. They also prove a
Beale--Kato--Majda type breakdown criterion involving the mean curvature,
the Lipschitz norm of the velocity, the Taylor sign, and a positive lower
bound on the fluid depth. This provides a different approach to
low-regularity incompressible free-boundary flows; the compressible liquid
problem studied here additionally requires control of a nondegenerate
acoustic evolution and its nonlinear boundary compatibility conditions.

Recent Eulerian approaches provide a framework for proving Hadamard
well-posedness at low regularity. For incompressible free-boundary Euler,
Ifrim, Pineau, Tataru, and Taylor \cite{IE} established existence,
uniqueness, and strong continuous dependence for the velocity and the
free surface in $H^s$, with $s>5/2$ in three dimensions. Their analysis
combines geometric energy estimates, comparison of solutions on different
domains, and approximation by regular solutions. Related methods were
developed for the incompressible free-boundary MHD equations in
\cite{MHD}. For physical-vacuum gases, Ifrim and Tataru \cite{PV}
established a Hadamard theory in weighted Sobolev spaces adapted to the
vanishing sound speed. Liu and Luo \cite{LL} developed the corresponding
Eulerian theory for the full compressible Euler equations, including the
transport and coupling of nonconstant entropy. These works motivate our
use of regularized time discretization, distance functionals, and
frequency envelopes. Applying this framework to liquids requires a
simultaneous treatment of the free-surface geometry and a nondegenerate
acoustic boundary problem. In particular, the nonlinear compatibility
conditions must be incorporated both into approximation of the initial
state and into construction of the evolution. We discuss these issues below.

\subsection{Technical difficulties and main ideas}\label{sub:technical}

\paragraph{Coupling acoustic waves to the free surface.}
In the incompressible problem, the pressure is determined at each time by
the elliptic equation
\[
 -\Delta p=\operatorname{tr}\bigl((\nabla u)^2\bigr),
 \qquad p|_{\Gamma_t}=0,
\]
together with $\divg u=0$. Thus pressure is recovered from the velocity
and the domain, and is not an independent initial variable. The
free-surface dynamics and the transport of vorticity are coupled through
this elliptic problem; see \cite{IE}. For a compressible liquid, the
enthalpy and the divergence instead evolve as an acoustic system. The
term $D_t^2h$ in the wave equation displayed above prevents recovery of
$h$ from an elliptic equation involving only the instantaneous velocity.
At the same time, the acoustic normal derivative determines the Taylor
coefficient driving the boundary motion. Consequently, estimates for the
wave equation and estimates for the moving surface must close together.

This coupling is visible already in the linearized energy. If $(z,w)$
denotes the variation of $(h,u)$ and $\vartheta$ the normal displacement
of the boundary, then the linearized boundary condition is
$z=a\vartheta$. The natural energy is
\[
 \int_{\Omega_t}\bigl(z^2+|w|^2\bigr)\,dx
 +\int_{\Gamma_t}a\vartheta^2\,dS.
\]
Its derivative contains $D_ta$, whose leading term is
$\nabla_n\divg u$; see Proposition~\ref{prop:linEid}. Controlling only
the Lipschitz norm of $u$ does not control this term. In three dimensions,
the additional regularity $\divg u\in H^{s-1/2}$, $s>3$, supplies the
needed pointwise bound. The corresponding acoustic regularity is
$h\in H^{s+1/2}$. This explains the choice of the hybrid space
$\mathbf H^s_{1/2}$ and why a common Sobolev index for all variables
does not reflect the estimates used here.

At high order, we combine material derivatives of $(h,u)$ with
tangential derivatives, boundary curvature, and vorticity. The material
derivatives capture the acoustic evolution, while the tangential energy
cancels the boundary terms through the evolution of the curvature.
Div--curl and elliptic estimates then recover full spatial regularity.
We also propagate one additional acoustic derivative in
$\mathbf H^k_1$, without requiring that additional derivative of the
entire velocity or the surface. These two energy levels are needed to
recover the half-derivative separation in the final solution space.

\paragraph{The distinction from a physical vacuum.}
In the physical-vacuum setting, the squared sound speed vanishes linearly
with the distance to the boundary. The acoustic operator therefore
degenerates there, and the analysis uses weighted Sobolev spaces and
degenerate elliptic operators; see \cite{PV,LL}. For example, for the
polytropic law $P(\rho)=\rho^\gamma$, $\gamma>1$, the variable
$r=P'(\rho)=\gamma\rho^{\gamma-1}$ satisfies
\[
 D_tr+(\gamma-1)r\divg u=0.
\]
Since $r=0$ at the gas boundary, this equation propagates its zero set
without requiring $\divg u=0$ there. The vanishing acoustic speed and the
boundary motion are thus built into the same degenerate evolution.

For the liquid considered here, the boundary density and sound speed
remain positive. The acoustic equation is nondegenerate and must satisfy
the Dirichlet condition $h|_{\Gamma_t}=0$ on the moving surface. In
particular, $D_th=-\divg u$ forces $\divg u|_{\Gamma_t}=0$ for regular
solutions. Further material differentiation gives additional constraints.
The absence of acoustic degeneracy therefore brings a boundary
compatibility problem that is not present in the same form in the
physical-vacuum formulation. The liquid estimates must control the
resulting boundary traces in ordinary Sobolev spaces.

\paragraph{Regularization within the compatibility constraints.}
The first conditions in the hierarchy are
\[
 h|_\Gamma=0,\qquad \divg u|_\Gamma=0,\qquad
 \left.\left(\Delta h+
 \operatorname{tr}\bigl((\nabla u)^2\bigr)\right)\right|_\Gamma=0.
\]
Higher conditions involve increasingly many derivatives and nonlinear
products. Convolution does not preserve these identities, and smoothing
the domain also changes the surface on which they are evaluated.
Moreover, the rough initial state satisfies only the conditions meaningful
at its regularity, whereas the smoother approximations need additional
conditions to generate regular solutions. Thus one must both preserve
the existing hierarchy and extend it to higher order, with errors small
enough for the convergence argument.

The operators $\Psi_m$ in Proposition~\ref{prop:Psim} accomplish this by
combining smoothing of the domain and the fluid variables with elliptic
boundary corrections. The decomposition $z_j=Z_j+N_j(h,u)$ separates
the principal linear boundary operators from their nonlinear errors and
allows an iterative correction. The resulting states satisfy any
prescribed higher finite order of compatibility, while retaining the
regularity and approximation bounds at scale $2^{-m}$. These bounds
must also distinguish $(\nabla h,\divg u)$ from the remaining variables.
At integer $s$, the last compatibility expression belongs only to
$H^{1/2}$ and need not have a classical trace. The Lions--Magenes condition on
$z_s$ supplies the endpoint boundary control needed for the same
regularization argument and for convergence of its zero extension.

\paragraph{Compatibility during the discrete evolution.}
The regularized forward Euler construction introduces a second use of
boundary correction. A direct time step loses a spatial derivative and
does not preserve the higher compatibility conditions. Smoothing the
input can restore regularity, but it also perturbs the boundary traces
and the energy. It is therefore insufficient to produce a smooth state
with a uniformly bounded energy at each step: over a fixed time interval
there are of order $\epsilon^{-1}$ steps, so the energy increment must be
of order $\epsilon$.

Our construction regularizes the domain, the acoustic component, and the
tangential velocity in separate stages using elliptic resolvents. Boundary
lifting operators with prescribed normal derivatives then reduce the
compatibility errors before transport. The dissipation supplied by the
resolvents absorbs the highest-order errors in the energy calculation.
The one-step estimate in Theorem~\ref{thm:1step} propagates quantitative
approximate compatibility, rather than assuming exact compatibility at
every intermediate state. Together with the $O(\epsilon)$ energy
increment and the $O(\epsilon^2)$ consistency error, this provides bounds
uniform in the time step and permits passage to a regular solution.

\paragraph{Distances on different domains and strong convergence.}
To establish uniqueness and continuous dependence, one must compare
solutions whose domains move differently. We use an Eulerian distance
on their intersection, following the approach in \cite{IE,PV,LL}.
For liquids, the linearized boundary energy identifies the necessary
surface term: the square of the enthalpy difference is weighted by the
inverse Taylor coefficient. Near either free surface, the Taylor sign
relates that difference to the normal separation of the boundaries.
This allows the $L^2$ distance to control both the bulk variables and
the interface. Its propagation yields uniqueness under the pointwise
bounds in Theorem~\ref{thm:L2dist}.

The $L^2$ distance alone does not provide the acoustic approximation
rate required by the hybrid topology. We therefore introduce a partial
$H^2$ distance based on the second material derivatives of $(h,u)$.
It controls the differences of $\Delta h$ and $\nabla\divg u$, modulo
lower-order terms. There is a further geometric difficulty: although
each enthalpy vanishes on its own boundary, their difference does not
satisfy a homogeneous Dirichlet condition on the boundary of the
intersection. Thus an estimate for its Laplacian cannot be converted
directly into an estimate for its full Hessian. We align the domains
by a diffeomorphism, estimate the errors introduced in the elliptic
operator, and then transfer the resulting acoustic bounds back to the
common region; see Proposition~\ref{prop:H2dist} and
Corollary~\ref{cor:PH2}.

Finally, the two distance estimates are paired with the energy bounds
in $\mathbf H^k$ and $\mathbf H^k_1$ through frequency envelopes.
Compactness by itself only gives convergence in weaker spaces; strong
convergence in $\mathbf H^s_{1/2}$ also requires control of the
high-frequency tails. The envelope estimates provide this control for
the surface, velocity, and acoustic variables at their respective
regularities. Uniform control of the tails for convergent initial states,
together with continuity at a fixed regularization scale, yields
continuity of the data-to-solution map in the full hybrid topology.
The endpoint compatibility quantity is included in this convergence
when $s$ is an integer.

\subsection{Outline of the proof}

Section~\ref{sec:prel} derives the linearized energy identity, constructs
the compatibility-preserving regularization operators $\Psi_m$, and
proves their frequency-envelope bounds. Section~\ref{sec:dist} develops
the $L^2$ and partial $H^2$ distances for solutions on different domains.
These estimates provide the stability bounds used in the approximation
argument.

Section~\ref{sec:apest} proves coercivity and propagation of the energies
in $\mathbf H^k$ and $\mathbf H^k_1$. Section~\ref{sec:regsol}
constructs regular solutions by the discrete scheme, including the
regularization, compatibility correction, transport, and convergence
steps. In Section~\ref{sec:rough}, the energy and distance estimates are
combined with frequency envelopes to establish a common lifespan for
regularized solutions, strong convergence in $C([0,T],\mathbf H^s_{1/2})$,
and continuous dependence. This completes the proof of
Theorem~\ref{thm:LWP}.

We write $E\lesssim_A F$ for $E\leq C(A)F$. Dependence on the fixed
collar $\Lambda_*$ and on the lower bound of the Taylor coefficient is
usually suppressed. For a spacetime operator $T$ and a function $G$,
$T\,f(G)$ means that composition is performed before applying $T$, while
$Tf\circ G$ means that $T$ is applied before composition.

\section{Preliminaries}\label{sec:prel}

\subsection{Linearization and basic energy identity}\label{sec:LinEq}

We linearize \eqref{eq:CE} together with the free-boundary conditions. Let $(h^\epsilon,u^\epsilon,\Omega^\epsilon)$ be a one-parameter family of solutions with $(h^\epsilon,u^\epsilon,\Omega^\epsilon)|_{\epsilon=0}=(h,u,\Omega)$. Write $(z,w)=\partial_\epsilon|_{\epsilon=0}(h^\epsilon,u^\epsilon)$ and let $s=(\partial_\epsilon|_{\epsilon=0}\psi^\epsilon)\cdot n_\Gamma$ denote the normal variation of the interface. Here $\psi^\epsilon$ parameterizes the varying boundary. The linearized variables satisfy
\begin{equation}\label{eq:LinCE}
\left\{\begin{aligned}
&D_tz+\divg w=-w\cdot\nabla h\quad&\text{in }\Omega,\\
&D_tw+\nabla z=-w\cdot\nabla u\quad&\text{in }\Omega,\\
&D_ts-w\cdot n=s(\nabla_nu\cdot n)\quad&\text{on }\Gamma,\\
&z-as=0\quad&\text{on }\Gamma.
\end{aligned}\right.
\end{equation} 
The geometric linearization is obtained as in \cite{IE}. When $h,u$ are Lipschitz and $\Gamma$ is $C^1$, the right-hand sides are lower-order terms. We therefore consider the following inhomogeneous system:
\begin{equation}\label{eq:LinCE2}
\left\{\begin{aligned}
&D_tz+\divg w=f\quad&\text{in }\Omega,\\
&D_tw+\nabla z=g\quad&\text{in }\Omega,\\
&D_ts-w\cdot n=q\quad&\text{on }\Gamma,\\
&z-as=r\quad&\text{on }\Gamma.
\end{aligned}\right.
\end{equation} 

\begin{proposition}\label{prop:linEid}
Suppose that $(z,w,s)$ is a solution to \eqref{eq:LinCE2}. Define the linearized energy as
\[
E_\text{lin}(z,w,s)=\int_\Omega z^2+|w|^2\,dx+\int_\Gamma as^2\,dS. 
\]
Then
\[\begin{aligned}
\frac{1}{2}\frac{d}{dt}E_\text{lin}(z,w,s)&=\int_\Omega fz+g\cdot w\,dx+\int_\Gamma qs\, adS-\int_\Gamma rw\cdot n\,dS\\
&\quad+\int_\Omega\frac{1}{2}(z^2+|w|^2)\divg u\,dx+\int_\Gamma\frac{s^2}{2}(D_ta+a\divg_\Gamma u)\,dS. 
\end{aligned}\]
Consequently, if $\Vert\nabla u\Vert_{L^\infty(\Omega)}+\Vert D_ta\Vert_{L^\infty(\Gamma)}\leq B$, $\Vert(f,g)\Vert_{L^2(\Omega)}+\Vert q\Vert_{L^2(\Gamma)}\leq M$, and the Taylor condition holds on $[0,T]$, then we have
\[
\frac{d}{dt}E_\text{lin}(z,w,s)\lesssim BE_\text{lin}(z,w,s)+ME_\text{lin}(z,w,s)^\frac{1}{2}+\Vert r\Vert_{H^\frac{1}{2}(\Gamma)}\Vert w\cdot n\Vert_{H^{-\frac{1}{2}}(\Gamma)}. 
\]
\end{proposition}
\par\smallskip\noindent\textbf{Remark.}\ The inhomogeneous system \eqref{eq:LinCE2} also accommodates commuted equations with $s\equiv0$. In particular, the material derivatives $D_t^jh$ vanish on the moving boundary whenever their traces are defined. The corresponding commuted equations generally contain forcing terms.
\proofheading{Proof}Multiply the bulk equations in \eqref{eq:LinCE2} by $z,w$ and the boundary evolution equation by $as$. Integrating, using the boundary relation $z=as+r$, and applying the divergence theorem and Lemma~\ref{lem:tpid} gives the identity. The estimate follows by Cauchy--Schwarz and trace duality. \proofqed

\subsection{Compatibility-preserving regularization}\label{sub:regop}

We construct smooth approximations that retain quantitative control of the state and satisfy higher-order compatibility conditions.

\begin{proposition}\label{prop:Psim}
Let $s\geq2$, $s'=s+\sigma$ with $\sigma\in\{0,\frac{1}{2}\}$ and let $K>s'-\frac{1}{2}$ be odd integer. Suppose
\[
A\geq\Vert\Gamma\Vert_{C^{1,\delta}}+\Vert(h,u)\Vert_{C^{\frac{1}{2}+\delta}(\Omega)}+\mathbf 1_{\sigma=\frac{1}{2}}\Vert(\divg u,\nabla h)\Vert_{C^\delta(\Omega)}
\]
and
\[
B\geq\Vert\Gamma\Vert_{H^s}+\Vert(h,u)\Vert_{H^s(\Omega)}+\Vert(\divg u,\nabla h)\Vert_{H^{s'-1}(\Omega)}. 
\]
Suppose 
\[
z_j|_\Gamma=0,\quad 0\leq j<s'-\frac{1}{2},\quad\Vert\mathbf1_{s'-\frac{1}{2}\in\mathbb N}\ z_{s'-\frac{1}{2}}\Vert_{H_{00}^{1/2}(\Omega)}\leq B. 
\]
Then for $m\gg1$ there is a $(h_m,u_m,\Omega_m)=\Psi_m(h,u,\Omega)$ with the following properties: \\
(P1) Regularity bound: 
\[
\Vert\Gamma_m\Vert_{C^{1,\delta}}+\Vert(h_m,u_m)\Vert_{C^{\frac{1}{2}+\delta}(\Omega_m)}+\mathbf 1_{\sigma=\frac{1}{2}}\Vert(\divg u_m,\nabla h_m)\Vert_{C^\delta(\Omega_m)}\lesssim_A1,
\]
and
\[
\Vert\Gamma_m\Vert_{H^{s+r}}+\Vert(h_m,u_m)\Vert_{H^{s+r}(\Omega_m)}+\Vert(\divg u_m,\nabla h_m)\Vert_{H^{s'-1+r}(\Omega_m)}\lesssim_A2^{mr}B.
\]
(P2) Error bound: 
\[
\Vert\eta-\eta_m\Vert_{L^2(\Gamma_*)}+\Vert(h-h_m,u-u_m)\Vert_{L^2(\Omega\cap\Omega_m)}\lesssim_A2^{-ms}B,
\]
and
\[
\Vert(\divg u-\divg u_m,\nabla h-\nabla h_m)\Vert_{L^2(\Omega\cap\Omega_m)}\lesssim_A2^{-m(s'-1)}B. 
\]
(P3) Compatibility: 
\[
z_{j,m}|_{\Gamma_m}=0,\quad j=0,1,...,K.
\]
\end{proposition}
We begin with the domain-enlarging regularization from \cite{IE}. Choose $\chi_0\in C_c^\infty(\Omega_*)$, $\chi_i\in C_c^\infty(\mathbb R^3)$, and directions $e_i\in S^2$ such that
\[
\sum_{i=0}^L\chi_i=1\text{ on }\Omega_*,\quad\sum_{i=1}^L\chi_i=1\text{ near }\Gamma_*,\quad|e_i-n(x)|<\frac{1}{10},\ x\in\supp\chi_i\cap\Gamma_*. 
\]
Next select 
\[
\phi_0\in C_c^\infty(\mathbb R^3),\quad\phi_i\in C_c^\infty\left(B_\frac{1}{10}(e_i)\right),\quad i=1,...,L
\]
such that
\[
\int\phi_i(x)x^\alpha\,dx=\begin{cases}1&|\alpha|=0\\0&|\alpha|=1,...,N\end{cases}
\]
and define $\phi_i^m(x)=2^{3m}\phi_i(2^mx)$. Here $N$ is a sufficiently large integer. The condition on moments of $\phi_i$ implies that the Fourier transform of $\phi_i$ is $1+O(|\xi|^{N+1})$ as $\xi\to0$. Then define
\[
\Phi_mv=\sum_{i=0}^L\phi_i^m*(\chi_iv). 
\]

\begin{lemma}\label{lem:Phim}
Let $m\gg1$ and $\Omega_m$ be a domain bounded by $\Gamma_m\in\Lambda_*$ with $d_H(\Omega,\Omega_m)\leq c2^{-m}$. Then we have: (i)
\[
\Vert\Phi_mf\Vert_{C^{\alpha+r}(\Omega_m)}\lesssim_A2^{mr}\Vert f\Vert_{C^\alpha(\Omega)},\quad r\geq0,
\]
\[
\Vert\Phi_mf\Vert_{H^{s+r}(\Omega_m)}\lesssim_A2^{mr}\Vert f\Vert_{H^s(\Omega)},\quad r\geq0,
\]
\[
\Vert(I-\Phi_m)f\Vert_{H^{s-r}(\Omega)}\lesssim_A2^{-mr}\Vert f\Vert_{H^s(\Omega)},\quad r\in[0,s].
\]
(ii) Define $\Vert v\Vert_{C^\alpha_\text{div}(\Omega)}=\Vert v\Vert_{C^\alpha(\Omega)}+\Vert\divg v\Vert_{C^\alpha(\Omega)}$ and $\Vert v\Vert_{H^s_\text{div}(\Omega)}=\Vert v\Vert_{H^s(\Omega)}+\Vert\divg v\Vert_{H^s(\Omega)}$. Then
\[
\Vert\Phi_mv\Vert_{C^{\alpha+r}_\text{div}(\Omega_m)}\lesssim_A2^{mr}\Vert v\Vert_{C^\alpha_\text{div}(\Omega)},\quad r\geq0,
\]
\[
\Vert\Phi_mv\Vert_{H^{s+r}_\text{div}(\Omega_m)}\lesssim_A2^{mr}\Vert v\Vert_{H^s_\text{div}(\Omega)},\quad r\geq0,
\]
\[
\Vert(I-\Phi_m)v\Vert_{H^{s-r}_\text{div}(\Omega)}\lesssim_A2^{-mr}\Vert v\Vert_{H^s_\text{div}(\Omega)},\quad r\in[0,s].
\]
\end{lemma}
\proofheading{Proof}(i) can be found in Proposition A.32 in \cite{MHD} and in the proof of Proposition 6.3 in \cite{IE}. \\
(ii) It suffices to establish the following commutator estimates: 
\begin{equation}\label{eq:commpt}
\Vert[\partial_k,\Phi_m]v\Vert_{C^{\alpha+r}(\Omega_m)}\lesssim_A2^{mr}\Vert v\Vert_{C^\alpha(\Omega)},\quad r\geq0,
\end{equation}
\begin{equation}\label{eq:commhf}
\Vert[\partial_k,\Phi_m]v\Vert_{H^{s+r}(\Omega_m)}\lesssim_A2^{mr}\Vert v\Vert_{H^s(\Omega)},\quad r\geq0,
\end{equation}
\begin{equation}\label{eq:commlf}
\Vert[\partial_k,\Phi_m]v\Vert_{H^{s-r}(\Omega_m)}\lesssim_A2^{-mr}\Vert v\Vert_{H^s(\Omega)},\quad r\in[0,s]. 
\end{equation}
To this end, we use the identity 
\begin{align}
[\partial_k,\Phi_m]v&=\sum_i\phi_i^m*\left((\partial_k\chi_i)\mathcal E_\Omega v\right)&\text{ on }\Omega_m\cup\Omega\label{eq:hfid}\\
&=\sum_i(\phi_i^m*\cdot-I)\left((\partial_k\chi_i)\mathcal E_\Omega v\right)&\text{ on }\Omega_m\cup\Omega\label{eq:lfid}
\end{align}
where in the second line $\phi_i^m*\cdot$ denotes the convolution operator and we have used that $\sum_i\partial_k\chi_i=0$ near $\overline\Omega$. By \eqref{eq:hfid} and the estimate of $\phi_i^m*\cdot$ we obtain \eqref{eq:commpt}--\eqref{eq:commhf}. By \eqref{eq:lfid} and the estimate of $\phi_i^m*\cdot-I$ we obtain \eqref{eq:commlf}. \proofqed
The proof also uses the following composition, product, trace, and elliptic estimates. These are weaker forms of the estimates in Appendix~\ref{app:estimates}, with constants suited to the regularized domains.

\begin{proposition}\label{prop:Moser3}
Let $r\geq0$ and $G:\mathbb R^3\to\mathbb R^3$ be a diffeomorphism with $\Vert\nabla G\Vert_{C^\delta(\mathbb R^3)}+\Vert\nabla G^{-1}\Vert_{C^\delta(\mathbb R^3)}\leq A$. Then we have
\[
\Vert f(G)\Vert_{H^s(\mathbb R^3)}\lesssim_A\Vert f\Vert_{H^s(\mathbb R^3)}+\Vert G\Vert_{C^{s+r+\delta}(\mathbb R^3)}\Vert f\Vert_{H^{1-r}(\mathbb R^3)}. 
\]
\end{proposition}
\proofheading{Proof}This is Proposition 5.8 in \cite{IE}. \proofqed

\begin{lemma}\label{lem:biltrdiri2}
(i) For $s\geq0$ we have
\begin{equation}\label{eq:bil2}
\Vert fg\Vert_{H^s(\Omega)}\lesssim_A\Vert f\Vert_{L^\infty(\Omega)}\Vert g\Vert_{H^s(\Omega)}+\Vert f\Vert_{C^{s+\delta}(\Omega)}\Vert g\Vert_{L^2(\Omega)}. 
\end{equation}
(ii) For $s>0$ we have
\begin{equation}\label{eq:tr2}
\Vert f|_\Gamma\Vert_{H^s(\Gamma)}\lesssim_A\Vert f\Vert_{H^{s+\frac{1}{2}}(\Omega)}+\Vert\Gamma\Vert_{C^{s+\frac{1}{2}+\delta}}\Vert f\Vert_{H^1(\Omega)}. 
\end{equation}
(iii) For $s\geq2$ we have
\begin{equation}\label{eq:diri2}
\Vert f\Vert_{H^s(\Omega)}\lesssim_A\Vert\Delta f\Vert_{H^{s-2}(\Omega)}+\Vert f|_\Gamma\Vert_{H^{s-\frac{1}{2}}(\Gamma)}+\Vert\Gamma\Vert_{C^{s+\delta}}\Vert f\Vert_{H^1(\Omega)}. 
\end{equation}
\end{lemma}
\proofheading{Proof}For (i), use a paraproduct decomposition and estimate $P_jf$ in $L^\infty$ and $P_{<j}g$ in $L^2$. For (ii), combine the argument of Proposition~\ref{prop:tr} with the composition estimate in Proposition~\ref{prop:Moser3}. Assertion (iii) is Proposition 5.20 of \cite{IE}. \proofqed
\proofheading{Proof of Proposition~\ref{prop:Psim}}Define $\Omega_m$ by
\[
\eta_m=P_{<m}\eta. 
\] 
Set $(\tilde h,\tilde u)=(\Phi_mh,\Phi_mu)$. By Lemma~\ref{lem:Phim} $(\tilde h,\tilde u,\Omega_m)$ satisfies (P1)-(P2) and $(\tilde h,\tilde u)$ is defined on a neighborhood $\Omega^{[m]}$ enlarged by $2^{-m}$. Using that
\begin{equation}\label{eq:trint}
\Vert f|_\Gamma\Vert_{L^2(\Gamma)}\lesssim_A\Vert f\Vert_{H^{r_0}(\Omega)}^{1-\theta}\Vert f\Vert_{H^{r_1}(\Omega)}^\theta,\quad0\leq r_0<\frac{1}{2}<r_1\leq1,\ \theta=\frac{\frac{1}{2}-r_0}{r_1-r_0},
\end{equation}
using \eqref{eq:bil2} and using the bound on $h-\tilde h,u-\tilde u,\divg u-\divg\tilde u$ we derive
\[
\Vert\tilde z_j|_\Gamma\Vert_{L^2(\Gamma)}\lesssim_A2^{-m(s'-j-\frac{1}{2})}B,\quad j<s'-\frac{1}{2},
\]
and
\[
\Vert\tilde h|_\Gamma\Vert_{L^\infty(\Gamma)}+\mathbf1_{\sigma=\frac{1}{2}}2^{-m}\Vert\divg\tilde u|_\Gamma\Vert_{L^\infty(\Gamma)}\lesssim_A2^{-m(\frac{1}{2}+\sigma+\delta)}.
\]
If $s'-\frac{1}{2}=j_c\in\mathbb N$, we can also bound $\Vert\tilde z_{j_c}|_\Gamma\Vert_{L^2(\Gamma)}$. To do this, notice that the zero extension $z_{j_c}\mathbf 1_\Omega$ belongs to $H^\frac{1}{2}(\mathbb R^3)$ with $\Vert z_{j_c}\mathbf 1_\Omega\Vert_{H^{1/2}(\mathbb R^3)}\lesssim\Vert z_{j_c}\Vert_{H_{00}^{1/2}(\Omega)}$. Thus, we have
\[
2^\frac{m}{2}\Vert\Phi_m(z_{j_c}\mathbf 1_\Omega)|_{\mathbb R^3\backslash\Omega}\Vert_{L^2(\mathbb R^3\backslash\Omega)}+2^{-\frac{m}{2}}\Vert\Phi_m(z_{j_c}\mathbf 1_\Omega)\Vert_{H^1(\mathbb R^3)}\lesssim\Vert z_{j_c}\Vert_{H^{1/2}_{00}(\Omega)}\lesssim_AB,
\]
and by \eqref{eq:trint} on the complement domain we have $\Vert\Phi_mz_{j_c}|_\Gamma\Vert_{L^2(\Gamma)}\lesssim_AB$. We still need to control the difference $\Phi_mz_{j_c}-\tilde z_{j_c}$. Recall \eqref{eq:zjwj}. For the linear part we invoke \eqref{eq:commlf} and trace estimate to obtain $\Vert\Phi_mZ_{j_c}-\tilde Z_{j_c}\Vert_{H^{1/2}(\Gamma)}\lesssim_AB$. For the nonlinear part, from $s\geq2$ we note that bilinear terms in $N_{j_c}$ must contain $h,\divg u$ as a factors. Then by tame estimate we have $\Vert N_{j_c}\Vert_{H^1(\Omega)}+\Vert\tilde N_{j_c}\Vert_{H^1(\Omega)}\lesssim_AB$ and hence $\Vert\Phi_mN_{j_c}-\tilde N_{j_c}\Vert_{H^{1/2}(\Gamma)}\lesssim_AB$. Combining these we get
\[
\Vert\tilde z_{j_c}|_\Gamma\Vert_{L^2(\Gamma)}\lesssim_AB. 
\]
With the gradient estimate
\[
\Vert\nabla\tilde z_j\Vert_{L^\infty(\Omega^{[m]})}\lesssim_A2^{m(j+\frac{1}{2}-\sigma)}
\]
on some neighborhood $\Omega^{[m]}$ enlarged by roughly $2^{-m},$ the distance estimate 
\[
\Vert\eta-\eta_m\Vert_{L^2(\Gamma_*)}\lesssim_A2^{-ms}B,\quad\Vert\eta-\eta_m\Vert_{L^\infty(\Gamma_*)}\lesssim_A2^{-m(1+\delta)}
\]
and fundamental theorem of calculus, we further derive
\[%\begin{equation}\label{eq:trL2}
\Vert\tilde z_j|_{\Gamma_m}\Vert_{L^2(\Gamma_m)}\lesssim_A2^{-m(s'-j-\frac{1}{2})}B 
\]%\end{equation}
for $j\leq s'-\frac{1}{2}$, and derive
\[
\Vert\tilde h|_{\Gamma_m}\Vert_{L^\infty(\Gamma_m)}+\mathbf1_{\sigma=\frac{1}{2}}2^{-m}\Vert\divg\tilde u|_{\Gamma_m}\Vert_{L^\infty(\Gamma_m)}\lesssim_A2^{-m(\frac{1}{2}+\sigma+\delta)}.
\]
In fact, the first estimate is true for $s'-\frac{1}{2}<j\leq K$ as well; the second can be extended to
\[
\Vert\tilde z_j|_\Gamma\Vert_{L^\infty(\Gamma)}\lesssim_A2^{m(j-\frac{1}{2}-\sigma-\delta)},\quad j\leq K
\]
thanks to the regularity bound and \eqref{eq:trint}. Likewise, by \eqref{eq:tr2} we have
\[
\Vert\tilde z_j|_{\Gamma_m}\Vert_{H^{K-j+\frac{1}{2}}(\Gamma_m)}\lesssim_A2^{m(K+1-s')}B,\quad j\leq K
\]
and
\[
\Vert\tilde z_j|_{\Gamma_m}\Vert_{C^{K-j+\delta}(\Gamma_m)}\lesssim_A2^{m(K-\frac{1}{2}-\sigma)},\quad j\leq K. 
\]

We now correct the compatibility errors using Lemma~\ref{lem:calG}. The estimates above give
\[
P(\tilde h,\tilde u)\lesssim_A1,\quad Q(\tilde h,\tilde u)\lesssim_AB. 
\]
Starting with $(h^0,u^0)=(\tilde h,\tilde u)$, define $(h^{j+1},u^{j+1})=\mathcal G(h^j,u^j)$. For sufficiently large $m$, Lemma~\ref{lem:calG} makes $P(h^j,u^j)$ and $Q(h^j,u^j)$ decay geometrically. The correction bounds then give geometric convergence of $(h^j,u^j)$ in $H^{K+1}\cap C^{K+\delta}$.
%Moreover, using Lemma~\ref{lem:ImJm} (i), Proposition~\ref{prop:Schauder} and interpolation, we are able to show $\Vert(\divg u^j,\nabla h^j)\Vert_{C^\delta(\Omega_m)}\lesssim_A1$ when $\sigma=\frac{1}{2}$. 
Hence, the estimates (P1)-(P2) remain true for $(h^j,u^j)$ uniformly as $j\to\infty$. 
By setting $(h_m,u_m)=\lim(h^j,u^j)$ we conclude Proposition~\ref{prop:Psim}. \proofqed

\begin{lemma}\label{lem:ImJm}
With $\Omega_m$ above define $J_m=(\Delta_{\Omega_m}-2^{2m})^{-1}$ and define $v=I_m\psi$ via
\[\left\{\begin{aligned}
&v-2^{-2m}\Delta v=0\quad&\text{in }\Omega_m,\\
&v=\psi\quad&\text{on }\Gamma_m.
\end{aligned}\right.\]
Then we have: (i) 
\[
\Vert I_m\psi\Vert_{L^\infty(\Omega_m)}\leq\Vert\psi\Vert_{L^\infty(\Gamma_m)},\quad\Vert J_mf\Vert_{L^\infty(\Omega_m)}\leq2^{-2m}\Vert f\Vert_{L^\infty(\Omega_m)}. 
\]
(ii) 
\[
\Vert I_m\psi\Vert_{L^2(\Omega_m)}+2^{-m}\Vert I_m\psi\Vert_{H^1(\Omega_m)}\lesssim_A2^{-\frac{m}{2}}\Vert\psi\Vert_{L^2(\Gamma_m)}+2^{-m}\Vert\psi\Vert_{H^\frac{1}{2}(\Gamma_m)},
\]
and
\[
\Vert J_mf\Vert_{L^2(\Omega_m)}+2^{-m}\Vert J_mf\Vert_{H^1(\Omega_m)}+2^{-2m}\Vert J_mf\Vert_{H^2(\Omega_m)}\lesssim_A2^{-2m}\Vert f\Vert_{L^2(\Omega_m)}. 
\]
\end{lemma}
\par\smallskip\noindent\textbf{Remark.}\ The identity $I_m\psi=\mathcal H_m\psi-(I-2^{-2m}\Delta)^{-1}\mathcal H_m\psi$ will be useful below.
\proofheading{Proof}(i) is just the maximum principle. \\
(ii) The $L^2$ estimate of $J_m$ follows from the elliptic energy estimate. Then by \eqref{eq:diri2}, by
\begin{equation}\label{eq:Gammpt}
\Vert\Gamma\Vert_{C^{1+\delta+r}}\lesssim_A2^{mr}
\end{equation}
and by interpolation we obtain the $H^2$ estimate. Next, write $v=v_\text{ho}+v_\text{in}$, where $v_\text{ho}=\Delta_{\Omega_m}^{-1}\Delta v$ and $v_\text{in}=\mathcal H_m\psi$. The elliptic energy identity for $v$ gives
\[
\Vert v\Vert_{L^2(\Omega_m)}^2+2^{-2m}\Vert\nabla v_\text{ho}\Vert_{L^2(\Omega_m)}^2=2^{-2m}\int_{\Gamma_m}\nabla_nv_\text{ho}\,\psi\,dS. 
\]
The trace estimate gives
\[
\Vert\nabla v_\text{ho}\Vert_{L^2(\Gamma_m)}\lesssim_A\Vert\nabla v_\text{ho}\Vert_{L^2(\Omega_m)}^\frac{1}{2}\Vert\nabla v_\text{ho}\Vert_{H^1(\Omega_m)}^\frac{1}{2}. 
\]
By \eqref{eq:diri2}, \eqref{eq:Gammpt} we have
\[
\Vert\nabla v_\text{ho}\Vert_{H^1(\Omega_m)}\lesssim_A2^{2m}\Vert v\Vert_{L^2(\Omega_m)}+2^{-m}\Vert\nabla v_\text{ho}\Vert_{L^2(\Omega_m)}. 
\]
Now using Young's inequality we obtain 
\[
\Vert v\Vert_{L^2(\Omega_m)}+2^{-m}\Vert\nabla v_\text{ho}\Vert_{L^2(\Omega_m)}\lesssim_A2^{-\frac{m}{2}}\Vert\psi\Vert_{L^2(\Gamma_m)}+2^{-m}\Vert\psi\Vert_{H^\frac{1}{2}(\Gamma_m)}. 
\]
Together with $\Vert v_\text{in}\Vert_{H^1(\Omega_m)}\lesssim_A\Vert\psi\Vert_{H^{1/2}(\Gamma_m)}$, this proves the estimate for $I_m$. \proofqed

\begin{corollary}\label{cor:ImJm}
(i) Suppose $k$ is an even integer and
\[
\Vert\psi\Vert_{C^{k,\delta}(\Gamma_m)}+2^{m(k+\delta)}\Vert\psi\Vert_{L^\infty(\Gamma_m)}\leq P_I,
\]
\[
\Vert f\Vert_{C^{k,\delta}(\Omega_m)}+2^{m(k+\delta)}\Vert f\Vert_{L^\infty(\Omega_m)}\leq P_J.
\]
Then we have
\[
\Vert I_m\psi\Vert_{C^\alpha(\Omega_m)}\lesssim_A2^{-m(k+\delta-\alpha)}P_I,\quad\alpha\in[0,k+\delta],
\]
\[
\Vert J_mf\Vert_{C^\alpha(\Omega_m)}\lesssim_A2^{-m(k+2+\delta-\alpha)}P_J,\quad\alpha\in[0,k+2+\delta]. 
\]
(ii) Suppose $s\geq0$ and
\[
\Vert\psi\Vert_{H^{s+\frac{1}{2}}(\Gamma_m)}+2^{m(s+\frac{1}{2})}\Vert\psi\Vert_{L^2(\Gamma_m)}\leq Q_I,
\]
\[
\Vert f\Vert_{H^s(\Omega_m)}+2^{ms}\Vert f\Vert_{L^2(\Omega_m)}\leq Q_J.
\]
Then we have
\[
\Vert I_m\psi\Vert_{H^r(\Omega_m)}\lesssim_A2^{-m(s+1-r)}Q_I,\quad r\in[0,s+1],
\]
\[
\Vert J_mf\Vert_{H^r(\Omega_m)}\lesssim_A2^{-m(s+2-r)}Q_J,\quad r\in[0,s+2]. 
\]
\end{corollary}
\proofheading{Proof}(i) We mainly deal with the estimate of $I_m$ for $k=0$. Write $I_m\psi=w+\mathcal H_m\psi$, where $w$ is a solution to
\[
w-2^{-2m}\Delta w=-\mathcal H_m\psi\text{ in }\Omega_m,\quad w=0\text{ on }\Gamma_m. 
\]
By \eqref{eq:calHCalp} and the maximum principle we have $\Vert\mathcal H_m\psi\Vert_{C^\delta(\Omega_m)}+2^{m\delta}\Vert\mathcal H_m\psi\Vert_{L^\infty(\Omega_m)}\lesssim_AP_I$ and 
\[
\Vert w\Vert_{L^\infty(\Omega_m)}\leq2^{-\delta m}P_I. 
\]
By \eqref{eq:DiriCkalp} and \eqref{eq:Gammpt} we have
\[
\Vert w\Vert_{C^{2,\delta}(\Omega_m)}\lesssim_A2^{2m}(\Vert w\Vert_{C^\delta(\Omega_m)}+\Vert\mathcal H_m\psi\Vert_{C^\delta(\Omega_m)})+2^m\Vert w\Vert_{W^{1,\infty}(\Omega_m)}. 
\]
Interpolating among H\"older spaces we obtain $\Vert w\Vert_{C^\delta(\Omega_m)}\lesssim_AP_I$ and consequently, $\Vert v\Vert_{C^\delta(\Omega_m)}\lesssim_AP_I$. Interpolating it with the $L^\infty$ bound in Lemma~\ref{lem:ImJm} we conclude the estimate of $I_m$ for $k=0$. For $k\geq2$ the result follows from \eqref{eq:DiriCkalp}, \eqref{eq:Gammpt} and induction. The proof of the estimate of $J_m$ is almost parallel. \\
(ii) is a consequence of Lemma~\ref{lem:ImJm} (ii) together with \eqref{eq:diri2}, \eqref{eq:Gammpt} and interpolation among Sobolev spaces. \proofqed

\begin{lemma}\label{lem:calG}
With $s'=s+\sigma$ and $\Omega_m$ above define a map $\mathcal G:(h,u)\mapsto(h^{(K)},u^{(K)})$ via: 
\[
h^{(1)}=h-I_m(h|_{\Gamma_m}),\quad u^{(1)}=u-\nabla J_mI_m(\divg u|_{\Gamma_m})
\]
and
\[
h^{(j)}=h^{(j-2)}-J_m^\frac{j-1}{2}I_m\left(z_{j-1}^{(j-2)}|_{\Gamma_m}\right),\quad u^{(j)}=u^{(j-2)}-\nabla J_m^\frac{j+1}{2}I_m\left(z_j^{(j-2)}|_{\Gamma_m}\right)
\]
for $j=3,5,...,K$. Suppose
\[
\Vert(h,u)\Vert_{C^{\frac{1}{2}+\delta+r}(\Omega_m)}\leq 2^{mr}A_0,\quad r\in[0,K-\frac{1}{2}].
\]
Define
\[
P(h,u):=\sum_{i=0}^K2^{m(\frac{1}{2}+\sigma-i+\delta)}\left(\Vert z_i|_{\Gamma_m}\Vert_{L^\infty(\Gamma_m)}+2^{-m(K-i+\delta)}\Vert z_i|_{\Gamma_m}\Vert_{C^{K-i+\delta}(\Gamma_m)}\right)
\]
and
\[
Q(h,u):=\sum_{i=0}^K2^{m(s'-i-\frac{1}{2})}\left(\Vert z_i|_{\Gamma_m}\Vert_{L^2(\Gamma_m)}+2^{-m(K-i+\frac{1}{2})}\Vert z_i|_{\Gamma_m}\Vert_{H^{K-i+\frac{1}{2}}(\Gamma_m)}\right). 
\]
Assume $P(h,u),Q(h,u)\leq C$. Then we have
\[
\begin{aligned}&2^{m(\frac{1}{2}+\sigma+\delta)}\left(\Vert(h-h^{(K)},u-u^{(K)})\Vert_{L^\infty(\Omega_m)}+2^{-m(K+\delta)}\Vert(h-h^{(K)},u-u^{(K)})\Vert_{C^{K,\delta}(\Omega_m)}\right)\\ &\qquad\lesssim_{A,A_0}P(h,u),
\end{aligned}
\]
\[
2^{ms'}\left(\Vert(h-h^{(K)},u-u^{(K)})\Vert_{L^2(\Omega_m)}+2^{-m(K+1)}\Vert(h-h^{(K)},u-u^{(K)})\Vert_{H^{K+1}(\Omega_m)}\right)\lesssim_{A,A_0}Q(h,u)
\]
and
\[
P(h^{(K)},u^{(K)})\lesssim_{A,A_0}2^{-\frac{m}{2}}P(h,u),\quad Q(h^{(K)},u^{(K)})\lesssim_{A,A_0}2^{-\frac{m}{2}}Q(h,u). 
\]
\end{lemma}
\proofheading{Proof}We first examine the pointwise estimates. By Corollary~\ref{cor:ImJm} (i) we have
\begin{equation}\label{eq:ptdiff}
\begin{aligned}&2^{m(\frac{1}{2}+\sigma+\delta)}\left(\Vert(h-h^{(1)},u-u^{(1)})\Vert_{L^\infty(\Omega_m)}+2^{-m(K+\delta)}\Vert(h-h^{(1)},u-u^{(1)})\Vert_{C^{K,\delta}(\Omega_m)}\right)\\ &\qquad\lesssim_{A,A_0}P(h,u).
\end{aligned}
\end{equation}
It follows that
\begin{equation}\label{eq:ptbd}
\Vert(h^{(1)},u^{(1)})\Vert_{C^{\frac{1}{2}+\delta+r}(\Omega_m)}\lesssim_{A,A_0}2^{mr},\quad r\in[0,K-\frac{1}{2}]
\end{equation}
and that
\begin{equation}\label{eq:pttr}
2^{m(\frac{1}{2}+\sigma-i+\delta)}\left(\Vert z_i^{(1)}|_{\Gamma_m}\Vert_{L^\infty(\Gamma_m)}+2^{-m(K-i+\delta)}\Vert z_i^{(1)}|_{\Gamma_m}\Vert_{C^{K-i+\delta}(\Gamma_m)}\right)\lesssim_{A,A_0}P(h,u) 
\end{equation}
for $i=2,3,\ldots,K$. By construction, $z_0^{(1)}$ and $z_1^{(1)}$ vanish on $\Gamma_m$.

Now assume $j\geq3$ and $(h^{(j-2)},u^{(j-2)})$ satisfies \eqref{eq:ptbd}, \eqref{eq:pttr} for $i=j-1,j,...,K$ 
and verifies
\[
2^{m(\frac{1}{2}+\sigma-i+\delta)}\left(\Vert z_i^{(j-2)}\Vert_{L^\infty(\Gamma_m)}+2^{-m(K-i+\delta)}\Vert z_i^{(j-2)}\Vert_{C^{K-i+\delta}(\Gamma_m)}\right)\lesssim_{A,A_0}2^{-\frac{m}{2}}P(h,u), \ i\leq j-2. 
\]
Following the above arguments we derive \eqref{eq:ptdiff} for $(h^{(j-2)}-h^{(j)},u^{(j-2)}-u^{(j)})$ and \eqref{eq:ptbd}, \eqref{eq:pttr} for $(h^{(j)},u^{(j)})$ and $i=j-1,j,...,K$. Notice that 
\begin{equation}\label{eq:Zi(j)1}
Z_i^{(j-2)}=Z_i^{(j)}\text{ on }\Gamma_m,\quad i\leq j-2,
\end{equation}
where $Z_i^{(j)}$ is the linear part of $z_i^{(j)}$ according to \eqref{eq:zjwj}, and that 
\begin{equation}\label{eq:Zi(j)2}
Z_i^{(j)}+N_i^{(j-2)}=0\text{ on }\Gamma_m,\quad i=j-1,j. 
\end{equation}
Thus the remaining boundary errors $z_i^{(j)}|_{\Gamma_m}$ come only from lower-order terms. Consequently,
\[
2^{m(\frac{1}{2}+\sigma-i+\delta)}\left(\Vert z_i^{(j)}\Vert_{L^\infty(\Gamma_m)}+2^{-m(K-i+\delta)}\Vert z_i^{(j)}\Vert_{C^{K-i+\delta}(\Gamma_m)}\right)\lesssim_{A,A_0}2^{-\frac{m}{2}}P(h,u), \ i\leq j. 
\]
This closes the induction and proves both the bound for $P(h^{(K)},u^{(K)})$ and the pointwise bound for the total correction.

The $L^2$-based estimates are parallel. %This time we use Corollary~\ref{cor:ImJm} (ii), \eqref{eq:bil2}--\eqref{eq:diri2} and \eqref{eq:trint}. Details will not be repeated. \proofqed
Indeed, by Corollary~\ref{cor:ImJm} (ii) we have
\[
2^{ms'}\left(\Vert(h-h^{(1)},u-u^{(1)})\Vert_{L^2(\Omega_m)}+2^{-m(K+1)}\Vert(h-h^{(1)},u-u^{(1)})\Vert_{H^{K+1}(\Omega_m)}\right)\lesssim_AQ(h,u). 
\]
By \eqref{eq:trint} and \eqref{eq:ptbd}, it follows that
\begin{equation}\label{eq:Hstr}
2^{m(s'-i-\frac{1}{2})}\left(\Vert z_i^{(1)}|_{\Gamma_m}\Vert_{L^2(\Gamma_m)}+2^{-m(K-i+\frac{1}{2})}\Vert z_i^{(1)}|_{\Gamma_m}\Vert_{H^{K-i+\frac{1}{2}}(\Gamma_m)}\right)\lesssim_{A,A_0}Q(h,u),\quad i=2,...,K
\end{equation}
and, as mentioned above, $z_0^{(1)},z_1^{(1)}$ exactly vanish on $\Gamma_m$. Then by an induction as above, using \eqref{eq:Zi(j)1}--\eqref{eq:Zi(j)2}, we can show for $j=3,5,...,K$ that
\[\begin{aligned}
2^{ms'}\left(\Vert(h^{(j-2)}-h^{(j)},u^{(j-2)}-u^{(j)})\Vert_{L^2(\Omega_m)}+2^{-m(K+1)}\Vert(h^{(j-2)}-h^{(j)},u^{(j-2)}-u^{(j)})\Vert_{H^{K+1}(\Omega_m)}\right)\\
\lesssim_{A,A_0}Q(h,u),
\end{aligned}\]
that $z_i^{(j)}$ satisfies \eqref{eq:Hstr} for $i=j+1,j+2,...,K$ and that
\[
\begin{aligned}&2^{m(s'-i-\frac{1}{2})}\left(\Vert z_i^{(j)}|_{\Gamma_m}\Vert_{L^2(\Gamma_m)}+2^{-m(K-i+\frac{1}{2})}\Vert z_i^{(j)}|_{\Gamma_m}\Vert_{H^{K-i+\frac{1}{2}}(\Gamma_m)}\right)\\ &\qquad\lesssim_{A,A_0}2^{-\frac{m}{2}}Q(h,u),\ i=0,1,...,j.
\end{aligned}
\]
This completes the proof. \proofqed

We conclude by proving continuity of $\Psi_m$ at a fixed regularization scale. Suppose $(h,u,\Omega)$ and $(h_\epsilon,u_\epsilon,\Omega_\epsilon)$ satisfy Proposition~\ref{prop:Psim} with $s=s'$ and uniform bounds, and that
\[
d\left((h,u,\Omega),(h_\epsilon,u_\epsilon,\Omega_\epsilon)\right):=\Vert(h-h_\epsilon,u-u_\epsilon)\Vert_{L^\infty(\Omega\cap\Omega_\epsilon)}+\Vert\eta_\Gamma-\eta_{\Gamma_\epsilon}\Vert_{L^\infty(\Gamma_*)}\to0. 
\]
For the continuity argument, we may use $\epsilon$ as an upper bound on the displayed distance, with $\epsilon\to0$. Then, for fixed $m$, we have
\begin{equation}\label{eq:Psimconti}
\lim_{\epsilon\to0}d\left(\Psi_m(h,u,\Omega),\Psi_m(h_\epsilon,u_\epsilon,\Omega_\epsilon)\right)=0. 
\end{equation}
The convergence of $\Gamma_m^\epsilon$ follows from the boundedness of the surface frequency cutoff. To prove convergence of $h_m^\epsilon,u_m^\epsilon$, let $I_m^\epsilon,J_m^\epsilon,\mathcal G_\epsilon$ denote the counterparts of $I_m,J_m,\mathcal G$ on $\Omega_m^\epsilon$, and set $\tilde\Omega_m=\Omega_m^\epsilon\cap\Omega_m$. We use the next two lemmas.

\begin{lemma}\label{lem:IJconti}
(i) Suppose $f,f_\epsilon$ are functions on $\Omega_m,\Omega_m^\epsilon$. We have
\[\begin{aligned}
\Vert J_mf-J_m^\epsilon f_\epsilon\Vert_{L^\infty(\tilde\Omega_m)}\lesssim_A2^{-2m}\Vert f-f_\epsilon\Vert_{L^\infty(\tilde\Omega_m)}+\epsilon\big(2^{-m}\Vert f\Vert_{L^\infty(\Omega_m)}+2^{-(1+\delta)m}\Vert f\Vert_{C^\delta(\Omega_m)}\big)\\
+\epsilon\big(2^{-m}\Vert f_\epsilon\Vert_{L^\infty(\Omega_m^\epsilon)}+2^{-(1+\delta)m}\Vert f_\epsilon\Vert_{C^\delta(\Omega_m^\epsilon)}\big). 
\end{aligned}\]
(ii) Suppose $g,g_\epsilon$ are functions on $\Omega_m,\Omega_m^\epsilon$ and $\psi=g|_{\Gamma_m}$, $\psi_\epsilon=g_\epsilon|_{\Gamma_m^\epsilon}$. We have
\[\begin{aligned}
\Vert I_m\psi-I_m^\epsilon\psi_\epsilon\Vert_{L^\infty(\tilde\Omega_m)}&\lesssim_A\Vert g-g_\epsilon\Vert_{L^\infty(\tilde\Omega_m)}+\epsilon^\delta\big(\Vert g\Vert_{C^\delta(\Omega_m)}+\Vert g_\epsilon\Vert_{C^\delta(\Omega_m^\epsilon)}\big)\\
&\quad+\epsilon^\delta\big(2^{m\delta}\Vert\psi\Vert_{L^\infty(\Gamma_m)}+\Vert\psi\Vert_{C^\delta(\Gamma_m)}\big)+\epsilon^\delta\big(2^{m\delta}\Vert\psi_\epsilon\Vert_{L^\infty(\Gamma_m^\epsilon)}+\Vert\psi_\epsilon\Vert_{C^\delta(\Gamma_m^\epsilon)}\big).
\end{aligned}\]
\end{lemma}

\begin{lemma}\label{lem:Gconti}
Suppose $(h,u)$ satisfies the condition in Lemma~\ref{lem:calG} and $(h_\epsilon,u_\epsilon)$ satisfies the analogous condition for $\Omega_m^\epsilon$. Suppose $\Vert(h-h_\epsilon,u-u_\epsilon)\Vert_{L^\infty(\tilde\Omega_m)}\leq\epsilon$. Then for some function $c(\epsilon)$ vanishing as $\epsilon\to0$ we have
\[
\Vert\mathcal G(h,u)-\mathcal G_\epsilon(h_\epsilon,u_\epsilon)\Vert_{L^\infty(\tilde\Omega_m)}\lesssim_{m,A,A_0}c(\epsilon). 
\]
\end{lemma}

We return to \eqref{eq:Psimconti}. Use the notation $\tilde h,\tilde u,h^j,u^j$ from the proof of Proposition~\ref{prop:Psim}, and define the corresponding quantities with subscript $\epsilon$. Lemma~\ref{lem:Phim} gives $\Vert(\tilde h-\tilde h_\epsilon,\tilde u-\tilde u_\epsilon)\Vert_{L^\infty(\tilde\Omega_m)}\lesssim_A\epsilon$. Lemma~\ref{lem:Gconti} then yields
\[
\lim_{\epsilon\to0}\Vert(h^{(N)}-h^{(N)}_\epsilon,u^{(N)}-u^{(N)}_\epsilon)\Vert_{L^\infty(\tilde\Omega_m)}=0
\]
for every $N\in\mathbb N^*$. By Lemma~\ref{lem:calG} we have
\[
\lim_{N\to\infty}\sum_{j\geq N}\Vert(h^{(j)}-h^{(j+1)},u^{(j)}-u^{(j+1)})\Vert_{L^\infty(\Omega_m)}=0
\]
and similar estimate for $h^{(j)}_\epsilon,u^{(j)}_\epsilon$. Letting $\epsilon\to0$ and then $N\to\infty$ we conclude \eqref{eq:Psimconti}. \par\smallskip
\proofheading{Proof of Lemma~\ref{lem:IJconti}}(i) By the maximum principle we have
\[\begin{aligned}
\Vert J_mf-J_m^\epsilon f_\epsilon\Vert_{L^\infty(\tilde\Omega_m)}&\leq2^{-2m}\Vert f-f_\epsilon\Vert_{L^\infty(\tilde\Omega_m)}+\Vert J_mf|_{\tilde\Gamma_m}\Vert_{L^\infty(\tilde\Gamma_m)}+\Vert J_m^\epsilon f_\epsilon|_{\tilde\Gamma_m}\Vert_{L^\infty(\tilde\Gamma_m)}\\
&\lesssim2^{-2m}\Vert f-f_\epsilon\Vert_{L^\infty(\tilde\Omega_m)}+\epsilon\big(\Vert\nabla J_mf\Vert_{L^\infty(\Omega_m)}+\Vert\nabla J^\epsilon_mf_\epsilon\Vert_{L^\infty(\Omega_m^\epsilon)}\big),
\end{aligned}\]
and then by Corollary~\ref{cor:ImJm} the desired estimate follows. \\
(ii) Again by the maximum principle we have
\[
\Vert I_m\psi-I_m^\epsilon\psi_\epsilon\Vert_{L^\infty(\tilde\Omega_m)}\leq\Vert I_m\psi-I_m^\epsilon\psi_\epsilon\Vert_{L^\infty(\tilde\Gamma_m)}. 
\]
Write $I_m\psi-I_m^\epsilon\psi_\epsilon=(I_m\psi-g)+(g-g_\epsilon)+(g_\epsilon-I_m^\epsilon\psi_\epsilon)$ and notice that
\[
\Vert I_m\psi-g\Vert_{L^\infty(\tilde\Gamma_m)}\lesssim\epsilon\big(\Vert\nabla I_m\psi\Vert_{L^\infty(\Omega_m)}+\Vert\nabla I_mg\Vert_{L^\infty(\Omega_m)}\big). 
\]
The same estimate holds for $g_\epsilon-I_m^\epsilon\psi_\epsilon$. Now using Corollary~\ref{cor:ImJm} we complete the proof. \proofqed
\proofheading{Proof of Lemma~\ref{lem:Gconti}}Invoking Lemma~\ref{lem:IJconti} with $g=h,\divg u$, $g_\epsilon=h_\epsilon,\divg u_\epsilon$ we derive
\[
\Vert(h^{(1)}-h^{(1)}_\epsilon,u^{(1)}-u^{(1)}_\epsilon)\Vert_{L^\infty(\tilde\Omega_m)}\lesssim_{m,A,A_0}c(\epsilon). 
\]
Interpolating this with high-order bound using \eqref{eq:Holdinter}, we obtain
\[
\Vert z_j^{(1)}-z_{j,\epsilon}^{(1)}\Vert_{L^\infty(\tilde\Omega_m)}\lesssim_{m,A,A_0}c(\epsilon),\quad j\leq K. 
\]
Repeating the argument we complete the proof. \proofqed

\subsection{Frequency envelopes}\label{sec:fe}

A frequency envelope is a slowly varying, square-summable sequence that controls the energy at each frequency. For $v\in H^s(\Omega)$, set
\begin{equation}\label{eq:feHs}
c_j\left(v;H^s(\Omega)\right)=\sup_l2^{-\delta_e|j-l|}\Vert P_lv\Vert_{H^s(\Omega)},
\end{equation}
where $\delta_e\in(0,\frac{1}{4})$ will be fixed from now on. For $s\geq\frac{1}{2}$ and $(h,u,\Omega)\in\mathbf H^s_{1/2}$, the frequency envelope is defined as
\begin{equation}\label{eq:feHs1/2}
\begin{aligned}
c_j((h,u,\Omega);\mathbf H^s_{1/2})=\sup_l2^{-\delta_e|j-l|}\left(\Vert P_{*l}\eta\Vert_{H^s(\Gamma_*)}+\Vert P_lu\Vert_{H^s(\Omega)}+\Vert(P_lh,P_l\divg u)\Vert_{H^{s+\frac{1}{2}}\times H^{s-\frac{1}{2}}(\Omega)}\right)\\
\mathbf1_{s\in\mathbb N}\left(c_j(z_s\mathbf1_\Omega;H^\frac{1}{2}(\mathbb R^3))+2^{-\delta j}\Vert(h,u,\Omega)\Vert_{\mathbf H^s_{1/2}}\right). 
\end{aligned}
\end{equation}
In these definitions $P_l$, $P_{*l}$ are respectively the frequency localization operators on $\Omega$, $\Gamma_*$. 

\begin{lemma}\label{lem:Phim2}
Let $s\geq\frac{1}{2}$ and $v\in H^s(\Omega)$ with $\divg v\in H^{s-\frac{1}{2}}(\Omega)$. Let $c_j$ be the frequency envelope of $(v,\divg v)\in H^s\times H^{s-\frac{1}{2}}(\Omega)$ according to \eqref{eq:feHs}. Finally, let $v^m=\Phi_mv$ and $\Omega_m$ be a domain bounded by $\Gamma_m\in\Lambda_*$ with $d_H(\Omega,\Omega_m)\leq c2^{-m}$. Then we have: \\
(i)
\[
\Vert(v^m,\divg v^m)\Vert_{H^{s+r}\times H^{s-\frac{1}{2}+r}(\Omega_m)}\lesssim_A2^{mr}c_m,\quad r>\delta_e. 
\]
(ii) 
\[
\Vert v-v^m\Vert_{H^{s-r}(\Omega)}+\mathbf1_{r\leq s-\frac{1}{2}}\Vert\divg v-\divg v^m\Vert_{H^{s-\frac{1}{2}-r}(\Omega)}\lesssim_A2^{-mr}c_m,\quad r\in(\delta_e,s]. 
\]
(iii)
\[
\Vert v^{m+1}-v^m\Vert_{H^{s-r}(\Omega_m)}+\mathbf1_{r\leq s-\frac{1}{2}}\Vert\divg v^{m+1}-\divg v^m\Vert_{H^{s-\frac{1}{2}-r}(\Omega_m)}\lesssim_A2^{-mr}c_m,\quad r\in(\delta_e,s]. 
\]
(iv) If $k-\delta_e>s>\frac{1}{2}+\delta_e$, define
\begin{equation}\label{eq:cjtd1/2}
\begin{aligned}\tilde c_j={}&\sup_l2^{-\delta_e|j-l|}
\Bigl(2^{l(s-\frac{1}{2})}\Vert(v^l-v^{l+1},\divg(v^l-v^{l+1}))\Vert_{H^\frac{1}{2}\times L^2(\Omega)}\\ &\quad+2^{-l(k-s)}\Vert(v^l,\divg v^l)\Vert_{H^k\times H^{k-\frac{1}{2}}(\Omega)}\Bigr). 
\end{aligned}
\end{equation}
Then we have
\[
c_j\lesssim_A\tilde c_j.
\]
\end{lemma}
\proofheading{Proof}The estimates of $v^m,v-v^m,v^m-v^{m+1}$ in (i)-(iii) can be inferred from the proof of Proposition 6.6 in \cite{IE}. For the estimates of the divergence, we use \eqref{eq:hfid}, \eqref{eq:lfid} and note that the convolution operator $\phi_i^m*\cdot$ behaves like $\Phi_m$: 
\[
\Vert[\partial_k,\Phi_m]v\Vert_{H^{s+r}(\Omega_m)}\lesssim_A2^{mr}c_m(v;H^s(\Omega)),\quad r>\delta_e,
\]
\begin{equation}\label{eq:commlfenv}
\Vert[\partial_k,\Phi_m]v\Vert_{H^{s-r}(\Omega_m)}\lesssim_A2^{-mr}c_m(v;H^s(\Omega)),\quad r\in(\delta_e,s].
\end{equation}

To prove (iv), we write $P_jv=P_jv^j+\sum_{l>j}P_j(v^l-v^{l-1})$. By Bernstein's estimate we have
\[
\Vert P_jv\Vert_{H^s(\Omega)}\lesssim_A2^{-j(k-s)}\Vert v^j\Vert_{H^k(\Omega)}+\sum_{l>j}2^{j(s-\frac{1}{2})}\Vert v^l-v^{l-1}\Vert_{H^\frac{1}{2}(\Omega)}\lesssim_A\tilde c_j. 
\]
Similarly, we have
\[
\Vert P_j\divg v\Vert_{H^{s-\frac{1}{2}}(\Omega)}\lesssim_A2^{-j(k-s)}\Vert\divg v^j\Vert_{H^{k-\frac{1}{2}}(\Omega)}+\sum_{l>j}2^{j(s-\frac{1}{2})}\Vert\divg(v^l-v^{l-1})\Vert_{L^2(\Omega)}\lesssim_A\tilde c_j. 
\]
The slow variation of $\tilde c_j$ completes the estimate. \proofqed

\begin{proposition}\label{prop:Psim2}
Suppose $s\geq2$ and $([s-\delta_e,s)\cup(s,s+\delta_e])\cap\mathbb N=\emptyset$. Suppose $(h,u,\Omega)\in\mathbf H^s_{1/2}$ satisfies the pointwise bound and the boundary conditions in Proposition~\ref{prop:Psim} with $s'=s+\frac{1}{2}$. Then for $(h_m,u_m,\Omega_m)=\Psi_m(h,u,\Omega)$, $c_m=c_m((h,u,\Omega);\mathbf H^s_{1/2})$ we have
\[
\Vert\Gamma_m\Vert_{H^{s+r}}+\Vert(h_m,u_m)\Vert_{H^{s+r}(\Omega_m)}+\Vert(\divg u_m,\nabla h_m)\Vert_{H^{s-\frac{1}{2}+r}(\Omega_m)}\lesssim_A2^{mr}c_m,
\]
\[
\Vert\eta-\eta_m\Vert_{L^2(\Gamma_*)}+\Vert(h-h_m,u-u_m)\Vert_{L^2(\Omega\cap\Omega_m)}\lesssim_A2^{-ms}c_m,
\]
and
\[
\Vert(\divg u-\divg u_m,\nabla h-\nabla h_m)\Vert_{L^2(\Omega\cap\Omega_m)}\lesssim_A2^{-m(s-\frac{1}{2})}c_m. 
\]
\end{proposition}
\proofheading{Proof}We refine the proof of Proposition~\ref{prop:Psim}, replacing the global bound $B$ by the frequency envelope $c_m$. Bernstein's estimates give the required bounds for $\eta_m$ and $\eta-\eta_m$. For the interior variables, Lemma~\ref{lem:Phim2} yields
\begin{equation}\label{eq:inthfenv}
\Vert(\tilde h,\tilde u)\Vert_{H^{s+r}(\Omega^{[m]})}+\Vert(\divg\tilde u,\nabla\tilde h)\Vert_{H^{s-\frac{1}{2}+r}(\Omega^{[m]})}\lesssim_A2^{mr}c_m,\quad r>\delta_e
\end{equation}
\begin{equation}\label{eq:intlfenv}
\Vert(h-\tilde h,u-\tilde u)\Vert_{H^{s-r}(\Omega)}+\mathbf1_{r\leq s-\frac{1}{2}}\Vert(\divg u-\divg\tilde u,\nabla h-\nabla\tilde h)\Vert_{H^{s-\frac{1}{2}-r}(\Omega)}\lesssim_A2^{-mr}c_m,\quad r\in(\delta_e,s]
\end{equation}
where $\Omega^{[m]}$ is some neighborhood of $\Omega\cup\Omega_m$. Using \eqref{eq:trint}, carefully choosing $r_1$ to avoid the $\delta_e$-gap in above estimates, we get
\[
\Vert\tilde z_j|_\Gamma\Vert_{L^2(\Gamma)}\lesssim_A2^{-m(s-j)}c_m,\quad j<s. 
\]
If $s\in\mathbb N$, by Lemma~\ref{lem:Phim2}, \eqref{eq:commlfenv}, \eqref{eq:trint} and \eqref{eq:commlf} we have
\[
\Vert\Phi_mz_s|_\Gamma\Vert_{L^2(\Gamma)}\lesssim_A\Vert(\Phi_m-I)(z_s\mathbf1_\Omega)|_{\mathbb R^3\backslash\Omega}\Vert_{L^2(\mathbb R^3\backslash\Omega)}^\frac{1}{2}\Vert\Phi_m(z_s\mathbf1_\Omega)|_{\mathbb R^3\backslash\Omega}\Vert_{H^1(\mathbb R^3\backslash\Omega)}^\frac{1}{2}\lesssim_Ac_m,
\]
\[
\Vert\Phi_mZ_s-\tilde Z_s\Vert_{L^2(\Gamma)}\lesssim_A\Vert\Phi_mZ_s-\tilde Z_s\Vert_{L^2(\Omega)}^\frac{1}{2}\Vert\Phi_mZ_s-\tilde Z_s\Vert_{H^1(\Omega)}^\frac{1}{2}\lesssim_A2^{-\frac{m}{2}}\Vert(h,u,\Omega)\Vert_{\mathbf H^s_{1/2}},
\]
\[
\Vert\Phi_mN_s-N_s\Vert_{L^2(\Gamma)}\lesssim_A\Vert\Phi_mN_s-N_s\Vert_{L^2(\Omega)}^\frac{1}{2}\Vert\Phi_mN_s-N_s\Vert_{H^1(\Omega)}^\frac{1}{2}\lesssim_A2^{-\frac{m}{2}}\Vert(h,u,\Omega)\Vert_{\mathbf H^s_{1/2}}.
\]
We also prove 
\begin{equation}\label{eq:NsminusNstd}
\Vert N_s-\tilde N_s\Vert_{H^1(\Omega)}\lesssim_A2^{-\delta m}\Vert(h,u,\Omega)\Vert_{\mathbf H^s_{1/2}}.
\end{equation}
For $s\geq2$, expand $N_s-\tilde N_s$ into three types of terms. The first has the form $\partial^a\divg(u-\tilde u)\partial^b(h,u,\tilde h,\tilde u)$, with $a+b=s-1$, or the same expression with $\divg(u-\tilde u)$ replaced by $\nabla(h-\tilde h)$. Estimate \eqref{eq:trilest} gives
\[\begin{aligned}
\Vert\partial^a\divg(u-\tilde u)\partial^b(h,u,\tilde h,\tilde u)\Vert_{H^{1-\delta}(\Omega)}\lesssim_A\Vert\divg(u-\tilde u)\Vert_{H^{s-\frac{1}{2}-\delta}(\Omega)}\Vert(h,u,\tilde h,\tilde u)\Vert_{C^{\frac{1}{2}+\delta}(\Omega)}\\
+\Vert\divg(u-\tilde u)\Vert_{L^\infty(\Omega)}\Vert(h,u,\tilde h,\tilde u)\Vert_{H^s(\Omega)}. 
\end{aligned}\]
The right-hand side is dominated by $2^{-\delta m}\Vert(h,u,\Omega)\Vert_{\mathbf H^s_{1/2}}$ thanks to Lemma~\ref{lem:Phim} and the estimate
\[
\Vert v-\Phi_mv\Vert_{C^{\alpha-r}(\Omega)}\lesssim_A2^{-mr}\Vert v\Vert_{C^\alpha(\Omega)},\quad r\in[0,\alpha],
\]
which can be found in \cite{IE}. The second kind is $\partial^a(\divg u,\nabla h,\divg\tilde u,\nabla\tilde h)\partial^b(h-\tilde h,u-\tilde u)$. The same argument gives: 
\[\begin{aligned}&
\Vert\partial^a(\divg u,\nabla h,\divg\tilde u,\nabla\tilde h)\partial^b(h-\tilde h,u-\tilde u)\Vert_{H^1(\Omega)}\\ &\lesssim_A\Vert(\divg u,\nabla h,\divg\tilde u,\nabla\tilde h)\Vert_{H^{s-\frac{1}{2}}(\Omega)}\Vert(h-\tilde h,u-\tilde u)\Vert_{C^\frac{1}{2}(\Omega)}\\
&\quad\ +\Vert(\divg u,\nabla h,\divg\tilde u,\nabla\tilde h)\Vert_{C^\delta(\Omega)}\Vert(h-\tilde h,u-\tilde u)\Vert_{H^{s-\delta}(\Omega)}\\
&\lesssim_A2^{-\delta m}\Vert(h,u,\Omega)\Vert_{\mathbf H^s_{1/2}}.
\end{aligned}\]
The third kind is $\partial^a(h-\tilde h,u-\tilde u)\partial^{b_1}(h,u,\tilde h,\tilde u)\cdots\partial^{b_k}(h,u,\tilde h,\tilde u)$ with $k\geq2$, $a+b_1+\cdots+b_k=s$. The estimate is as follows:
\[\begin{aligned}
&\Vert\partial^a(h-\tilde h,u-\tilde u)\partial^{b_1}(h,u,\tilde h,\tilde u)\cdots\partial^{b_k}(h,u,\tilde h,\tilde u)\Vert_{H^1(\Omega)}\\
&\lesssim_A\Vert(h-\tilde h,u-\tilde u)\Vert_{H^{s+1-\frac{k}{2}-k\delta}(\Omega)}\Vert(h,u,\tilde h,\tilde u)\Vert_{C^{\frac{1}{2}+\delta}(\Omega)}^k\\
&\quad\ +\Vert(h-\tilde h,u-\tilde u)\Vert_{C^\frac{1}{2}(\Omega)}\Vert(h,u,\tilde h,\tilde u)\Vert_{C^\frac{1}{2}(\Omega)}^{k-1}\Vert(h,u,\tilde h,\tilde u)\Vert_{H^{s+1-\frac{k}{2}}(\Omega)}\\
&\lesssim_A2^{-\delta m}\Vert(h,u,\Omega)\Vert_{\mathbf H^s_{1/2}}, 
\end{aligned}\]
where in the last line we have crucially used that $k\geq2$. This completes the justification of \eqref{eq:NsminusNstd}. Combining these we obtain
\[
\Vert\tilde z_s|_\Gamma\Vert_{L^2(\Gamma)}\lesssim_Ac_m. 
\]
Applying the fundamental theorem of calculus with the refined bound \eqref{eq:intlfenv} gives
\[
\Vert\tilde z_j|_{\Gamma_m}\Vert_{L^2(\Gamma_m)}\lesssim_A2^{-m(s-j)}c_m,\quad j\leq s. 
\]
Using \eqref{eq:trint}, carefully choosing $r_0$ to avoid the $\delta_e$-gap in \eqref{eq:inthfenv}, we extend above estimate to $j\leq K$. Using \eqref{eq:inthfenv}, \eqref{eq:tr2} we obtain
\[
\Vert\tilde z_j|_{\Gamma_m}\Vert_{H^{K-j+\frac{1}{2}}(\Gamma_m)}\lesssim_A2^{m(K+\frac{1}{2}-s)}c_m,\quad j\leq K. 
\]
These together imply $Q(\tilde h,\tilde u)\lesssim_Ac_m$, and invoking Lemma~\ref{lem:calG} we complete the proof. \proofqed

\begin{corollary}\label{cor:unil2}
Suppose $s\geq2+\delta$, $(h^l,u^l,\Omega^l)\in\mathbf H^s_{1/2}$ satisfies the compatibility condition in Proposition~\ref{prop:Psim} with $s'=s+\frac{1}{2}$, and $(h^l,u^l,\Omega^l)\to(h,u,\Omega)$ in $\mathbf H^s_{1/2}$. Then $c^l_m=c_m((h^l,u^l,\Omega^l);\mathbf H^s_{1/2})$ is uniformly $\ell^2_m$: 
\[
\lim_{m\to\infty}\sup_l\sum_{j>m}(c^l_j)^2=0. 
\]
\end{corollary}
\proofheading{Proof}The convergence of $\eta_{\Gamma^l}$, $\mathbf1_{s\in\mathbb N}z_s^l\mathbf1_{\Omega^l}$ imply the uniform $\ell^2$-summability of $c_m(\eta_{\Gamma^l};H^s(\Gamma_*))$, $c_m(\mathbf1_{s\in\mathbb N}z_s^l\mathbf1_{\Omega^l};H^\frac{1}{2}(\mathbb R^3))$. Similarly, if $\tilde h^l,\tilde u^l$ are extensions of $h^l,u^l$ verifying
\[
\Vert(\tilde h^l-\tilde h,\tilde u^l-\tilde u)\Vert_{H^s(\mathbb R^3)}+\Vert(\tilde h^l-\tilde h,\divg\tilde u^l-\divg\tilde u)\Vert_{H^{s+\frac{1}{2}}\times H^{s-\frac{1}{2}}(\mathbb R^3)}\to0, 
\]
we have that $c_m\big(\tilde h^l,\tilde u^l,\divg\tilde u^l;H^{s+\frac{1}{2}}\times H^s\times H^{s-\frac{1}{2}}(\mathbb R^3)\big)$ are uniformly $\ell^2_m$. At this stage we cannot conclude the uniform $\ell^2$-summability of $c_m\big(h^l,u^l,\divg u^l;H^{s+\frac{1}{2}}\times H^s\times H^{s-\frac{1}{2}}(\Omega^l)\big)$, because the two envelopes involve different frequency localization operators. To compare these envelopes, we set $(\tilde h^l_m,\tilde u^l_m)=\Phi_m(\tilde h^l,\tilde u^l)$. Define $\tilde c_m(\tilde u^l)$ via \eqref{eq:cjtd1/2}. Similarly, define
\[
\tilde c_m(\tilde h^l)=\sup_j2^{-\delta_e|j-m|}\left(2^{j(s+\frac{1}{2})}\Vert\tilde h^l_j-\tilde h^l_{j+1}\Vert_{L^2(\Omega^l)}+2^{-j(k-s-\frac{1}{2})}\Vert\tilde h^l_j\Vert_{H^k(\Omega^l)}\right). 
\]
Then Lemma~\ref{lem:Phim2} (i), (iii) imply $\tilde c_m(\tilde h^l)+\tilde c_m(\tilde u^l)\lesssim_Ac_m\big(\tilde h^l,\tilde u^l,\divg\tilde u^l;H^{s+\frac{1}{2}}\times H^s\times H^{s-\frac{1}{2}}(\mathbb R^3)\big)$, while Lemma~\ref{lem:Phim2} (iv) implies $c_m\big(h^l,u^l,\divg u^l;H^{s+\frac{1}{2}}\times H^s\times H^{s-\frac{1}{2}}(\Omega^l)\big)\lesssim_A\tilde c_m(\tilde h^l)+\tilde c_m(\tilde u^l)$. This completes the proof. \proofqed

\section{Distance estimates}\label{sec:dist}

We compare solutions whose fluid domains need not coincide. The distance $D$ controls the bulk state and the interface at the $L^2$ level and yields uniqueness. The partial $H^2$ distance $D_2$ retains the additional acoustic derivatives needed in the rough-data limit.

\subsection{The \texorpdfstring{$L^2$}{L2} distance}\label{sub:l2distance}

For two states $S=(h,u,\Omega)$ and $S^m=(h^m,u^m,\Omega^m)$, define
\begin{equation}\label{eq:L2distfunc}
D(S,S^m)=\int_{\tilde\Omega}(h-h^m)^2+|u-u^m|^2dx+\int_{\tilde\Gamma}\tilde a^{-1}(h-h^m)^2dS,
\end{equation}
where
\[
\tilde\Omega=\Omega\cap\Omega^m,\quad\tilde\Gamma=\partial\tilde\Omega,\quad\tilde a=a\mathbf1_{\Gamma\cap\Omega^m}+a^m\mathbf1_{\Gamma^m\cap\Omega}. 
\]
\par\smallskip\noindent\textbf{Remark.}\ The definition makes sense because $\tilde\Gamma$, being the graph of the Lipschitz continuous function $\eta_\Gamma\wedge\eta_{\Gamma^m}$, is itself a Lipschitz boundary. 

\begin{theorem}[Propagation of the $L^2$ distance]\label{thm:L2dist}
Suppose $S,S^m$ are two solutions to \eqref{eq:CE} in $C^2$. Then we have
\[
\frac{d}{dt}D(S,S^m)\lesssim_{A,A^m}(B+B^m)D(S,S^m), 
\]
where
\[
A=\Vert\Gamma\Vert_{C^{1,\delta}},\quad B=\Vert(u,\divg u,\nabla h)\Vert_{C^1(\Omega)},
\]
and $A^m,B^m$ are defined accordingly. 
\end{theorem}
The proof uses the following transport identities. They follow from the Reynolds transport formula.

\begin{lemma}\label{lem:tpid}
(i) Given a velocity field $v$ defined on a time-dependent domain $\mathcal D_t$ with Lipschitz boundary flowing with normal velocity $v_b$, we have
\[
\frac{d}{dt}\int_{\mathcal D_t}f\,dx=\int_{\mathcal D_t}D_tf+f\divg v\,dx+\int_{\partial\mathcal D_t}f(v_b-v\cdot n)\,dS. 
\]
(ii) If the time-dependent surface $\mathcal S_t$ flows with velocity $v$, then we have
\[
\frac{d}{dt}\int_{\mathcal S_t}f\,dS=\int_{\mathcal S_t}D_tf+f\divg_{\mathcal S_t}v\,dS. 
\]
\end{lemma}
\par\smallskip\noindent\textbf{Remark.}\ Part (i) applies to the intersection domain because $\tilde\Gamma_t$ is Lipschitz. Its outward normal is
\[
\tilde n(x)=\begin{cases}n(x)&x\in\Gamma\cap\overline{\Omega^m}\\n^m(x)&x\in\Gamma^m\cap\Omega\end{cases}
\]
for almost every $x\in\tilde\Gamma$. If we choose $v=u$ (or equally, $u^m$), then we have
\begin{equation}\label{eq:vbvdiff}
|v_b-v\cdot n|\leq|u-u^m|. 
\end{equation}
Part (ii) requires a Lipschitz flow velocity on $\mathcal S_t$, which need not exist for the intersection boundary $\tilde\Gamma_t$. We instead split the boundary integral in \eqref{eq:L2distfunc} as
\[
\int_\Gamma a^{-1}(h^m-h)^2\mathbf1_{\Gamma\cap\Omega^m}\,dS+\int_{\Gamma^m}(a^m)^{-1}(h-h^m)^2\mathbf1_{\Gamma^m\cap\Omega}\,dS.
\]
The boundary conditions imply that $(h^m-h)^2\mathbf1_{\Gamma\cap\Omega^m}$ and $(h-h^m)^2\mathbf1_{\Gamma^m\cap\Omega}$ are Lipschitz. We may therefore apply part (ii) to each integral on its original moving boundary.
\proofheading{Proof of Theorem~\ref{thm:L2dist}}We first calculate the temporal derivative of $I=\int_{\tilde\Omega}(h-h^m)^2+|u-u^m|^2dx$. By Lemma~\ref{lem:tpid} (i) we have
\[\begin{aligned}
\frac{1}{2}\frac{dI}{dt}=\int_{\tilde\Omega}(h-h^m)(\partial_t+u\cdot\nabla)(h-h^m)+(u-u^m)\cdot(\partial_t+u\cdot\nabla)(u-u^m)\,dx\\
+\int_{\tilde\Gamma}\frac{1}{2}((h-h^m)^2+|u-u^m|^2)(v_b-u\cdot\tilde n)\,dS+\mathcal O(D). 
\end{aligned}\]
Throughout the proof $\mathcal O(D)$ stands for $O((B+B^m)D(S,S^m))$. By \eqref{eq:vbvdiff} and Lemma~\ref{lem:trpqr}, the boundary integral can be dominated by $\Vert u-u^m\Vert_{C^1(\tilde\Omega)}\Vert u-u^m\Vert_{L^2(\tilde\Omega)}^2$, which is $\mathcal O(D)$. For the interior integral, we write $(\partial_t+u\cdot\nabla)h^m=(\partial_t+u^m\cdot\nabla)h^m+(u-u^m)\cdot\nabla h^m$,  $(\partial_t+u\cdot\nabla)u^m=(\partial_t+u^m\cdot\nabla)u^m+(u-u^m)\cdot\nabla u^m$, use \eqref{eq:CE} and use the divergence theorem, yielding that
\[\begin{aligned}
\frac{1}{2}\frac{dI}{dt}&=\int_{\tilde\Omega}(h-h^m)\divg(u^m-u)+(u-u^m)\cdot\nabla(h^m-h)\,dx+\mathcal O(D)\\
&=\int_{\Gamma\cap\Omega^m}(h-h^m)(u^m-u)\cdot n\,dS+\int_{\Gamma^m\cap\Omega}(h-h^m)(u^m-u)\cdot n^m\,dS+\mathcal O(D). 
\end{aligned}\]

We next calculate the temporal derivative of $J=\int_{\tilde\Gamma}\tilde a^{-1}(h-h^m)^2\,dS$. By Lemma~\ref{lem:tpid} and the remark following it, we have
\[\begin{aligned}
\frac{1}{2}\frac{dJ}{dt}=&\int_{\Gamma\cap\Omega^m}a^{-1}(h-h^m)(\partial_t+u\cdot\nabla)(h-h^m)\,dS\\
&+\int_{\Gamma^m\cap\Omega}(a^m)^{-1}(h-h^m)(\partial_t+u^m\cdot\nabla)(h-h^m)\,dS+\mathcal O(D).
\end{aligned}\]
From \eqref{eq:CE} we see
\[\begin{aligned}
&(\partial_t+u\cdot\nabla)(h-h^m)=\divg u^m-\divg u+(u^m-u)\cdot\nabla h^m,\\
&(\partial_t+u^m\cdot\nabla)(h-h^m)=\divg u^m-\divg u+(u^m-u)\cdot\nabla h. 
\end{aligned}\]
From the boundary conditions $\divg u|_\Gamma=0$, $\divg u^m|_{\Gamma^m}=0$ and the Taylor condition we see $|\divg u^m-\divg u|\lesssim B|h-h^m|$ on $\tilde\Gamma$. Thus, we have
\[
\frac{1}{2}\frac{dJ}{dt}=\int_{\Gamma\cap\Omega^m}a^{-1}(h-h^m)(u^m-u)\cdot\nabla h^m\,dS+\int_{\Gamma^m\cap\Omega}(a^m)^{-1}(h-h^m)(u^m-u)\cdot\nabla h\,dS+\mathcal O(D).
\]
By $a^{-1}\nabla h^m=a^{-1}\nabla(h^m-h)-n$ on $\Gamma\cap\Omega^m$ and $(a^m)^{-1}\nabla h=(a^m)^{-1}\nabla(h-h^m)-n^m$ on $\Gamma^m\cap\Omega$, we obtain
\[\begin{aligned}
\frac{1}{2}\frac{dI}{dt}+\frac{1}{2}\frac{dJ}{dt}=&\int_{\Gamma\cap\Omega^m}a^{-1}(h-h^m)(u^m-u)\cdot\nabla(h^m-h)\,dS\\
&+\int_{\Gamma^m\cap\Omega}(a^m)^{-1}(h-h^m)(u^m-u)\cdot\nabla(h-h^m)\,dS+\mathcal O(D). 
\end{aligned}\]
Therefore, using H\"older's inequality, we derive
\[
\frac{d}{dt}\frac{1}{2}D(S,S^m)\lesssim_A\Vert h-h^m\Vert_{L^2(\tilde\Gamma)}\Vert u-u^m\Vert_{L^3(\tilde\Gamma)}\Vert\nabla(h-h^m)\Vert_{L^5(\tilde\Gamma)}^\frac{5}{6}\Vert\nabla(h-h^m)\Vert_{L^\infty(\tilde\Gamma)}^\frac{1}{6}. 
\]
By Lemma~\ref{lem:trpqr} we have $\Vert u-u^m\Vert_{L^3(\tilde\Gamma)}\lesssim_AB^\frac{1}{3}\Vert u-u^m\Vert_{L^2(\tilde\Omega)}^\frac{2}{3}$ and
\[\begin{aligned}
\Vert\nabla(h-h^m)\Vert_{L^5(\tilde\Gamma)}\lesssim_AB^\frac{1}{5}\Vert\nabla(h-h^m)\Vert_{L^4(\tilde\Omega)}^\frac{4}{5}\lesssim_AB^\frac{3}{5}\Vert h-h^m\Vert_{L^2(\tilde\Omega)}^\frac{2}{5}. 
\end{aligned}\]
This completes the proof of the theorem. \proofqed

\begin{lemma}\label{lem:trpqr}
(i) Let $1<p<\infty$, $1<q\leq\infty$, and $r=(p-1)q/(q-1)\geq1$, with $r=p-1$ when $q=\infty$. Then
\[
\Vert f\Vert_{L^p(\tilde\Gamma)}\lesssim_A\Vert f\Vert_{L^r(\tilde\Omega)}^{1-\frac{1}{p}}\Vert f\Vert_{W^{1,q}(\tilde\Omega)}^\frac{1}{p}. 
\]
(ii) We have
\[
\Vert\nabla f\Vert_{L^4(\tilde\Omega)}\lesssim_A\Vert f\Vert_{L^2(\tilde\Omega)}^\frac{1}{2}\Vert f\Vert_{C^2(\tilde\Omega)}^\frac{1}{2}. 
\]
\end{lemma}
\proofheading{Proof}(i) Since $\Gamma,\Gamma^m\in\Lambda_*$, there is a vector field $X$, depending only on the collar $\Lambda_*$, with $X\cdot\tilde n\sim1$ and $|\nabla X|\lesssim_A1$. The divergence theorem gives
\[
\Vert f\Vert_{L^p(\tilde\Gamma)}^p\sim\int_{\tilde\Gamma}|f|^pX\cdot\tilde n\,dS=\int_{\tilde\Omega}(\divg X)|f|^p+p|f|^{p-2}fX\cdot\nabla f\,dx. 
\]
By H\"older's inequality the right-hand side is dominated by
\[
\Vert f\Vert_{L^r(\tilde\Omega)}^{p-1}\Vert f\Vert_{W^{1,q}(\tilde\Omega)}. 
\]
(ii) By the divergence theorem we have
\[
\Vert\nabla f\Vert_{L^4(\tilde\Omega)}^4=\int_{\tilde\Omega}|\nabla f|^2\nabla f\cdot\nabla f\,dx=\int_{\tilde\Gamma}f|\nabla f|^2\nabla_{\tilde n}f\,dS-\int_{\tilde\Omega}f\divg(|\nabla f|^2\nabla f)\,dx. 
\]
The interior integral is bounded by $3\Vert f\Vert_{L^2(\tilde\Omega)}\Vert f\Vert_{C^2(\tilde\Omega)}\Vert\nabla f\Vert_{L^4(\tilde\Omega)}^2$. By (i), the boundary integral is bounded by 
\[
\Vert\nabla f\Vert_{L^5(\tilde\Gamma)}^3\Vert f\Vert_{L^\frac{5}{2}(\tilde\Gamma)}\lesssim\Vert\nabla f\Vert_{L^4(\tilde\Omega)}^\frac{12}{5}\Vert\nabla f\Vert_{C^1(\tilde\Omega)}^\frac{3}{5}\Vert f\Vert_{L^2(\tilde\Omega)}^\frac{3}{5}\Vert\nabla f\Vert_{L^4(\tilde\Omega)}^\frac{2}{5}. 
\]
Dividing by the common factors gives the desired estimate; the case in which a factor vanishes is immediate. \proofqed

\subsection{Partial \texorpdfstring{$H^2$}{H2} distance}

Let $\mathbf v$ denote the pair $(h,u)$. With the notations in the last subsection, define
\begin{equation}\label{eq:H2distfunc}
D_2(\mathbf v,\mathbf v^m)=\int_{\tilde\Omega}(z_2-z_2^m)^2+|w_2-w_2^m|^2\,dx. 
\end{equation}
%Of course, $z_2=z_2(\mathbf v),w_2=w_2(\mathbf v)$ and $z_2^m=z_2(\mathbf v^m),w_2^m=w_2(\mathbf v^m)$. 

\begin{proposition}\label{prop:H2dist}
Let $k>s>3$, $M>0$ and $m\gg1$. Suppose $(h,u,\Omega)$ is a solution in $\mathbf H^k_1$ with the following properties. \\
(i) Pointwise bound: 
\[
\Vert\Gamma\Vert_{C^{1,\delta}}+\Vert u\Vert_{C^{\frac{1}{2}+\delta}(\Omega)}+\Vert h\Vert_{C^{1,\delta}(\Omega)}\leq A\text{ on }[0,T],
\]
and
\[
\Vert\Gamma_t\Vert_{C^2}+\Vert u\Vert_{C^1(\Omega_t)}+\Vert(\divg u,\nabla h)\Vert_{C^1(\Omega_t)}+2^{-m}\left(\Vert\nabla u\Vert_{C^1(\Omega_t)}+\Vert(\Delta h,\nabla\divg u)\Vert_{C^1(\Omega_t)}\right)\leq B(t).
\]
(ii) Regularity bound: 
\[
\Vert(h,u,\Omega)\Vert_{\mathbf H^k}\leq 2^{m(k-s)}M,\quad\Vert(h,u,\Omega)\Vert_{\mathbf H_1^k}\leq 2^{m(k+\frac{1}{2}-s)}M\quad\text{on }[0,T]. 
\]

Suppose $(h^m,u^m,\Omega^m)$ is another solution satisfying (i), (ii). Further assume\\
(iii) Distance estimate: 
\[
D(S,S^m)\leq2^{-2ms}M^2\text{ on }[0,T].
\]
Then we have
\[
\frac{d}{dt}D_2(\mathbf v,\mathbf v^m)\lesssim_AB^2\left(D_2(\mathbf v,\mathbf v^m)+2^{-m(2s-3)}M^2\right). 
\]
\end{proposition}
\proofheading{Proof}From the distance estimate we see
\begin{equation}\label{eq:vvmL2}
\Vert\mathbf v-\mathbf v^m\Vert_{L^2(\tilde\Omega)}\leq2^{-ms}M,
\end{equation}
\begin{equation}\label{eq:etaL2}
\Vert\eta_\Gamma-\eta_{\Gamma^m}\Vert_{L^2(\Gamma_*)}\lesssim_A2^{-ms}M. 
\end{equation}
By Lemma~\ref{lem:Phim} and the regularity bound, we have
\begin{align}
\Vert\nabla(\mathbf v-\mathbf v^m)\Vert_{L^2(\tilde\Omega)}&\leq\Vert\nabla(I-\Phi_m)\mathbf v\Vert_{L^2(\Omega)}+\Vert\nabla(I-\Phi_m)\mathbf v^m\Vert_{L^2(\Omega^m)}+\Vert\nabla\Phi_m(\mathbf v-\mathbf v^m)\Vert_{L^2(\tilde\Omega)}\nonumber\\
&\lesssim_A2^{-m(s-1)}M. \label{eq:vvmH1}
\end{align}
It follows that
\begin{equation}\label{eq:A2v}
\Vert\Delta(h-h^m)\Vert_{L^2(\tilde\Omega)}+\Vert\nabla\divg(u-u^m)\Vert_{L^2(\tilde\Omega)}\lesssim_AD_2(\mathbf v,\mathbf v^m)^\frac{1}{2}+2^{-m(s-1)}BM. 
\end{equation}
This estimate does not directly control the full Hessian because the difference lacks a common Dirichlet boundary condition. We therefore align the domains. The distance bound \eqref{eq:etaL2}, the regularity bound, and interpolation give a diffeomorphism $\psi^m:\overline\Omega\to\overline{\Omega^m}$ satisfying
\begin{equation}\label{eq:Dpsim}
\Vert\psi^m-\text{Id}\Vert_{H^r(\Omega)}+\Vert(\psi^m)^{-1}-\text{Id}\Vert_{H^r(\Omega^m)}\lesssim_A2^{-m(s+\frac{1}{2}-r)}M,\quad r\in(\frac{1}{2},k+\frac{1}{2}],
\end{equation}
and satisfying $\Vert\nabla\psi^m\Vert_{C^\delta(\Omega)}+\Vert\nabla(\psi^m)^{-1}\Vert_{C^\delta(\Omega^m)}\lesssim_A1$. 

We prove that 
\begin{equation}\label{eq:cochdiff}
\begin{aligned}
&\Vert\nabla(h-h^m(\psi^m))\Vert_{L^2(\Omega)}\lesssim_A2^{-m(s-1)}BM,\\
&\Vert\Delta(h-h^m(\psi^m))\Vert_{L^2(\Omega)}\lesssim_AD_2(\mathbf v,\mathbf v^m)^\frac{1}{2}+2^{-m(s-\frac{3}{2})}M+2^{-m(s-1-\delta)}BM. 
\end{aligned}
\end{equation}
%By \eqref{eq:Dpsim} and the chain rule we have $\Vert\nabla h^m(\psi^m)-\nabla h^m\circ\psi^m\Vert_{L^2(\Omega)}\lesssim_A2^{-m(s-1)}BM$. Next, 
Write
\begin{equation}\label{eq:hhmdec}
\begin{aligned}
\nabla h-\nabla h^m(\psi^m)=(I-\Phi_m)\nabla h-(I-\Phi_m)(\nabla h^m\circ\psi^m)-\Phi_m(\nabla h^m\circ\psi^m-\nabla h^m)\\
+\Phi_m\nabla(h-h^m)-(\nabla\psi^m-I)\nabla h^m\circ\psi^m. 
\end{aligned}
\end{equation}
The first term is $O_{L^2(\Omega)}(2^{-m(s-1)}M)$, due to Lemma~\ref{lem:Phim} and the regularity bound. The control of the second term is similar, involving Moser's estimate and \eqref{eq:Dpsim}. The third and the fourth terms, observed on $\Omega$, depends only on $\nabla h^m\circ\psi^m-\nabla h^m,\nabla(h-h^m)$ restricted on a subset $\Omega'$ which is roughly $2^{-m}$ away from $\Gamma\cup\Gamma^m$. Thus, for the third term we use the fundamental theorem of calculus, yielding that
\[
\Vert\Phi_m(\nabla h^m\circ\psi^m-\nabla h^m)\Vert_{L^2(\Omega)}\lesssim_AB\Vert\psi^m-\text{Id}\Vert_{L^2(\Omega)}\lesssim_A2^{-m(s-\delta)}BM, 
\] 
and for the fourth term we read from \eqref{eq:vvmH1} that $\Vert\Phi_m\nabla(h-h^m)\Vert_{L^2(\Omega)}\lesssim_A2^{-m(s-1)}M$. The fifth term is $O_{L^2(\Omega)}(2^{-m(s-1)}BM)$ by the chain rule and \eqref{eq:Dpsim}. Combining these we justify the first statement in \eqref{eq:cochdiff}. For the second statement, again use the chain rule: 
\[
\Delta h^m(\psi^m)=\nabla\psi^m_i\cdot\nabla\psi^m_j\,\partial_i\partial_jh^m\circ\psi^m+\Delta\psi^m\cdot\nabla h^m\circ\psi^m. 
\]
By \eqref{eq:Dpsim} we have
\[
\Vert\nabla\psi^m_i\cdot\nabla\psi^m_j-\delta_{ij}\Vert_{L^2(\Omega)}\lesssim_A2^{-m(s-\frac{1}{2})}M,\quad\Vert\Delta\psi^m\Vert_{L^2(\Omega)}\lesssim_A2^{-m(s-\frac{3}{2})}M. 
\]
It follows that $\Vert\Delta h^m(\psi^m)-\Delta h^m\circ\psi^m\Vert_{L^2(\Omega)}\lesssim_A2^{-m(s-\frac{3}{2})}M+2^{-m(s-\frac{1}{2})}BM$. Now replace $\nabla h,\nabla h^m$ in \eqref{eq:hhmdec} by $\Delta h,\Delta h^m$. Arguing as above and using \eqref{eq:A2v}, we obtain \eqref{eq:cochdiff}; this is the first place the pointwise bound on $2^{-m}\nabla\Delta h^m$ is used. 

We now estimate $\nabla^2(h-h^m(\psi^m))$ using \eqref{eq:cochdiff} and the following refinement of \eqref{eq:diri2}:
\begin{equation}\label{eq:diri3}
\Vert f\Vert_{H^2(\Omega)}\lesssim_A\Vert\Delta f\Vert_{L^2(\Omega)}+\Vert\Gamma\Vert_{C^2}\Vert f\Vert_{H^1(\Omega)}, 
\end{equation}
This follows by flattening the boundary and freezing the coefficients. Hence,
\[
\Vert\nabla^2(h-h^m(\psi^m))\Vert_{L^2(\Omega)}\lesssim_AD_2(\mathbf v,\mathbf v^m)^\frac{1}{2}+2^{-m(s-\frac{3}{2})}M+2^{-m(s-1-\delta)}B^2M. 
\]
Inverting the proof of \eqref{eq:cochdiff}, we see
\begin{equation}\label{eq:D2hdiff}
\Vert\nabla^2(h-h^m)\Vert_{L^2(\tilde\Omega)}\lesssim_AD_2(\mathbf v,\mathbf v^m)^\frac{1}{2}+2^{-m(s-\frac{3}{2})}M+2^{-m(s-1-\delta)}B^2M. 
\end{equation}

One last ingredient is the following estimate: 
\begin{equation}\label{eq:w2H1diff}
\Vert\nabla(w_2-w_2^m)\Vert_{L^2(\tilde\Omega)}\lesssim_A2^mD_2(\mathbf v,\mathbf v^m)^\frac{1}{2}+2^{-m(s-\frac{5}{2})}M+2^{-m(s-2)}BM. 
\end{equation}
The nonlinear terms are estimated as in \eqref{eq:vvmH1}. For the linear part, note the regularity bound
\[
\Vert\nabla\divg u\Vert_{H^{k-1}(\Omega)}+\Vert\nabla\divg u^m\Vert_{H^{k-1}(\Omega^m)}\lesssim_A2^{m(k+\frac{1}{2}-s)}M. 
\]
With the decomposition
\[
\nabla\divg(u-u^m)=(I-\Phi_m)\nabla\divg u-(I-\Phi_m)\nabla\divg u^m+\Phi_m\nabla\divg(u-u^m), 
\]
by Lemma~\ref{lem:Phim} and by \eqref{eq:A2v}, we obtain
\[
\Vert\nabla^2\divg(u-u^m)\Vert_{L^2(\tilde\Omega)}\lesssim_A2^mD_2(\mathbf v,\mathbf v^m)^\frac{1}{2}+2^{-m(s-\frac{5}{2})}M+2^{-m(s-2)}BM, 
\]
which leads to \eqref{eq:w2H1diff}. 

We can now estimate $d\,D_2(\mathbf v,\mathbf v^m)/dt$. Set $\mathcal A\mathbf v=(\divg u,\nabla h)$. Direct calculation gives
\[\begin{aligned}
&D_tz_2+\divg w_2=F(\nabla\mathbf v,\nabla\mathcal A\mathbf v)+f(\nabla\mathbf v),&\text{ in }\Omega,\\
&D_tw_2+\nabla z_2=G(\nabla\mathbf v,\nabla\mathcal A\mathbf v)+g(\nabla\mathbf v),&\text{ in }\Omega,
\end{aligned}\]
where $F(\cdot,\cdot),G(\cdot,\cdot)$ are bilinear forms and $f,g$ are cubic polynomials. The variables $z_2^m,w_2^m$ satisfy the corresponding system in $\Omega^m$ (with $D_t$ replaced by $\partial_t+u^m\cdot\nabla$). Subtracting the second system from the first, we derive the following difference system in $\tilde\Omega$: 
\begin{equation}\label{eq:diffsys}
\begin{aligned}
&D_t(z_2-z_2^m)+\divg(w_2-w_2^m)=(u^m-u)\cdot\nabla z_2^m+F_1+F_2+\tilde f,\\
&D_t(w_2-w_2^m)+\nabla(z_2-z_2^m)=(u^m-u)\cdot\nabla w_2^m+G_1+G_2+\tilde g,
\end{aligned}
\end{equation}
where
\[
F_1=F(\nabla(\mathbf v-\mathbf v^m),\nabla\mathcal A\mathbf v^m),\quad F_2=F(\nabla\mathbf v,\nabla\mathcal A(\mathbf v-\mathbf v^m)),
\]
\[
G_1=G(\nabla(\mathbf v-\mathbf v^m),\nabla\mathcal A\mathbf v^m),\quad G_2=G(\nabla\mathbf v,\nabla\mathcal A(\mathbf v-\mathbf v^m)),
\]
and $\tilde f,\tilde g$ are quadratic in $\nabla\mathbf v,\nabla\mathbf v^m$ while linear in $\nabla(\mathbf v-\mathbf v^m)$. Pair \eqref{eq:diffsys} with $(z_2-z_2^m,w_2-w_2^m)$ and use Lemma~\ref{lem:tpid} on $\tilde\Omega$. The left-hand side is
\[
\frac{d}{dt}\frac{1}{2}D_2(\mathbf v,\mathbf v^m)+\mathcal I+O(BD_2),
\]
where $\mathcal I$ is bounded by
\begin{equation}\label{eq:calIbd}
\int_{\tilde\Gamma}(|u-u^m|+1)\,|w_2-w_2^m|\,|z_2-z_2^m|\,dS.
\end{equation}
The right-hand side is dominated by
\begin{equation}\label{eq:rhsbd}
\begin{aligned}&D_2(\mathbf v,\mathbf v^m)^\frac{1}{2}\Bigl(\Vert(z_2^m,w_2^m)\Vert_{C^1(\tilde\Omega)}\Vert(u-u^m)\Vert_{L^2(\tilde\Omega)}\\ &\quad+B\Vert\nabla\mathcal A(\mathbf v-\mathbf v^m)\Vert_{L^2(\tilde\Omega)}+B^2\Vert\nabla(\mathbf v-\mathbf v^m)\Vert_{L^2(\tilde\Omega)}\Bigr).
\end{aligned}
\end{equation}
By the pointwise bound on $\mathbf v,\mathbf v^m$, and by \eqref{eq:vvmL2}--\eqref{eq:A2v}, \eqref{eq:D2hdiff}, we bound \eqref{eq:rhsbd} by $B^2(D_2(\mathbf v,\mathbf v^m)+2^{-m(2s-3)}M^2)$. By \eqref{eq:calIbd} and H\"older, we bound $\mathcal I$ by
\[
\Vert w_2-w_2^m\Vert_{L^2(\tilde\Gamma)}\Vert z_2-z_2^m\Vert_{L^2(\tilde\Gamma)}. 
\]
By the conditions $z_2|_\Gamma=0$, $z_2^m|_{\Gamma^m}=0$, by fundamental theorem of calculus and by \eqref{eq:etaL2}, we have
\[
\Vert z_2-z_2^m\Vert_{L^2(\tilde\Gamma)}\lesssim_A2^mB\Vert\eta_\Gamma-\eta_{\Gamma^m}\Vert_{L^2(\Gamma_*)}\lesssim_A2^{-m(s-1)}BM. 
\]  
By the trace estimate and by \eqref{eq:A2v}, \eqref{eq:w2H1diff} we have
\[
\Vert w_2-w_2^m\Vert_{L^2(\tilde\Gamma)}\lesssim_A\Vert w_2-w_2^m\Vert_{L^2(\tilde\Omega)}^\frac{1}{2}\Vert w_2-w_2^m\Vert_{W^{1,2}(\tilde\Omega)}^\frac{1}{2}\lesssim_A2^\frac{m}{2}D_2(\mathbf v,\mathbf v^m)^\frac{1}{2}+2^{-m(s-2)}BM. 
\]
Combining these we obtain
\[
|\mathcal I|\lesssim_A2^{-m(s-\frac{3}{2})}BD_2(\mathbf v,\mathbf v^m)+2^{-m(2s-3)}B^2M^2. 
\]
This completes the proof of the proposition. \proofqed

\begin{corollary}\label{cor:PH2}
With the assumptions above, we have
\[
\Vert(\nabla\divg(u-u^m),\ \nabla^2(h-h^m))\Vert_{L^2(\tilde\Omega)}\lesssim_A2^{-m(s-\frac{3}{2})}M+2^{-m(s-1-\delta)}B^2M\text{ on }[0,T],
\]
provided that $\int_0^TB(t)^2\,dt<1$ and that
\[
D_2(\mathbf v,\mathbf v^m)|_{t=0}\leq2^{-m(2s-3)}M^2. 
\]
\end{corollary}

\section{A priori estimates}\label{sec:apest}

For even integers $k\geq4$, we first construct an energy equivalent to the $\mathbf H^k$ norm and prove its propagation. We then add one acoustic derivative to obtain the $\mathbf H^k_1$ estimate. The estimates separate the geometric, acoustic, and vorticity contributions.

\subsection{Energy functional and coercivity}\label{sub:coer}

For $(h,u,\Omega)\in\mathbf H^k$, define the energy as the sum of vorticity, acoustic, and tangential/free-surface contributions:
\begin{equation}\label{eq:Evwb}
E_k(h,u,\Omega)=E_k^v+E_k^w+E_k^b. 
\end{equation}
The three contributions are
\[
E_k^v=\Vert\curl u\Vert_{W^{k-1,2}(\Omega)}^2,
\]
\[
E_k^w=\sum_{j=0}^k\Vert(z_j,w_j)\Vert_{L^2(\Omega)}^2,
\]
and
\[
E_k^b=\sum_{j=1}^{k/2}\left(\Vert(L_\zeta^jh,L_\zeta^ju)\Vert_{L^2(\Omega)}^2+\Vert a^\frac{1}{2}\Delta_\Gamma^{j-1}\kappa\Vert_{L^2(\Gamma)}^2\right),
\]
where $w_j=(-D_t)^ju$ is the multilinear expression defined via \eqref{eq:CE}, $L_\zeta$ is a second-order differential operator
\[
L_\zeta=\zeta^{ij}\partial_i\partial_j,
\]
and $\zeta^{ij}$ is the harmonic extension of the tangential projection:
\[
\zeta^{ij}=\mathcal H(\delta^{ij}-n^in^j). 
\]

\begin{theorem}[Coercivity of the energy]\label{thm:coer}
Suppose $(h,u,\Omega)\in\mathbf H^k$ with $\Gamma\in\Lambda_*$ satisfies
\[
z_j|_\Gamma=0,\quad j=0,...,k-1,
\] 
the Taylor sign condition, and the bound
\begin{equation}\label{eq:AinLinft}
\Vert(h,u)\Vert_{C^{\frac{1}{2}+\delta}(\Omega)}+\Vert\Gamma\Vert_{C^{1,\delta}}\leq A. 
\end{equation}
Then we have
\[
E_k(h,u,\Omega)\sim_A\Vert(h,u,\Omega)\Vert_{\mathbf H^k}^2. 
\]
\end{theorem}
\proofheading{Proof of Theorem~\ref{thm:coer}, upper bound}First, we have $E_k^v\lesssim_A\Vert u\Vert_{H^k(\Omega)}^2$. 

Next, induction shows that $z_j,w_j$ are multilinear expressions in $(h,u)$ involving $j$ derivatives. Thus, Proposition~\ref{prop:trilest} implies $E_k^w\lesssim_A\Vert(h,u)\Vert_{H^k(\Omega)}^2$. Similarly, with 
\begin{equation}\label{eq:zetabd}
\Vert\zeta\Vert_{C^\delta(\Omega)}\lesssim_A1,\quad\Vert\zeta\Vert_{H^{k-\frac{1}{2}}(\Omega)}\lesssim_A\Vert\Gamma\Vert_{H^k},
\end{equation}
Proposition~\ref{prop:trilest} implies $\sum_{j=1}^{k/2}\Vert(L_\zeta^jh,L_\zeta^ju)\Vert_{L^2(\Omega)}\lesssim_A\Vert(h,u,\Omega)\Vert_{\mathbf H^k}$, while \eqref{eq:zetabd} can be deduced from Propositions \ref{prop:Schauder}, \ref{prop:LapDiri}. 

Finally, with $\kappa=\divg_\Gamma n$, Proposition~\ref{prop:mlinb} implies $\sum_{j=1}^{k/2}\Vert\Delta_\Gamma^{j-1}\kappa\Vert_{L^2(\Gamma)}\lesssim_A\Vert\Gamma\Vert_{H^k}$. This completes the proof of the upper bound. \proofqed

The reverse inequality uses the following Hodge estimate.

\begin{proposition}\label{prop:Lzeta}
Let $s\geq0$ and $k\geq1$ be integers. For $r\geq0$, $\alpha,\beta\in[0,1]$ and partition $v=v_j^1+v_j^2$ we have
\begin{equation}\label{eq:Lzetabd}
\begin{aligned}
\Vert L_\zeta^kv\Vert_{H^s(\Omega)}\lesssim_A\Vert v\Vert_{H^{s+2k}(\Omega)}+\Vert\Gamma\Vert_{H^{s+2k}}\Vert v\Vert_{C^{\frac{1}{2}+\delta}}. 
\end{aligned}
\end{equation}
and
\begin{equation}\label{eq:Lzetacoer}
\begin{aligned}
\Vert v\Vert_{H^{s+2k}(\Omega)}\lesssim_A&\Vert \divg v\Vert_{H^{s+2k-1}(\Omega)}+\Vert \curl v\Vert_{H^{s+2k-1}(\Omega)}+\Vert L_\zeta^kv\Vert_{H^s(\Omega)}\\
&+\Vert\Gamma\Vert_{H^{s+2k-\frac{1}{2}+r}}\sup_j2^{-j(r+\alpha-1)}\Vert v_j^1\Vert_{C^\alpha(\Omega)}+\sup_j2^{j(s+2k-1+\beta-\frac{\delta}{2})}\Vert v_j^2\Vert_{H^{1-\beta}(\Omega)}. 
\end{aligned}
\end{equation}
\end{proposition}
\proofheading{Proof}By Proposition~\ref{prop:LapDiri} we have
\begin{equation}\label{eq:zetaregularity}
\Vert\zeta\Vert_{H^{s+2k+r-\frac{1}{2}}(\Omega)}\lesssim_A\Vert\Gamma\Vert_{H^{s+2k+r}}. 
\end{equation}
Then \eqref{eq:Lzetabd} follows from Proposition~\ref{prop:trilest}. For \eqref{eq:Lzetacoer} we perform induction on $k,s$. For $k=1,s=0$, note the pointwise bound
\begin{equation}\label{eq:ptHodge}
|\nabla^mv|^2\lesssim|\nabla^{m-1}\divg v|^2+|\nabla^{m-1}\curl v|^2+\sum_{i_1,\dots,i_m}\zeta^{i_1j_1}\cdots\zeta^{i_mj_m}\partial^m_{i_1\cdots i_m}v\partial^m_{j_1\cdots j_m}v,
\end{equation}
proved in \cite{CL}. Then it suffices to control $\int\zeta^{ij}\zeta^{kl}\partial_i\partial_kv\cdot\partial_j\partial_lv\,dx$. For this, integrating by parts twice we derive
\[\begin{aligned}
\Vert L_\zeta v\Vert_{L^2(\Omega)}^2=&\int_\Omega\zeta^{ij}\zeta^{kl}\partial_i\partial_jv\cdot\partial_k\partial_lv\,dx\\
=&\int_\Omega\zeta^{ij}\zeta^{kl}\partial_i\partial_kv\cdot\partial_j\partial_lv\,dx+\int_\Omega\partial_j(\zeta^{ij}\zeta^{kl})\partial_k\partial_iv\cdot\partial_lv-\partial_k(\zeta^{ij}\zeta^{kl})\partial_j\partial_iv\cdot\partial_lv\,dx. 
\end{aligned}\]
The last integral can be controlled by using \eqref{eq:zetaregularity}, Proposition~\ref{prop:bilest2} and Young's inequality, which leads to \eqref{eq:Lzetacoer} with $k=1,s=0$. For the case $k=1,s>0$, we commute $L_\zeta$ with $\partial^s$. The commutator $[\partial^s,L_\zeta]v$ can be controlled by using \eqref{eq:zetaregularity}, Proposition~\ref{prop:bilest2}. Denoting $w=\partial^sv$, we have
\[\begin{aligned}
\Vert w\Vert_{H^2(\Omega)}\lesssim_A&\Vert\divg w\Vert_{H^1(\Omega)}+\Vert\curl w\Vert_{H^1(\Omega)}+\Vert L_\zeta w\Vert_{L^2(\Omega)}\\
&+\Vert\Gamma\Vert_{H^{\frac{3}{2}+r+s}}\sup_j2^{-j(r+s+\alpha-1)}\Vert w_j^1\Vert_{C^\alpha(\Omega)}+\sup_j2^{-j(1+\beta-\frac{\delta}{2})}\Vert w_j^2\Vert_{H^{1-\beta}(\Omega)}, 
\end{aligned}\]
where $w=w_j^1+w_j^2$ is specified as
\[
w_j^1=P_{<j}v_j^1,\quad w_j^2=P_{<j}v_j^2+P_{\geq j}v. 
\]
By Bernstein's estimates we obtain \eqref{eq:Lzetacoer} for $k=1,s>0$. Then the general case follows by induction on $k$. \proofqed
\proofheading{Proof of Theorem~\ref{thm:coer}, lower bound}First, we recover the regularity of $\Gamma$. Consider the decomposition of mean curvature: 
\[
\kappa=\kappa_j^1+\kappa_j^2=\divg_\Gamma n_j^1+\divg_\Gamma n_j^2,\quad\text{where }n_j^1(Y)=P_{<j}\,n(Y),\  n_j^2(Y)=P_{\geq j}\,n(Y),
\]
$Y=\text{Id}+\eta n_*:\Gamma_*\to\Gamma$ and $P_{<j}$ is the frequency localization on $\Gamma_*$. Bernstein's estimate implies 
\[
\Vert\kappa_j^1\Vert_{C^0(\Gamma)}\lesssim_A2^j,\quad\Vert\kappa_j^2\Vert_{L^2(\Gamma)}\lesssim_A2^{-j(k-2-\frac{\delta}{2})}\Vert\Gamma\Vert_{H^{k-\frac{\delta}{2}}},\quad\forall j>0. 
\]
Then by \eqref{eq:LGk} we have 
\[
\Vert\kappa\Vert_{H^{k-2}(\Gamma)}\lesssim_A\left(E_k^b\right)^{1/2}+\Vert\Gamma\Vert_{H^{k-\frac{\delta}{2}}},
\]
and by Proposition~\ref{prop:kappa} and interpolation we conclude
\[
\Vert\Gamma\Vert_{H^k}\lesssim_A\left(E_k^b\right)^{1/2}. 
\]

Next, recover the regularity of $h$. The expression $z_k$ contains $k$ derivatives of $(h,u)$ and has leading term $\Delta^{k/2}h$. Its nonlinear part satisfies
\[
\Vert N_k\Vert_{L^2(\Omega)}\lesssim_A\Vert(h,u)\Vert_{H^{k-\frac{1}{2}}(\Omega)},
\]
thanks to Proposition~\ref{prop:trilest}. It follows that
\begin{equation}\label{eq:Ltoph}
\Vert\Delta^\frac{k}{2}h\Vert_{L^2(\Omega)}\lesssim_A(E^w_k)^{1/2}+\Vert(h,u)\Vert_{H^{k-\frac{1}{2}}(\Omega)}. 
\end{equation}
The compatibility conditions control the boundary traces of the lower powers of the Laplacian:
\begin{equation}\label{eq:Llowtr}
\Delta^{\frac{k}{2}-l}h|_\Gamma=-N_{k-2l},\quad l=1,...,\frac{k}{2}. 
\end{equation}
Propositions \ref{prop:trilest}, \ref{prop:tr} imply 
\begin{equation}\label{eq:Ntr}
\Vert N_{k-2l}\Vert_{H^{2l-\frac{1}{2}}(\Gamma)}\lesssim_A\Vert(h,u)\Vert_{H^{k-\frac{1}{2}}(\Omega)}. 
\end{equation}
Combining \eqref{eq:Ltoph}--\eqref{eq:Ntr} and \eqref{eq:LDk}, we conclude
\[
\Vert h\Vert_{H^k(\Omega)}\lesssim_A(E^w_k)^{1/2}+\Vert(h,u)\Vert_{H^{k-\frac{1}{2}}(\Omega)}.
\]

Finally, we recover the regularity of $u$. The estimate of $\divg u$ is very similar to that of $h$. By \eqref{eq:LDk}, 
\[
\Delta^j\divg u|_\Gamma=-N_{2j+1},\quad 2j+1<k, 
\]
and multilinear estimate, we obtain 
\[
\Vert\divg u\Vert_{H^{k-1}(\Omega)}\lesssim_A(E^w_k)^{1/2}+\Vert(h,u)\Vert_{H^{k-\frac{1}{2}}(\Omega)}. 
\]
The estimate of $\curl u$ directly follows from Proposition~\ref{prop:SteinExt}: $\Vert\curl u\Vert_{H^{k-1}(\Omega)}\lesssim_A(E^v_k)^{1/2}$. Then \eqref{eq:Lzetacoer} implies
\[
\Vert u\Vert_{H^k(\Omega)}\lesssim_AE_k(h,u,\Omega)^{1/2}+\Vert(h,u)\Vert_{H^{k-\frac{1}{2}}(\Omega)}.
\]
Combining the estimates and using interpolation, we complete the proof of Theorem~\ref{thm:coer}. \proofqed

\subsection{Propagation of the energy}\label{sub:propa}

\begin{theorem}[Propagation of the energy]\label{thm:apest}
Suppose $(h,u,\Omega_t)$ with $\Gamma_t\in\Lambda_*$ is a sufficiently regular solution to \eqref{eq:CE} satisfying the Taylor condition. Suppose \eqref{eq:AinLinft} holds on $[0,T]$ and 
\[
\Vert\nabla u\Vert_{L^\infty(\Omega)}+\Vert\nabla\zeta\Vert_{L^\infty(\Omega)}+\Vert\nabla_n\divg u\Vert_{L^\infty(\Gamma)}+1\leq B(t). 
\]
Then we have 
\[
\frac{d}{dt}E_k(h,u,\Omega_t)\lesssim_AB^2E_k(h,u,\Omega_t). 
\]
\end{theorem}

We use two algebraic identities and a boundary divergence estimate. Lemma~\ref{lem:zwsys} follows by induction from \eqref{eq:CE}.

\begin{lemma}\label{lem:zwsys}
For $j=1,2,...$ we have
\begin{equation}\label{eq:wave}
\left\{\begin{aligned}
&D_tz_j+\divg w_j=f_j\\
&D_tw_j+\nabla z_j=g_j
\end{aligned}\right.
\end{equation}
where $f_j,g_j$ are multilinear expressions without linear terms involving $j+1$ derivatives on $h,u$. 
\end{lemma}

\begin{lemma}\label{lem:Ebsys}
For $j=1,2,...$ we have
\begin{equation}\label{eq:Ebsys}
\left\{
\begin{aligned}
&D_tL_{\zeta}^jh+\divg L_{\zeta}^ju=f_{2j}^b,\quad&\text{in }\Omega,\\
&D_tL_{\zeta}^ju+\nabla L_{\zeta}^jh=g_{2j}^b,\quad&\text{in }\Omega,\\
&D_t\Delta_{\Gamma}^{j-1}\kappa+L_{\zeta}^ju\cdot n=q_{2j}+r^u_{2j},\quad&\text{on }\Gamma,\\
&L_\zeta^jh+a\Delta_{\Gamma}^{j-1}\kappa=r^h_{2j},\quad&\text{on }\Gamma. 
\end{aligned}
\right. 
\end{equation}
Here $f_{2j}^b,g_{2j}^b$ are polynomials in $\partial^m(h,u)$, $\partial^{m-1}\zeta,\partial^{m-2}D_t\zeta$, $1\leq m\leq2j$ involving $2j+1$ spacetime derivatives; $q_{2j}$ is a polynomial in $\nabla_\Gamma^m\partial u,\nabla_\Gamma^ln$, $m\leq2j-2$, $l\leq2j-1$ involving $2j-1$ covariant derivatives; and $r^f_{2j}$ is a polynomial in $\partial^mf,\nabla_\Gamma^ln,\partial^l\zeta$, $1\leq m\leq2j-1$, $0\leq l\leq2j-2$ involving $2j$ derivatives. 
\end{lemma}
\proofheading{Proof}Apply $L_\zeta^j$ to \eqref{eq:CE}, commute derivatives, and proceed by induction to obtain the first two equations.

To prove the third equation, we start from the basic relation
\begin{equation}\label{eq:Dtn}
D_tn=-\nabla_\Gamma u\cdot n. 
\end{equation}
Perform $\divg_\Gamma$ on both sides and use the formula
\[
[D_t,\zeta^{ij}\partial_j]f|_\Gamma=-\zeta^{ij}\partial_ju\cdot\nabla_\Gamma f+\nabla_\Gamma f\cdot\nabla_\Gamma u^kn_kn_i,
\]
we then derive
\[
D_t\kappa+\Delta_\Gamma u\cdot n=q_2. 
\]
Further performing $\Delta_{\Gamma}^{j-1}$ yields
\[
D_t\Delta_{\Gamma}^{j-1}\kappa+\Delta_\Gamma^ju\cdot n=q_{2j}. 
\]
It remains to show that 
\begin{equation}\label{eq:LzLGj}
L_\zeta^ju=\Delta_\Gamma^ju+q_{2j}+r^u_{2j}. 
\end{equation}
This follows from the formula
\begin{equation}\label{eq:Lztr}
L_\zeta f|_\Gamma=\Delta_\Gamma f+\kappa\nabla_nf
\end{equation}
and induction. 

The proof of (\ref{eq:Ebsys}d) is very similar. Indeed, $h|_\Gamma=0$ and \eqref{eq:Lztr} imply
\begin{equation}\label{eq:Lzhtr}
L_\zeta h|_\Gamma=-a\kappa. 
\end{equation}
As in \eqref{eq:LzLGj}, we have
\[
L_\zeta^jh|_\Gamma=-\Delta_\Gamma^{j-1}(a\kappa)+r_{2j}^h. 
\]
Finally, distributing derivatives in $\Delta_\Gamma^{j-1}(a\kappa)$ we conclude (\ref{eq:Ebsys}d). This completes the proof of the lemma. \proofqed

\begin{lemma}\label{lem:-1/2}
\[
\Vert\divg_\Gamma v\Vert_{H^{-\frac{1}{2}}(\Gamma)}\lesssim_A\Vert v\Vert_{H^\frac{1}{2}(\Gamma)}. 
\]
\end{lemma}
\proofheading{Proof}Let $g\in H^\frac{1}{2}(\Gamma)$. We have
\[
\int_\Gamma\divg_\Gamma v\,g\,dS=\sum_j\left(\int_\Gamma\divg_\Gamma P_{<j}v\,P_jg\,dS+\int_\Gamma\divg_\Gamma P_{\geq j}v\,P_jg\,dS\right):=\sum_j(I_j^1+I_j^2). 
\]
Bernstein's estimate implies
\[
|I_j^1|\lesssim_A\sum_{i<j}2^{-\frac{1}{2}(j-i)}\Vert P_iv\Vert_{H^\frac{1}{2}(\Gamma)}\Vert P_jg\Vert_{H^\frac{1}{2}(\Gamma)}. 
\]
Stokes' theorem gives
\[
I_j^2=\sum_{i\geq j}\int_\Gamma P_iv\cdot\nabla_\Gamma P_jg\,dS. 
\]
Again by Bernstein, we derive
\[
|I_j^2|\lesssim_A\sum_{i\geq j}2^{-\frac{1}{2}(i-j)}\Vert P_iv\Vert_{H^\frac{1}{2}(\Gamma)}\Vert P_jg\Vert_{H^\frac{1}{2}(\Gamma)}. 
\]
By the convolution inequality we deduce $|\int_\Gamma\divg_\Gamma v\,g\,dS|\lesssim_A\Vert v\Vert_{H^\frac{1}{2}(\Gamma)}\Vert g\Vert_{H^\frac{1}{2}(\Gamma)}$. This completes the proof. \proofqed
\proofheading{Proof of Theorem~\ref{thm:apest}}First, we control $\frac{d}{dt}E^v_k$. We note that $D_t\curl u$ is a quadratic form in $\nabla u$. Thus, Proposition~\ref{prop:trilest} gives
\[
\frac{d}{dt}E^v_k\lesssim_AB\Vert u\Vert_{H^k(\Omega)}^2. 
\]

Next, we control $\frac{d}{dt}E^w_k$ by using Lemma~\ref{lem:zwsys}. Since $(h,u,\Omega_t)$ is a sufficiently regular solution, from the boundary condition \eqref{eq:vacbc} we derive
\[
z_j|_\Gamma=0,\quad j=0,...,k. 
\]
Invoking Proposition~\ref{prop:linEid} with $s=0$ and using Proposition~\ref{prop:trilest} for handling $f_j,g_j$, we obtain
\[
\frac{d}{dt}E^w_k\lesssim_AB\Vert(h,u)\Vert_{H^k(\Omega)}^2.
\]

Finally, apply Lemma~\ref{lem:Ebsys} and Proposition~\ref{prop:linEid} to the tangential energy. We need the following bounds:
\begin{equation}\label{eq:fg}
\sum_{j=1}^{k/2}\Vert(f_{2j}^b,g_{2j}^b)\Vert_{L^2(\Omega)}\lesssim_AB^2\Vert(h,u,\Omega)\Vert_{\mathbf H^k},
\end{equation}
\begin{equation}\label{eq:q}
\sum_{j=1}^{k/2}\Vert q_{2j}\Vert_{L^2(\Gamma)}\lesssim_AB\Vert(h,u,\Omega)\Vert_{\mathbf H^k},
\end{equation}
\begin{equation}\label{eq:r}
\sum_{j=1}^{k/2}\Vert(r_{2j}^h,r_{2j}^u)\Vert_{H^\frac{1}{2}(\Gamma)}\lesssim_AB\Vert(h,u,\Omega)\Vert_{\mathbf H^k},
\end{equation}
and
\begin{equation}\label{eq:wn}
\sum_{j=1}^{k/2}\Vert L_\zeta^ju\cdot n\Vert_{H^{-\frac{1}{2}}(\Gamma)}\lesssim_AB\Vert(h,u,\Omega)\Vert_{\mathbf H^k}. 
\end{equation}

\proofheading{Proof of \eqref{eq:fg}}We need to control $D_t\zeta=\mathcal HD_t\zeta+[D_t,\mathcal H]\zeta$. For the harmonic part, by \eqref{eq:Dtn}, \eqref{eq:calHCalp}, \eqref{eq:LapDiri} we have
\begin{equation}\label{eq:Dtzetahar}
\Vert\mathcal HD_t\zeta\Vert_{L^\infty(\Omega)}\lesssim B,\quad \Vert\mathcal HD_t\zeta\Vert_{H^{k-1}(\Omega)}\lesssim_A\Vert(h,u,\Omega)\Vert_{\mathbf H^k}. 
\end{equation}
For the commutator, we use the formula: 
\begin{equation}\label{eq:[Dt,H]}
[D_t,\mathcal H]f=\Delta^{-1}(\Delta u\cdot\nabla\mathcal Hf+2\nabla u:\nabla^2\mathcal Hf)=\Delta^{-1}\divg F(\nabla u,\nabla\mathcal Hf) 
\end{equation}
where $F$ is some vector-valued bilinear form. The second equality can be seen from $\nabla u:\nabla^2\mathcal Hf=\divg(\nabla u\cdot\nabla\mathcal Hf)-\Delta u\cdot\nabla\mathcal Hf$ and
\[\begin{aligned}
\Delta u_j\partial_j\mathcal Hf&=\partial_i(\partial_iu_j-\partial_ju_i)\partial_j\mathcal Hf+\nabla\divg u\cdot\nabla\mathcal Hf\\
&=\partial_i[(\partial_iu_j-\partial_ju_i)\partial_j\mathcal Hf]+\divg(\divg u\nabla\mathcal Hf).
\end{aligned}\]
Thus, by $L^p$ estimate for elliptic system in divergence form and by Sobolev embedding, we derive
\[
\Vert[D_t,\mathcal H]\zeta\Vert_{L^\infty(\Omega)}\lesssim_AB^2. 
\]
By \eqref{eq:LapDiri} we derive
\begin{align}
\Vert[D_t,\mathcal H]\zeta\Vert_{H^{k-\frac{1}{2}}(\Omega)}&\lesssim_AB^2\Vert\Gamma\Vert_{H^k}+\Vert F(\nabla u,\nabla\zeta)\Vert_{H^{k-\frac{3}{2}}(\Omega)}\nonumber\\
&\lesssim_AB^2\Vert\Gamma\Vert_{H^k}+B(\Vert u\Vert_{H^{k-\frac{1}{2}}(\Omega)}+\Vert\Gamma\Vert_{H^k}).\label{eq:Dtzetacomm}
\end{align}
Combining \eqref{eq:zetabd}, \eqref{eq:Dtzetahar}, \eqref{eq:Dtzetacomm} and Proposition~\ref{prop:trilest}, we justify \eqref{eq:fg}. 

\proofheading{Proof of \eqref{eq:q}}Note that 
\begin{equation}\label{eq:Dnzeta}
\partial_i^\top n_j=-\zeta^{ik}\partial_k\zeta^{jl}n_l. 
\end{equation}
Thus, we have $\Vert\nabla_\Gamma n\Vert_{L^\infty(\Gamma)}\leq B(t)$. Then Propositions \ref{prop:mlinb}, \ref{prop:tr} implies
\[
\sum_{j=1}^{k/2}\Vert q_{2j}\Vert_{L^2(\Gamma)}\lesssim_AB(\Vert\Gamma\Vert_{H^k}+\Vert\nabla u\Vert_{H^{k-2}(\Gamma)})\lesssim_AB\Vert(h,u,\Omega)\Vert_{\mathbf H^k}. 
\]

\proofheading{Proof of \eqref{eq:r}}By \eqref{eq:Dnzeta}, $r_{2j}^{h,u}$ is a polynomial in $\partial^m(h,u)$ and $\partial^l\zeta$ and hence extends to $\Omega$. The trace estimate reduces the claim to the corresponding $H^1(\Omega)$ bounds, which follow from Proposition~\ref{prop:trilest}.

\proofheading{Proof of \eqref{eq:wn}}As shown in the proof of Lemma~\ref{lem:Ebsys}, we have $L_\zeta^ju=\Delta_\Gamma^ju+q_{2j}+r_{2j}^u$. Thus, it remains to control $\Delta_\Gamma^ju\cdot n$. Write $\Delta_\Gamma^ju\cdot n=\divg_\Gamma(\nabla_\Gamma\Delta_\Gamma^{j-1}u\cdot n)-\nabla_\Gamma\Delta_\Gamma^{j-1}u:\nabla_\Gamma n$. By Propositions \ref{prop:mlinb}, \ref{prop:tr} we have
\[
\sum_{j=1}^{k/2}\left(\Vert\nabla_\Gamma\Delta_\Gamma^{j-1}u\cdot n\Vert_{H^\frac{1}{2}(\Gamma)}+\Vert\nabla_\Gamma\Delta_\Gamma^{j-1}u:\nabla_\Gamma n\Vert_{L^2(\Gamma)}\right)\lesssim_AB(\Vert u\Vert_{H^k(\Omega)}+\Vert\Gamma\Vert_{H^k}). 
\]
Then the conclusion is clear from Lemma~\ref{lem:-1/2}. 

Now Proposition~\ref{prop:linEid} implies 
\[
\frac{d}{dt}E^b_k\lesssim_AB^2\Vert(h,u,\Omega)\Vert_{\mathbf H^k}^2. 
\]
Putting together the estimates of the temporal derivative of $E_k^{v,w,b}$ and using Theorem~\ref{thm:coer}, we complete the proof of Theorem~\ref{thm:apest}. \proofqed

\subsection{Energy estimate in \texorpdfstring{$\mathbf H^k_1$}{the hybrid space Hk1}}

To measure one additional acoustic derivative, define

\begin{equation}\label{eq:Ek1}
E_{k,1}(h,u,\Omega)=E^w_{k+1}+E^v_k+E^b_k=E_k(h,u,\Omega)+\Vert(z_{k+1},w_{k+1})\Vert_{L^2(\Omega)}^2. 
\end{equation}
We shall prove the following energy estimate. 

\begin{theorem}[One additional acoustic derivative]\label{thm:apinHk1}
With the assumptions in Theorem~\ref{thm:apest}, set
\begin{equation}\label{eq:B1}
B_1(t)=B(t)+\Vert(\Delta h,\nabla\divg u)\Vert_{L^\infty(\Omega_t)}. 
\end{equation}
Then we have
\begin{equation}\label{eq:coerEk1}
E_{k,1}(h,u,\Omega)\sim_A\Vert(h,u,\Omega)\Vert_{\mathbf H_1^k}^2+O\left(B_1^2E_k(h,u,\Omega)\right),
\end{equation}
and
\begin{equation}\label{eq:apEk1}
\frac{d}{dt}E_{k,1}(h,u,\Omega)\lesssim_AB_1^2E_{k,1}(h,u,\Omega). 
\end{equation}
\end{theorem}

The new term is the highest-order acoustic energy. Its estimate relies on the following structure of the bilinear terms.

\begin{lemma}\label{lem:bilchar}
A bilinear term in $f_j,g_j$, $j\geq4$ contains $\Delta h,\nabla\divg u$ as a factor. In other words, it looks like $\partial^{j-1-l}(h,u)\,\partial^l(\Delta h,\nabla\divg u)$. 
\end{lemma}
\proofheading{Proof}We only need to verify the conclusion for bilinear terms in $z_{\geq4},w_{\geq4}$. Iteratively perform $D_t$ four times on $(h,u)$ and use \eqref{eq:CE}. Since we are concerned with bilinear terms, exactly one $D_t$ contributes to the commutator with spatial derivatives, and the others act like the linear generator $\mathcal A: (z,w)\mapsto(-\divg w,-\nabla z)$. Thus, at least one factor is a spatial derivative of $\mathcal A^2(h,u)=(\Delta h,\nabla\divg u)$. \proofqed
\proofheading{Proof of Theorem~\ref{thm:apinHk1}}We first prove \eqref{eq:coerEk1}. If $\mathcal N_{k+1}$ is a nonlinear term in $z_{k+1},w_{k+1}$, Proposition~\ref{prop:trilest} implies 
\[
\Vert\mathcal N_{k+1}\Vert_{L^2(\Omega)}\lesssim_AB_1\Vert(h,u)\Vert_{H^k(\Omega)}. 
\]
Thus, we have $E_{k,1}(h,u,\Omega)\lesssim_A\Vert(h,u,\Omega)\Vert_{\mathbf H_1^k}^2+O\left(B_1^2E_k(h,u,\Omega)\right)$ and
\begin{equation}\label{eq:Ltophu}
\Vert(\Delta^\frac{k}{2}\divg u,\nabla\Delta^\frac{k}{2}h)\Vert_{L^2(\Omega)}\lesssim_A(E^w_{k+1})^\frac{1}{2}+B_1\Vert(h,u,\Omega)\Vert_{\mathbf H^k}. 
\end{equation}
Similarly, if $\mathcal N_{k+1-2l},\mathcal N_{k-2l}$ are respectively nonlinear terms in $w_{k+1-2l},z_{k-2l}$, Propositions \ref{prop:tr}, \ref{prop:trilest} imply
\[\begin{aligned}
\Vert\mathcal N_{k+1-2l}\Vert_{H^{2l-\frac{1}{2}}(\Gamma)}+\Vert\mathcal N_{k-2l}\Vert_{H^{2l+\frac{1}{2}}(\Gamma)}&\lesssim_AB_1\Vert\Gamma\Vert_{H^k(\Omega)}+\Vert\mathcal N_{k+1-2l}\Vert_{H^{2l}(\Omega)}+\Vert\mathcal N_{k-2l}\Vert_{H^{2l+1}(\Omega)}\\
&\lesssim_AB_1\Vert(h,u,\Omega)\Vert_{\mathbf H^k}. 
\end{aligned}\]
Thus, by the boundary condition $z_j|_\Gamma=0$, we have
\begin{equation}\label{eq:Linthu}
\sum_{l=1}^{k/2}\Vert\Delta^{\frac{k}{2}-l}\divg u\Vert_{H^{2l-\frac{1}{2}}(\Gamma)}+\sum_{l<k/2}\Vert\Delta^{\frac{k}{2}-l}h\Vert_{H^{2l+\frac{1}{2}}(\Gamma)}\lesssim_AB_1\Vert(h,u,\Omega)\Vert_{\mathbf H^k}. 
\end{equation}
Combining \eqref{eq:Ltophu}, \eqref{eq:Linthu}, by \eqref{eq:LDk} we obtain
\[
\Vert(\divg u,\Delta h)\Vert_{H^k\times H^{k-1}(\Omega)}\lesssim_A(E^w_{k+1})^\frac{1}{2}+B_1\Vert(h,u,\Omega)\Vert_{\mathbf H^k}. 
\]
This justifies \eqref{eq:coerEk1}. 

We next prove \eqref{eq:apEk1}. If $\mathcal T_{k+2}$ in $f_{k+2},g_{k+2}$ contains at least cubic nonlinearity, by Proposition~\ref{prop:trilest} we have
\[
\Vert\mathcal T_{k+2}\Vert_{L^2(\Omega)}\lesssim_AB_1^2\Vert(h,u)\Vert_{H^k(\Omega)}. 
\]
If $\mathcal B_{k+2}$ is a bilinear term in $f_{k+2},g_{k+2}$, by Lemma~\ref{lem:bilchar} and Proposition~\ref{prop:trilest} we have
\[
\Vert\mathcal B_{k+2}\Vert_{L^2(\Omega)}\lesssim_AB_1(\Vert(h,u)\Vert_{H^k(\Omega)}+\Vert(\divg u,\Delta h)\Vert_{H^k\times H^{k-1}(\Omega)}). 
\]
These combined with \eqref{eq:coerEk1} and Theorem~\ref{thm:coer} yield
\[
\frac{d}{dt}E_{k+1}^w\lesssim_AB_1^2E_{k,1}(h,u,\Omega). 
\]
This completes the proof of Theorem~\ref{thm:apinHk1}. \proofqed

\section{Construction of regular solutions}\label{sec:regsol}

We construct regular solutions by a discrete evolution that alternates smoothing, compatibility correction, and transport. Fix an even integer $k=2k_0>4$. To control the regularization errors, augment the energy by two weighted vorticity derivatives:
\begin{equation}\label{eq:calE}
\mathcal E(h,u,\Omega)=E_k(h,u,\Omega)+\epsilon^4\Vert\curl u\Vert_{W^{k+1,2}(\Omega)}^2.
\end{equation}

\begin{theorem}[One step of the regularized Euler scheme]\label{thm:1step}
Let $M,P,\epsilon$ be positive constants with
\begin{equation}\label{eq:numrel}
M\ll P\ll\frac{1}{\epsilon}. 
\end{equation}
Suppose $(h_0,u_0,\Omega_0)$ satisfies the following assumptions.

(A0) Taylor sign and geometry: \eqref{eq:TS} holds, $\Gamma_0\in\Lambda_*$, and
\begin{equation}\label{eq:ptass0}
\Vert\Gamma_0\Vert_{C^{1,\delta}}+\Vert(h_0,u_0)\Vert_{C^{\frac{1}{2}+\delta}(\Omega_0)}\leq A_0. 
\end{equation}

(A1) Energy bound: 
\begin{equation}\label{eq:Ebd0}
\mathcal E(h_0,u_0,\Omega_0)\leq M^2.
\end{equation}
(A2) Approximate compatibility: $h_0|_{\Gamma_0}=0$ and 
\begin{equation}\label{eq:zj0tr0}
\Vert z_{j,0}|_{\Gamma_0}\Vert_{H^{k-j-r}(\Gamma_0)}\leq P\epsilon^{r+\frac{1}{2}},\quad\frac{1}{2}\leq r\leq\min\{k-j,\frac{3}{2}\},\ j=1,...,k-1. 
\end{equation}
(A3) Initial regularity: 
\begin{equation}\label{eq:reg0}
\Vert(h_0,u_0,\Omega_0)\Vert_{\mathbf H^{k+1}}\leq P\epsilon^{-1}. 
\end{equation}
Then there is a $(h_1,u_1,\Omega_1)$ satisfying: \\
(i) One-step energy bound:
\begin{equation}\label{eq:Emono}
\mathcal E(h_1,u_1,\Omega_1)\leq\mathcal E(h_0,u_0,\Omega_0)+C(P)\epsilon.
\end{equation}
(ii) Compatibility and regularity: $(h_1,u_1)$ satisfies (A2) on $\Gamma_1$ and $(h_1,u_1,\Omega_1)$ satisfies (A3). \\
(iii) Approximate solution: 
\begin{equation}\label{eq:1stepappr}
\left\{\begin{aligned}
&h_1-h_0+\epsilon(u_0\cdot\nabla h_0+\divg u_0)=O_{C^1}(\epsilon^2),&\text{in }\Omega_1\cap\Omega_0\\
&u_1-u_0+\epsilon(u_0\cdot\nabla u_0+\nabla h_0)=O_{C^1}(\epsilon^2),&\text{in }\Omega_1\cap\Omega_0\\
&\Omega_1=(\text{Id}+\epsilon u_0)\Omega_0+O_{C^1}(\epsilon^2). 
\end{aligned}\right.
\end{equation}
\end{theorem}
\par\smallskip\noindent\textbf{Remark.}\ 
The estimates \eqref{eq:toterrinreg}, \eqref{eq:x1def}, and \eqref{eq:expsch}, together with Sobolev embedding, show that the surface changes by $O_{C^{1,\delta}}(\epsilon)$ and that the corresponding pointwise bound changes by $O(\epsilon)$. They also control the change in the Taylor coefficient. These bounds are tracked in the iteration below.
%In compensation, we will check on the pointwise bounds in subsection \ref{sub:cvgitersch} when we construct exact solutions. 

The proof of Theorem~\ref{thm:1step} occupies the next five subsections. Sections~\ref{sec:regofdom}--\ref{sub:tangreg} regularize the domain and the fluid variables using elliptic resolvents. Section~\ref{sub:compat_recovery} corrects the boundary traces using extension operators with prescribed normal derivatives. Section~\ref{sub:transport_step} analyzes transport and proves the one-step estimate. Section~\ref{sub:cvgitersch} then passes to the limit in the discrete scheme.

\subsection{Regularization of the domain}\label{sec:regofdom}

The regularized surface $\Gamma$ is defined via
\begin{equation}\label{eq:etaepsdef}
\eta=(I-\epsilon^2\Delta_{\Gamma_*})^{-1}\eta_0-C_0\epsilon^2,
\end{equation}
where $C_0$ is a constant depending only on $M$ such that $\Gamma\subset\Omega_0$. By spectral resolution we have
\begin{equation}\label{eq:etaregerr}
\Vert\Gamma\Vert_{H^{k+r}}\lesssim_M\epsilon^{-r},\quad\Vert\eta-\eta_0\Vert_{H^{k-r}(\Gamma_*)}\lesssim_M\epsilon^r,\quad\Vert\Gamma\Vert_{H^{k+3}}\lesssim_MP\epsilon^{-3},\quad r\in[0,2]. 
\end{equation}

\begin{lemma}\label{lem:zetaerr}
If $\zeta^{ij}=\delta^{ij}-\mathcal H(n^in^j)$ and $\zeta_0^{ij}=\delta^{ij}-\mathcal H_0(n_0^in_0^j)$, then we have
\begin{equation}\label{eq:zetaerr}
\Vert\zeta\Vert_{H^{k-\frac{1}{2}+r}(\Omega)}\lesssim_M\epsilon^{-r},\quad\Vert\zeta-\zeta_0\Vert_{H^{k-\frac{1}{2}-r}(\Omega)}\lesssim_M\epsilon^r,\quad\Vert\zeta\Vert_{H^{k+\frac{5}{2}}(\Omega)}\lesssim_MP\epsilon^{-3}\quad r\in[0,2]. 
\end{equation}
\end{lemma}
\proofheading{Proof}The regularity bound for $\zeta$ follows from the bound for $\eta$. To estimate the error, define a homotopy $Y_t:\Omega_*\to\mathbb R^3$ by
\begin{equation}\label{eq:Ytdef}
\Delta Y_t=0\quad\text{in }\Omega_*,\qquad Y_t|_{\Gamma_*}=\text{Id}+(t\eta+(1-t)\eta_0)n_*. 
\end{equation}
Then $Y_t$ is a diffeomorphism onto its image with $\Omega_0=Y_0(\Omega_*)\supset Y_t(\Omega_*)\supset Y_1(\Omega_*)=\Omega$. Using the relation 
\begin{equation}\label{eq:chainrule}
\partial_tv(Y_t)=\nabla v\circ Y_t\cdot\mathcal H_*[(\eta-\eta_0)n_*]
\end{equation}
and the bilinear estimate, we derive $\Vert\partial_t\zeta_0^{ij}(Y_t)\Vert_{H^{k-5/2}(\Omega_*)}\lesssim_M\epsilon^2$. Then by the trace estimate we have
\[
\Vert\zeta_0^{ij}(Y_1)-\zeta_0^{ij}(Y_0)\Vert_{H^{k-3}(\Gamma_*)}\lesssim_M\epsilon^2. 
\]
On the other hand, expressing $\zeta^{ij}(Y_1),\zeta_0^{ij}(Y_0)$ in terms of $\eta,\eta_0$, we find that
\[
\Vert\zeta^{ij}(Y_1)-\zeta_0^{ij}(Y_0)\Vert_{H^{k-3}(\Gamma_*)}\lesssim_M\epsilon^2. 
\]
Combining these gives
\[
\Vert\zeta^{ij}-\zeta_0^{ij}|_{\Gamma}\Vert_{H^{k-3}(\Gamma)}\lesssim_M\Vert\zeta^{ij}(Y_1)-\zeta_0^{ij}(Y_1)\Vert_{H^{k-3}(\Gamma_*)}\lesssim_M\epsilon^2.
\]
Interpolating this with the energy bound, we get $\Vert\zeta^{ij}-\zeta_0^{ij}|_{\Gamma}\Vert_{H^{k-1-r}(\Gamma)}\lesssim_M\epsilon^r$, $r\in[0,2]$. Finally, since $\zeta,\zeta_0$ are both harmonic in $\Omega$, by the elliptic estimate we conclude the lemma. \proofqed

The change from $\Gamma_0$ to $\Gamma$ alters the compatibility traces. Using \eqref{eq:chainrule}, the bilinear estimate, and the trace estimate gives
\begin{equation}\label{eq:zj0tr}
\Vert z_{j,0}|_\Gamma\Vert_{H^{k-r-j}(\Gamma)}\lesssim_MP\epsilon^{r+\frac{1}{2}},\quad \frac{1}{2}\leq r\leq\min\{k-j,\frac{3}{2}\},\ j=1,...,k-1. 
\end{equation}
For $j=0$, this argument gives an estimate with a loss of half a derivative:
\begin{equation}\label{eq:h0tr}
\Vert h_0|_\Gamma\Vert_{H^{k-2}(\Gamma)}\lesssim_M\epsilon^2,\quad\Vert h_0|_\Gamma\Vert_{H^{k-\frac{1}{2}}(\Gamma)}\lesssim_M\epsilon+\Vert\eta-\eta_0\Vert_{H^{k-\frac{1}{2}}(\Gamma_*)}\lesssim_M\epsilon+\Vert\Gamma\Vert_{H^{k+\frac{3}{2}}}\epsilon^2. 
\end{equation}

\begin{lemma}\label{lem:Emono0}
We have
\[
\mathcal E(h_0|_\Omega,u_0|_\Omega,\Omega)+\frac{\epsilon^2}{C(M)}\mathcal S\leq\mathcal E(h_0,u_0,\Omega_0)+C(M)\epsilon,
\]
where
\[
\mathcal S=\Vert\Gamma\Vert_{H^{k+1}}^2+\epsilon^2\Vert\Gamma\Vert_{H^{k+2}}^2. 
\]
\end{lemma}
\proofheading{Proof}We have
\[
E^w(h_0|_\Omega,u_0|_\Omega)+E^r(u_0|_\Omega)+\epsilon^4\Vert\curl u_0\Vert_{W^{k+1,2}(\Omega)}^2\leq E^w(h_0,u_0)+E^r(u_0)+\epsilon^4\Vert\curl u_0\Vert_{W^{k+1,2}(\Omega_0)}^2. 
\]
It remains to show
\begin{equation}\label{eq:kappadissip}
\Vert\tilde a^\frac{1}{2}\Delta_\Gamma^{k_0-1}\kappa\Vert_{L^2(\Gamma)}^2+\frac{\epsilon^2}{C(M)}\mathcal S\leq\Vert a_0^\frac{1}{2}\Delta_{\Gamma_0}^{k_0-1}\kappa_0\Vert_{L^2(\Gamma_0)}^2+C(M)\epsilon,
\end{equation}
where $\tilde a=-n\cdot\nabla h_0|_\Gamma$. Choose local coordinates $(\bar x_\alpha)_{\alpha=1}^2$ on $\Gamma_*$. The maps $Y_1$ and $Y_0$ induce coordinates on $\Gamma$ and $\Gamma_0$, respectively. In these coordinates, the mean curvature is
\[
\kappa=\sum_{\alpha,\beta=1}^2G^{\alpha\beta}(\eta,\partial\eta)\partial_\alpha\partial_\beta\eta:=L\eta,\qquad\kappa_0=\sum_{\alpha,\beta=1}^2G^{\alpha\beta}(\eta_0,\partial\eta_0)\partial_\alpha\partial_\beta\eta_0:=L_0\eta_0. 
\]
The Laplacian can be written as
\[
\Delta_\Gamma=\nu^{-1}\partial_\alpha(\nu g^{\alpha\beta}\partial_\beta),\qquad\Delta_{\Gamma_0}=\nu_0^{-1}\partial_\alpha(\nu_0g_0^{\alpha\beta}\partial_\beta),\qquad\Delta_{\Gamma_*}=\nu_*^{-1}\partial_\alpha(\nu_*g_*^{\alpha\beta}\partial_\beta). 
\]
Let $\chi$ be a cut-off function supported in the domain of $(\bar x_\alpha)_{\alpha=1}^2$. Then the boundary energy, localized by $\chi$, can be expressed as
\begin{align}
&I=\Vert\tilde a^\frac{1}{2}\chi\Delta_\Gamma^{k_0-1}\kappa\Vert_{L^2(\Gamma)}^2=\int_{\mathbb R^2}\tilde a\chi(\Delta_\Gamma^{k_0-1} L\eta)^2\,\nu d\bar x,\nonumber\\
&I_0=\Vert a_0^\frac{1}{2}\chi\Delta_{\Gamma_0}^{k_0-1}\kappa_0\Vert_{L^2(\Gamma_0)}^2=\int a_0\chi(\Delta_{\Gamma_0}^{k_0-1}L_0\eta_0)^2\,\nu_0d\bar x. \label{eq:Eb0loc}
\end{align}
The bounds $|a_0(\bar x)-\tilde a(\bar x)|\lesssim_A\epsilon$ and $|\nu-\nu_0|\lesssim_M\epsilon$ allow us to replace $a_0,\nu_0$ in \eqref{eq:Eb0loc} by $\tilde a,\nu$, up to the required error. By \eqref{eq:etaregerr}, we may also replace $\Delta_{\Gamma_0},L_0$ by $\Delta_\Gamma,L$. Consequently,
\[
\Vert a_0^\frac{1}{2}\chi\Delta_{\Gamma_0}^{k_0-1}\kappa_0\Vert_{L^2(\Gamma_0)}^2=\int_{\mathbb R^2}\tilde a\chi(\Delta_\Gamma^{k_0-1} L\eta_0)^2\,\nu d\bar x+O(\epsilon). 
\]
Next, by $\eta_0=(1-\epsilon^2\Delta_{\Gamma_*})\eta+C_0\epsilon^2$ we have
\[\begin{aligned}
I-I_0&=-\int(\Delta_\Gamma^{k_0-1} L(\eta-\eta_0))^2\,\chi\tilde a\nu d\bar x+2\epsilon^2\int(\Delta_\Gamma^{k_0-1} L\eta)\Delta_\Gamma^{k_0-1} L\Delta_{\Gamma_*}\eta\,\chi\tilde a\nu d\bar x+O(\epsilon)\\
&=-\epsilon^4\int(\Delta_\Gamma^{k_0-1} L\Delta_{\Gamma_*}\eta)^2\,\chi\tilde a\nu d\bar x-2\epsilon^2\int|\nabla_{\Gamma_*}\Delta_\Gamma^{k_0-1} L\eta|^2\,\chi\tilde a\nu d\bar x+J+O(\epsilon)
\end{aligned}\]
where $|\nabla_*f|^2=g_*^{\alpha\beta}\partial_\alpha f\partial_\beta f$ and
\[
J=2\epsilon^2\int(\Delta_\Gamma^{k_0-1} L\eta)[\Delta_\Gamma^{k_0-1} L,\Delta_{\Gamma_*}]\eta\,\chi\tilde a\nu d\bar x-2\epsilon^2\int(\Delta_\Gamma^{k_0-1} L\eta)g_*^{\alpha\beta}\partial_\beta(\Delta_\Gamma^{k_0-1} L\eta)\nu_*\partial_\alpha\left(\frac{\chi\tilde a\nu}{\nu_*}\right)d\bar x. 
\]
Bilinear estimate implies $|J|\lesssim_M\epsilon$. Moreover, as in Proposition~\ref{prop:LapGam} one can show that
\[
\Vert\chi\eta(\bar x)\Vert_{H^{k+1}(\mathbb R^2)}^2\lesssim_M\int|\nabla_{\Gamma_*}\Delta_\Gamma^{k_0-1} L\eta|^2\chi\tilde a\nu d\bar x+\epsilon^{-1},
\]
\[
\Vert\chi\eta(\bar x)\Vert_{H^{k+2}(\mathbb R^2)}^2\lesssim_M\int(\Delta_\Gamma^{k_0-1} L\Delta_{\Gamma_*}\eta)^2\chi\tilde a\nu d\bar x+\epsilon^{-3}. 
\]
Now using Young's inequality and by a partition of unity we conclude \eqref{eq:kappadissip} and complete the proof of the lemma. \proofqed

\subsection{Regularization of the acoustic component}\label{sub:acousticreg}

We regularize the highest acoustic derivatives by elliptic resolvents and then reconstruct the lower derivatives with boundary data determined by compatibility. Define $H_{k-2},U_{k-1}$ by
\begin{equation}\label{eq:Htop}
\left\{\begin{aligned}
&H_{k-2}-C_1\epsilon^2\left(\Delta H_{k-2}+\Phi_{\epsilon^{-1}}N_{k,0}\right)=\Delta^{k_0-1}h_0\quad&\text{in }\Omega\\
&H_{k-2}=P_{<\epsilon^{-1}}\left(\Delta^{k_0-1}h_0|_\Gamma\right)\quad&\text{on }\Gamma
\end{aligned}\right.
\end{equation}
and
\begin{equation}\label{eq:Utop}
\left\{\begin{aligned}
&U_{k-1}-C_1\epsilon^2\Delta U_{k-1}=\Delta^{k_0-1}\divg u_0\quad&\text{in }\Omega\\
&U_{k-1}=P_{<\epsilon^{-1}}\left(\Delta^{k_0-1}\divg u_0|_\Gamma\right)\quad&\text{on }\Gamma
\end{aligned}\right.
\end{equation}
where $C_1=O(1)$ is a sufficiently large constant to be fixed below. Define $H_{2j},U_{2j+1}$ recursively by
\[%\begin{equation}\label{eq:Hint}
\left\{\begin{aligned}
&\Delta H_{2j}=H_{2j+2}\quad&\text{in }\Omega\\
&H_{2j}=P_{<\epsilon^{-1}}\left(\Delta^jh_0|_\Gamma\right)\quad&\text{on }\Gamma
\end{aligned}\right.
\]%\end{equation}
for $0<2j<k-2$, 
\[%\begin{equation}\label{eq:Uint}
\left\{\begin{aligned}
&\Delta U_{2j+1}=U_{2j+3}\quad&\text{in }\Omega\\
&U_{2j+1}=P_{<\epsilon^{-1}}\left(\Delta^j\divg u_0|_\Gamma\right)\quad&\text{on }\Gamma
\end{aligned}\right.
\]%\end{equation}
for $2j+1<k-1$, and define $h_\epsilon=\Delta^{-1}H_2$, $u_\epsilon=u_0|_\Omega+\nabla\Delta^{-1}(U_1-\divg u_0)$. In other words, we have
\begin{equation}\label{eq:h0decomp}
h_0=h_\epsilon+\mathcal H(h_0|_\Gamma)+\sum_{j=1}^{k_0-2}\Delta^{-j}\mathcal HP_{\geq\epsilon^{-1}}(\Delta^jh_0|_\Gamma)-C_1\epsilon^2\Delta^{1-k_0}\left(\Delta H_{k-2}+\Phi_{\epsilon^{-1}}N_{k,0}\right),
\end{equation}
and
\begin{equation}\label{eq:u0decomp}
u_0=u_\epsilon-\sum_{j=0}^{k_0-2}\nabla\Delta^{-1-j}\mathcal HP_{\geq\epsilon^{-1}}(\Delta^j\divg u_0|_\Gamma)-C_1\epsilon^2\nabla\Delta^{-k_0}\Delta U_{k-1}. 
\end{equation}

By the balanced elliptic estimate, Bernstein's estimate and \eqref{eq:numrel} we have
\[
\Vert H_2\Vert_{H^{k-2+r}(\Omega)}+\Vert U_1\Vert_{H^{k-1+r}(\Omega)}\lesssim_M\epsilon^{-r},\quad r\in[0,2],
\]
\[
\Vert H_2-\Delta h_0\Vert_{H^{k-2-r}(\Omega)}+\Vert U_1-\divg u_0\Vert_{H^{k-1-r}(\Omega)}\lesssim_M\epsilon^r,\quad r\in[0,2]. 
\]
Combining \eqref{eq:zj0tr}--\eqref{eq:h0tr} we derive
\begin{equation}\label{eq:wavereg}
\Vert h_\epsilon\Vert_{H^{k+r}(\Omega)}+\Vert\divg u_\epsilon\Vert_{H^{k-1+r}(\Omega)}\lesssim_M\epsilon^{-r},\quad r\in[0,2],
\end{equation}
\begin{equation}\label{eq:uepsreg}
\Vert u_\epsilon\Vert_{H^{k+1}(\Omega)}\lesssim_MP\epsilon^{-1}, 
\end{equation}
and
\begin{equation}\label{eq:huepserr}
\Vert h_\epsilon-h_0|_\Omega\Vert_{H^{k-r}(\Omega)}+\Vert u_\epsilon-u_0|_\Omega\Vert_{H^{k-r}(\Omega)}\lesssim_M\epsilon^r,\quad r\in[0,2]. 
\end{equation}
From the boundary conditions we see 
\[
z_{i,\epsilon}|_\Gamma=P_{<\epsilon^{-1}}(z_{i,0}|_\Gamma)+P_{\geq\epsilon^{-1}}(N_{i,0}|_\Gamma)+N_{i,\epsilon}-N_{i,0},\quad i=1,2,...,k-1. 
\]
From \eqref{eq:Htop}--\eqref{eq:Utop} we see
\[\begin{aligned}
z_{k,\epsilon}|_\Gamma=C_1^{-1}\epsilon^{-2}P_{\geq\epsilon^{-1}}(N_{k-2,0}|_\Gamma)-C_1^{-1}\epsilon^{-2}P_{\geq\epsilon^{-1}}(z_{k-2,0}|_\Gamma)+N_{k,\epsilon}-\Phi_{\epsilon^{-1}}N_{k,0},
\end{aligned}\]
\[\begin{aligned}
\Delta^{k_0}\divg u_\epsilon|_\Gamma=C_1^{-1}\epsilon^{-2}P_{\geq\epsilon^{-1}}(N_{k-1,0}|_\Gamma)-C_1^{-1}\epsilon^{-2}P_{\geq\epsilon^{-1}}(z_{k-1,0}|_\Gamma). 
\end{aligned}\]
These together with \eqref{eq:zj0tr}, \eqref{eq:wavereg}--\eqref{eq:huepserr} imply 
\begin{equation}\label{eq:zjepstr}
\Vert z_{j,\epsilon}|_\Gamma\Vert_{H^{k-r-j}(\Gamma)}\lesssim_MP\epsilon^{r+\frac{1}{2}},\quad -\frac{3}{2}\leq r\leq\min\{k-j,\frac{3}{2}\},\ j=1,...,k+1.
\end{equation}

\begin{lemma}\label{lem:Emonoeps}
\[
\Vert(L_\zeta^{k_0}h_\epsilon, L_\zeta^{k_0}u_\epsilon)\Vert_{L^2(\Omega)}^2\leq \Vert(L_\zeta^{k_0}h_0, L_\zeta^{k_0}u_0)\Vert_{L^2(\Omega)}^2+C(P)\mathcal S^\frac{1}{2}\epsilon^\frac{3}{2}+C(P)\epsilon. 
\]
\end{lemma}
\proofheading{Proof}\textit{Step 1: compare the enthalpy with its homogeneous Dirichlet part.} We prove
\begin{equation}\label{eq:h0htd}
\Vert L_\zeta^{k_0}\tilde h\Vert_{L^2(\Omega)}^2\leq\Vert L_\zeta^{k_0}h_0\Vert_{L^2(\Omega)}^2+C(P)\mathcal S^\frac{1}{2}\epsilon^\frac{3}{2}+C(P)\epsilon,
\end{equation}
where $\tilde h=h_0-\mathcal H(h_0|_\Gamma)$. Write
\[
L_\zeta\tilde h=\Delta^{-1}\Delta L_\zeta\tilde h+\mathcal H\left(L_\zeta\tilde h|_\Gamma\right),
\]
\[
L_\zeta h_0=\Delta^{-1}\Delta L_\zeta h_0+\mathcal H\left(L_{\zeta_0}h_0|_\Gamma\right)+\mathcal H\left((\zeta^{ij}-\zeta_0^{ij})\partial_i\partial_jh_0\right).
\]
By Lemma~\ref{lem:zetaerr} we have
\[
\Vert\mathcal H\left((\zeta^{ij}-\zeta_0^{ij})\partial_i\partial_jh_0\right)\Vert_{H^{k-2}(\Omega)}\lesssim_M\epsilon. 
\]
By $\Delta\tilde h=\Delta h_0$ we have
\[
\Vert\Delta^{-1}\Delta L_\zeta(\tilde h-h_0)\Vert_{H^{k-2}(\Omega)}=\Vert\Delta^{-1}[\Delta,L_\zeta](\tilde h-h_0)\Vert_{H^{k-2}(\Omega)}\lesssim_M\epsilon. 
\]
It then remains to show
\begin{equation}\label{eq:h0htd2}
\Vert L_\zeta\tilde h-L_{\zeta_0}h_0\Vert_{H^{k-\frac{5}{2}}(\Gamma)}\lesssim_P\mathcal S^\frac{1}{2}\epsilon^\frac{3}{2}+\epsilon. 
\end{equation}
Recall the homotopy $Y_t$ from the last subsection. The left-hand side is controlled by $\Vert(L_\zeta\tilde h-L_{\zeta_0}h_0)\circ Y_1\Vert_{H^{k-\frac{5}{2}}(\Gamma_*)}$. By \eqref{eq:chainrule}, \eqref{eq:reg0} and the bilinear estimate we have
\[
\Vert L_{\zeta_0}h_0\circ Y_1-L_{\zeta_0}h_0\circ Y_0\Vert_{H^{k-\frac{5}{2}}(\Gamma_*)}\lesssim_MP\epsilon. 
\]
Next, notice that
\[
L_\zeta\tilde h\circ Y_1=-(\tilde a\kappa)\circ Y_1,\quad L_{\zeta_0}h_0\circ Y_0=-(a_0\kappa_0)\circ Y_0. 
\]
From \eqref{eq:h0tr}, \eqref{eq:chainrule} we see $\Vert\tilde a(Y_1)-a_0(Y_0)\Vert_{H^{k-\frac{5}{2}}(\Gamma_*)}\lesssim_M\epsilon$. Moreover, from the proof of \eqref{eq:kappadissip} we see
\[
\Vert\kappa(Y_1)-\kappa_0(Y_0)\Vert_{H^{k-\frac{5}{2}}(\Gamma_*)}\lesssim_M\Vert\Gamma\Vert_{H^{k+\frac{3}{2}}}\epsilon^2.
\]
By interpolation the right-hand side is controlled by $\mathcal S^\frac{1}{2}\epsilon^\frac{3}{2}$. Combining these we obtain \eqref{eq:h0htd2} and conclude \eqref{eq:h0htd}. 

\textit{Step 2: estimate the tangential enthalpy energy after acoustic smoothing.} We prove
\begin{equation}\label{eq:htdheps}
\Vert L_\zeta^{k_0}h_\epsilon\Vert_{L^2(\Omega)}^2\leq\Vert L_\zeta^{k_0}\tilde h\Vert_{L^2(\Omega)}^2+C(P)\mathcal S^\frac{1}{2}\epsilon^\frac{3}{2}+C(P)\epsilon. 
\end{equation}
Equation \eqref{eq:h0decomp} gives
\[
\tilde h=h_\epsilon+\sum_{j=1}^{k_0-2}\Delta^{-j}\mathcal HP_{\geq\epsilon^{-1}}(\Delta^jh_0|_\Gamma)-C_1\epsilon^2\Delta^{1-k_0}\left(\Delta H_{k-2}+\Phi_{\epsilon^{-1}}N_{k,0}\right)
\]
By \eqref{eq:zj0tr}, \eqref{eq:reg0}, Bernstein's estimate and Lemma~\ref{lem:Phim} we have %$\tilde h-h_\epsilon+C_1\epsilon^2\Delta^{1-k_0}\Delta H_{k-2}=O(\epsilon)$ in $H^k(\Omega)$. 
\begin{equation}\label{eq:midhar}
\Vert\tilde h-h_\epsilon+C_1\epsilon^2\Delta^{1-k_0}\Delta H_{k-2}\Vert_{H^k(\Omega)}\lesssim_MP\epsilon. 
\end{equation}
It follows that
\begin{equation}\label{eq:polar}
\Vert L_\zeta^{k_0}\tilde h\Vert_{L^2(\Omega)}^2=\Vert L_\zeta^{k_0}h_\epsilon\Vert_{L^2(\Omega)}^2+2I+\Vert L_\zeta^{k_0}(\tilde h-h_\epsilon)\Vert_{L^2(\Omega)}^2
\end{equation}
with
\[
I=\int_\Omega L_\zeta^{k_0}(\tilde h-h_\epsilon)\,L_\zeta^{k_0}h_\epsilon\,dx=-C_1\epsilon^2\int_\Omega L_\zeta^{k_0}\Delta^{1-k_0}\Delta H_{k-2}\,L_\zeta^{k_0}h_\epsilon\,dx+O(\epsilon). 
\]
Write $\Delta^{1-k_0}\Delta H_{k-2}=\Delta\Delta^{1-k_0}H_{k-2}-\Delta^{2-k_0}\mathcal H_\Gamma H_{k-2}$. By \eqref{eq:Htop}, \eqref{eq:zj0tr} we have $\Vert H_{k-2}\Vert_{H^\frac{7}{2}(\Gamma)}\lesssim_MP\epsilon^{-1}$, and by the elliptic estimate we have $\Vert\Delta^{2-k_0}\mathcal H_\Gamma H_{k-2}\Vert_{H^k(\Omega)}\lesssim_MP\epsilon^{-1}$. Thus, the integral $I$ can be reduced to 
\[\begin{aligned}
I&=-C_1\epsilon^2\int_\Omega L_\zeta^{k_0}\Delta\Delta^{1-k_0}H_{k-2}\,L_\zeta^{k_0}h_\epsilon\,dx+O(\epsilon)\\
&=-C_1\epsilon^2\int_\Omega\Delta L_\zeta^{k_0}\Delta^{1-k_0}H_{k-2}\,L_\zeta^{k_0}h_\epsilon\,dx+O(\epsilon),
\end{aligned}\]
where the second line is due to
\[
\Vert[L_\zeta^{k_0},\Delta]\Delta^{1-k_0}H_{k-2}\Vert_{L^2(\Omega)}\lesssim_M\Vert\Delta^{1-k_0}H_{k-2}\Vert_{H^{k+1}(\Omega)}\lesssim_M\epsilon^{-1}. 
\]
Next, $h_\epsilon$ is decomposed into
\[
h_\epsilon=h_\epsilon^0+h_\epsilon^1+h_\epsilon^2:=\Delta^{1-k_0}H_{k-2}-\sum_{j=1}^{k_0-2}\Delta^{-j}\mathcal H_\Gamma N_{2j,\epsilon}+\sum_{j=1}^{k_0-2}\Delta^{-j}\mathcal H_\Gamma z_{2j,\epsilon}.
\]
By \eqref{eq:zjepstr}, \eqref{eq:wavereg}--\eqref{eq:uepsreg}, elliptic and trace estimates we have
\[
\Vert h_\epsilon^2\Vert_{H^k(\Omega)}\lesssim_M\epsilon,\quad\Vert h_\epsilon^0\Vert_{H^{k+r}(\Omega)}\lesssim_M\epsilon^{-r},\quad r\in[0,2]
\]
and 
\begin{equation}\label{eq:heps1reg}
\Vert\Delta h_\epsilon^1\Vert_{H^{k-1}(\Omega)}\lesssim_M1,\quad\Vert\Delta h_\epsilon^1\Vert_{H^k(\Omega)}\lesssim_MP\epsilon^{-1},\quad\Vert h_\epsilon^1\Vert_{H^{k+r}(\Omega)}\lesssim_M\epsilon^{\frac{1}{2}-r},\ r\in[\frac{1}{2},2]. 
\end{equation}
Accordingly, $I$ is decomposed into
\[
I=\sum_{\alpha=0}^1I^\alpha+O(\epsilon)=\sum_{\alpha=0}^1(-C_1\epsilon^2)\int_\Omega\Delta L_\zeta^{k_0}h_\epsilon^0\,L_\zeta^{k_0}h_\epsilon^\alpha\,dx+O(\epsilon). 
\]
Green's formula rewrites $I^1$ as
\[\begin{aligned}
I^1=-C_1\epsilon^2\int_\Omega L_\zeta^{k_0}h_\epsilon^0\,\Delta L_\zeta^{k_0}h_\epsilon^1\,dx+C_1\epsilon^2\int_\Gamma L_\zeta^{k_0}h_\epsilon^0\,\nabla_nL_\zeta^{k_0}h_\epsilon^1-\nabla_nL_\zeta^{k_0}h_\epsilon^0\,L_\zeta^{k_0}h_\epsilon^1\,dS.
\end{aligned}\]
The interior integral on the right-hand side is $O(\epsilon)$ due to
\[
\Vert\Delta L_\zeta^{k_0}h_\epsilon^1\Vert_{L^2(\Omega)}\leq\Vert[\Delta,L_\zeta^{k_0}]h_\epsilon^1\Vert_{L^2(\Omega)}+\Vert L_\zeta^{k_0}\Delta h_\epsilon^1\Vert_{L^2(\Omega)}\lesssim_MP\epsilon^{-1},
\]
which is a consequence of \eqref{eq:heps1reg}. To control the boundary integral, we notice for every $f$ vanishing on $\Gamma$ that
\begin{equation}\label{eq:Lzk0tr}
\Vert L_\zeta^{k_0}f\Vert_{H^r(\Gamma)}\lesssim_M\Vert f\Vert_{H^{k-\frac{1}{2}+r}(\Omega)}+\Vert f\Vert_{C^2(\Omega)}\Vert\Gamma\Vert_{H^{k+r}},\quad r\in[0,1],
\end{equation}
which can be seen from \eqref{eq:Lztr}. Applying this to $h_0^\epsilon$ and using that
\[
\Vert\nabla_nL_\zeta^{k_0}h_\epsilon^1\Vert_{L^2(\Gamma)}\lesssim_M\Vert\nabla L_\zeta^{k_0}h_\epsilon^1\Vert_{L^2(\Omega)}^\frac{1}{2}\Vert\nabla L_\zeta^{k_0}h_\epsilon^1\Vert_{H^1(\Omega)}^\frac{1}{2}\lesssim_M\epsilon^{-1},
\]
we obtain $\big|C_1\epsilon^2\int_\Gamma L_\zeta^{k_0}h_\epsilon^0\,\nabla_nL_\zeta^{k_0}h_\epsilon^1\,dS\big|\lesssim_M\epsilon$. For the other term, commute $\nabla_n$ with $\zeta^{ij}\partial_i$ and integrate by parts: 
\[\begin{aligned}
\left|\int_\Gamma\nabla_nL_\zeta^{k_0}h_\epsilon^0\,L_\zeta^{k_0}h_\epsilon^1\,dS\right|&=\left|\int_\Gamma\zeta^{ij}\partial_i\nabla_n\partial_jL_\zeta^{k_0-1}h_\epsilon^0\,L_\zeta^{k_0}h_\epsilon^1\,dS\right|+O(\epsilon^{-1})\\
&\lesssim_M\Vert\nabla_n\partial_jL_\zeta^{k_0-1}h_\epsilon^0\Vert_{L^2(\Gamma)}\Vert L_\zeta^{k_0}h_\epsilon^1\Vert_{H^1(\Gamma)}\\
&\lesssim_M\epsilon^{-\frac{1}{2}}\left(1+\Vert\Gamma\Vert_{H^{k+1}}\right),
\end{aligned}\]
where the last line is due to \eqref{eq:Lzk0tr} and \eqref{eq:heps1reg}. Combining these we obtain 
\[
|I^1|\lesssim_P\mathcal S^\frac{1}{2}\epsilon^\frac{3}{2}+\epsilon.
\]
It remains to estimate $I^0$. Write $f=f_{ho}+f_{in}:=\Delta^{-1}\Delta f+\mathcal H_\Gamma f$ and notice that $\int_\Omega\nabla f_{ho}\cdot\nabla f_{in}\,dx=0$. We then derive
\[
I^0=C_1\epsilon^2\int_\Omega|\nabla(L_\zeta^{k_0}h_\epsilon^0)_{ho}|^2\,dx-C_1\epsilon^2\int_\Gamma\nabla_n(L_\zeta^{k_0}h_\epsilon^0)_{ho}\,L_\zeta^{k_0}h_\epsilon^0\,dS. 
\]
By \eqref{eq:Lzk0tr}, trace and elliptic estimate we obtain 
\[\begin{aligned}
\int_\Gamma\nabla_n(L_\zeta^{k_0}h_\epsilon^0)_{ho}\,L_\zeta^{k_0}h_\epsilon^0\,dS\lesssim_M\Vert\nabla(L_\zeta^{k_0}h_\epsilon^0)_{ho}\Vert_{L^2(\Omega)}^\frac{1}{2}\Vert\Delta L_\zeta^{k_0}h_\epsilon^0\Vert_{L^2(\Omega)}^\frac{1}{2}. 
\end{aligned}\]
We now use the positive third term on the right-hand side of \eqref{eq:polar}. By \eqref{eq:midhar},
\[
\int_\Omega\left(L_\zeta^{k_0}(\tilde h-h_\epsilon)\right)^2\,dx=C_1^2\epsilon^4\int_\Omega\left(L_\zeta^{k_0}\Delta h_\epsilon^0\right)^2\,dx+O(\epsilon)=C_1^2\epsilon^4\int_\Omega\left(\Delta L_\zeta^{k_0}h_\epsilon^0\right)^2\,dx+O(\epsilon). 
\]
Then by Young's inequality we conclude \eqref{eq:htdheps}. 

\textit{Step 3: estimate the tangential velocity energy.} We prove
\begin{equation}\label{eq:u0ueps}
\Vert L_\zeta^{k_0}u_\epsilon\Vert_{L^2(\Omega)}^2\leq\Vert L_\zeta^{k_0}u_0\Vert_{L^2(\Omega)}^2+C(P)\epsilon. 
\end{equation}
As in Step 2, we make use of the decomposition \eqref{eq:u0decomp}. As in \eqref{eq:midhar} we have
\[
\Vert\nabla\Delta^{-1-j}\mathcal HP_{\geq\epsilon^{-1}}(\Delta^j\divg u_0|_\Gamma)\Vert_{H^k(\Omega)}\lesssim_MP\epsilon,\quad j=0,...,k_0-2. 
\]
Thus, the polar identity now reads as
\[
\Vert L_\zeta^{k_0}u_0\Vert_{L^2(\Omega)}^2=\Vert L_\zeta^{k_0}u_\epsilon\Vert_{L^2(\Omega)}^2+2J+\Vert L_\zeta^{k_0}(u_0-u_\epsilon)\Vert_{L^2(\Omega)}^2
\]
with
\[
J=-C_1\epsilon^2\int_\Omega L_\zeta^{k_0}\nabla\Delta^{1-k_0}U_{k-1}\,L_\zeta^{k_0}u_\epsilon\,dx+O(\epsilon). 
\]
Commute $L_\zeta^{k_0}$ with $\nabla$, integrate by parts and commute derivatives again. We then get
\[
J=C_1\epsilon^2\int_\Omega L_\zeta^{k_0}\Delta^{1-k_0}U_{k-1}\,L_\zeta^{k_0}U_1\,dx-C_1\epsilon^2\int_\Gamma L_\zeta^{k_0}\Delta^{1-k_0}U_{k-1}\,L_\zeta^{k_0}u_\epsilon\cdot n\,dS+O(\epsilon),
\]
where we have noticed that
\[
\Vert[L_\zeta^{k_0},\nabla]\Delta^{1-k_0}U_{k-1}\Vert_{L^2(\Omega)}+\Vert[\divg,L_\zeta^{k_0}]u_\epsilon\Vert_{L^2(\Omega)}\lesssim_MP\epsilon^{-1}. 
\]
By \eqref{eq:Lzk0tr} and trace estimate, the boundary integral is bounded by
\[\begin{aligned}
\left|\int_\Gamma L_\zeta^{k_0}\Delta^{1-k_0}U_{k-1}\,L_\zeta^{k_0}u_\epsilon\cdot n\,dS\right|&\lesssim_M\Vert L_\zeta^{k_0}\Delta^{1-k_0}U_{k-1}\Vert_{L^2(\Gamma)}\Vert L_\zeta^{k_0}u_\epsilon\Vert_{L^2(\Omega)}^\frac{1}{2}\Vert L_\zeta^{k_0}u_\epsilon\Vert_{H^1(\Omega)}^\frac{1}{2}\\
&\lesssim_M(\epsilon^{-\frac{1}{2}}+1)\epsilon^{-\frac{1}{2}}\lesssim_M\epsilon^{-1}. 
\end{aligned}\]
Next, write
\[
U_1=U_1^0+U_1^1+U_1^2:=\Delta^{1-k_0}U_{k-1}-\sum_{j=0}^{k_0-2}\Delta^{-j}\mathcal H_\Gamma N_{2j+1,\epsilon}+\sum_{j=0}^{k_0-2}\Delta^{-j}\mathcal H_\Gamma z_{2j+1,\epsilon}. 
\]
From \eqref{eq:zjepstr} we see
\[
\Vert U_1^2\Vert_{H^{k-1}(\Omega)}\lesssim_MP\epsilon,\quad\Vert U_1^1\Vert_{H^k(\Omega)}\lesssim_M1. 
\]
Thus, using the divergence theorem and H\"older's inequality, we obtain
\[\begin{aligned}
\left|\int_\Omega L_\zeta^{k_0}U_1^0\,L_\zeta^{k_0}U_1^2\,dx\right|&=\left|\int_\Omega\partial_i\left(\zeta^{ij}L_\zeta^{k_0}U_1^0\right)\,\partial_jL_\zeta^{k_0-1}U_1^2\,dx\right|\\
&\lesssim_M\left(1+\Vert U_1^0\Vert_{H^{k+1}(\Omega)}\right)\epsilon\lesssim_M\epsilon^{-1}
\end{aligned}\]
and
\[
\left|\int_\Omega L_\zeta^{k_0}U_1^0\,L_\zeta^{k_0}U_1^1\,dx\right|\lesssim_M1+\Vert U_1^0\Vert_{H^k(\Omega)}\lesssim_M\epsilon^{-1}. 
\]
This implies \eqref{eq:u0ueps} and completes the proof of the lemma. \proofqed

\begin{lemma}\label{lem:wavedissip}
Define
\[
\mathcal D=\int_\Omega|\nabla\Delta H_{k-2}|^2+(\Delta U_{k-1})^2+\frac{\epsilon^2}{2}(\Delta^2H_{k-2})^2+\frac{\epsilon^2}{2}|\nabla\Delta U_{k-1}|^2\,dx. 
\]
Then we have
\[
\Vert(z_{k,\epsilon},w_{k,\epsilon})\Vert_{L^2(\Omega)}^2+C_1\epsilon^2\mathcal D\leq\Vert(z_{k,0},w_{k,0})\Vert_{L^2(\Omega)}^2+C(P)\epsilon. 
\]
\end{lemma}
\proofheading{Proof}By \eqref{eq:huepserr} we have
\[\begin{aligned}
\Vert(z_{k,0},w_{k,0})\Vert_{L^2(\Omega)}^2&-\Vert(z_{k,\epsilon},w_{k,\epsilon})\Vert_{L^2(\Omega)}^2\\
&=\Vert(\Delta^{k_0}h_0,\nabla\Delta^{k_0-1}\divg u_0)\Vert_{L^2(\Omega)}^2-\Vert(\Delta^{k_0}h_\epsilon,\nabla\Delta^{k_0-1}\divg u_\epsilon)\Vert_{L^2(\Omega)}^2+O(\epsilon). 
\end{aligned}\]
We first estimate the contribution of $h$. Write $\Vert\Delta^{k_0}h_0\Vert_{L^2(\Omega)}^2=\Vert\Delta^{k_0}h_\epsilon\Vert_{L^2(\Omega)}^2+2I+\Vert\Delta^{k_0}(h_0-h_\epsilon)\Vert_{L^2(\Omega)}^2$ with
\[\begin{aligned}
I&=\int_\Omega\Delta^{k_0}(h_0-h_\epsilon)\,\Delta^{k_0}h_\epsilon\,dx\\
&=-C_1\epsilon^2\int_\Omega\Delta^2H_{k-2}\Delta H_{k-2}\,dx+O(\epsilon)\\
&=C_1\epsilon^2\int_\Omega|\nabla\Delta H_{k-2}|^2\,dx-C_1\epsilon^2\int_\Gamma\nabla_n\Delta H_{k-2}\,\Delta H_{k-2}\,dS. 
\end{aligned}\]
By \eqref{eq:zjepstr}, the boundary integral is controlled by
\[\begin{aligned}&
\left|\int_\Gamma\nabla_n\Delta H_{k-2}\,\Delta H_{k-2}\,dS\right|\\ &\leq\left|\int_\Gamma\mathcal N\Delta H_{k-2}\,\Delta H_{k-2}\,dS\right|+\left|\int_\Gamma\nabla_n(\Delta H_{k-2})_{ho}\,\Delta H_{k-2}\,dS\right|\\
&\lesssim_M\Vert\Delta H_{k-2}\Vert_{H^\frac{1}{2}(\Gamma)}^2+\Vert\nabla(\Delta H_{k-2})_{ho}\Vert_{L^2(\Omega)}^\frac{1}{2}\Vert\Delta^2H_{k-2}\Vert_{L^2(\Omega)}^\frac{1}{2}\Vert\Delta H_{k-2}\Vert_{L^2(\Gamma)}\\
&\lesssim_M1+\epsilon^{-\frac{1}{2}}\mathcal D^\frac{1}{2},
\end{aligned}\]
where in the last line we have used that $\Vert\nabla f_{ho}\Vert_{L^2(\Omega)}^2\leq\Vert\nabla f\Vert_{L^2(\Omega)}^2$. 

We next estimate the contribution of $u$. Write 
\[
\Vert\nabla\Delta^{k_0-1}\divg u_0\Vert_{L^2(\Omega)}^2=\Vert\nabla\Delta^{k_0-1}\divg u_\epsilon\Vert_{L^2(\Omega)}^2+2J+\Vert\nabla\Delta^{k_0-1}\divg(u_0-u_\epsilon)\Vert_{L^2(\Omega)}^2
\]
with
\[\begin{aligned}
J&=\int_\Omega\nabla\Delta^{k_0-1}\divg(u_0-u_\epsilon)\cdot\nabla\Delta^{k_0-1}\divg u_\epsilon\,dx\\
&=-C_1\epsilon^2\int_\Omega\nabla\Delta U_{k-1}\cdot\nabla U_{k-1}\,dx\\
&=C_1\epsilon^2\int_\Omega(\Delta U_{k-1})^2\,dx-C_1\epsilon^2\int_\Omega\Delta U_{k-1}\,\nabla_nU_{k-1}\,dS. 
\end{aligned}\]
Again by \eqref{eq:zjepstr}, the boundary integral is controlled by
\[\begin{aligned}
\left|\int_\Omega\Delta U_{k-1}\,\nabla_nU_{k-1}\,dS\right|\lesssim_M\Vert\Delta U_{k-1}\Vert_{L^2(\Gamma)}\Vert U_{k-1}\Vert_{H^1(\Omega)}^\frac{1}{2}\Vert U_{k-1}\Vert_{H^2(\Omega)}^\frac{1}{2}\lesssim_MP\epsilon^{-1}. 
\end{aligned}\]

Finally, note that 
\[
\Vert\Delta^{k_0}(h_0-h_\epsilon)\Vert_{L^2(\Omega)}^2=C_1^2\epsilon^4\Vert\Delta^2H_{k-2}\Vert_{L^2(\Omega)}^2+O(\epsilon),
\]\[
\Vert\nabla\Delta^{k_0-1}\divg(u_0-u_\epsilon)\Vert_{L^2(\Omega)}^2=C_1^2\epsilon^4\Vert\nabla\Delta U_{k-1}\Vert_{L^2(\Omega)}^2. 
\]
Combining the estimates and using Young's inequality we complete the proof. \proofqed

We conclude this subsection with
\begin{equation}\label{eq:Emonoeps}
\mathcal E(h_\epsilon,u_\epsilon,\Omega)+C_1\epsilon^2\mathcal D+\frac{\epsilon^2}{C(M)}\mathcal S\leq\mathcal E(h_0,u_0,\Omega_0)+C(P)\epsilon. 
\end{equation}
This can be seen from $\curl u_\epsilon=\curl u_0$, Lemmas \ref{lem:Emono0}-\ref{lem:wavedissip} and Young's inequality. 

\subsection{Tangential regularization of the velocity}\label{sub:tangreg}

The acoustic regularization controls divergence while leaving vorticity unchanged. We now smooth tangentially to obtain the remaining velocity regularity. Define $X:\Omega\to\mathbb R^3$ by
\[
X=Y_1^{-1}\text{ on }\Gamma,\quad\Delta X=0\text{ in }\Omega.
\]
Then $X$ is a diffeomorphism onto $\Omega_*$ because $\Gamma\in\Lambda_*$ is sufficiently close to $\Gamma_*$. 
Let $\gamma:(-\delta_0,\delta_0)\times\Gamma_*\to\Omega_*$, $\gamma(s,x)=x-sn_*(x)$ be the normal exponential map. Define $\overline\Delta$ via
\begin{equation}\label{eq:Deltabardef}
\overline\Delta=\overline\nabla\cdot\overline\nabla,\quad\overline\nabla f\circ\gamma=\nabla_{\Gamma_*}f(\gamma). 
\end{equation}
Note that $\overline\Delta f|_{\Gamma_*}=\Delta_{\Gamma_*}f$. Also note that
\begin{equation}\label{eq:Deltabarcomm}
[\partial^j,\overline\Delta]=\overline\nabla\cdot P_j(\partial)+Q_j(\partial)
\end{equation}
where $P_j(\partial),Q_j(\partial)$ are $j$th order differential operators. To see this, choose a local orthogonal frame $e_1,e_2,\nu=\partial_s\gamma\circ\gamma^{-1}$. Note that $[\nabla_\nu,\overline\nabla_i]=0$ and that $[\nabla_{e_j},\overline\nabla_i]$ are again linear combinations of $\overline\nabla_i,i=1,2,3$. Finally, let $\chi,\tilde\chi$ be cut-off functions with 
\[
\Gamma_*\subset\subset\{\chi=1\}\subset\supp\chi\subset\{\tilde\chi=1\}\subset\supp\tilde\chi\subset\subset\text{image of }\gamma. 
\]
Choose the cutoffs to be constant on each surface $\gamma(s,\Gamma_*)$, so that $\overline\nabla\chi=\overline\nabla\tilde\chi=0$. Define the tangentially regularized velocity $\bar u_\epsilon$ by
\begin{equation}\label{eq:ubarepsdef}
\bar u_\epsilon(X^{-1})=\chi\bar w_\epsilon+(1-\chi)w_\epsilon,\quad \bar w_\epsilon=(I-\epsilon^2\overline\Delta)^{-1}(\tilde\chi w_\epsilon),\quad w_\epsilon=u_\epsilon(X^{-1}). 
\end{equation}

\begin{lemma}\label{lem:Deltabarellip}
With $J_\gamma=\det(\nabla\gamma^{-1})$, for $f$ supported near $\Gamma_*$ we have
\[
\int_{\Omega_*}|\overline\nabla f|^2\,J_\gamma dx=\int_{\Omega_*}-f\overline\Delta f\,J_\gamma dx
\]
and
\[
\sum_{j=0}^2\Vert\overline\nabla^jf\Vert_{L^2(\Omega_*)}\lesssim\Vert\overline\Delta f\Vert_{L^2(\Omega_*)}+\Vert f\Vert_{L^2(\Omega_*)}. 
\]
\end{lemma}
\proofheading{Proof}By a change of variables and Green's formula on $\Gamma_*$ we have
\[\begin{aligned}
\int_{\Omega_*}|\overline\nabla f|^2\,J_\gamma dx&=\int_{-\delta_0}^{\delta_0}\int_{\Gamma_*}|\nabla_{\Gamma_*}f(\gamma(s,\cdot))|^2\,dS\,ds\\
&=\int_{-\delta_0}^{\delta_0}\int_{\Gamma_*}-f(\gamma(s,\cdot))\Delta_{\Gamma_*}f(\gamma(s,\cdot))\,dS\,ds=\int_{\Omega_*}-f\overline\Delta f\,J_\gamma dx. 
\end{aligned}\]
By the elliptic estimate on $\Gamma_*$ we have
\[\begin{aligned}
\sum_{j=0}^2\Vert\overline\nabla^jf\Vert_{L^2(\Omega_*)}^2&\sim\int_{-\delta_0}^{\delta_0}\sum_{j=0}^2\Vert\nabla_{\Gamma_*}^jf(\gamma(s,\cdot))\Vert_{L^2(\Gamma_*)}^2\,ds\\
&\lesssim\int_{-\delta_0}^{\delta_0}\Vert\Delta_{\Gamma_*}f(\gamma(s,\cdot))\Vert_{L^2(\Gamma_*)}^2+\Vert f(\gamma(s,\cdot))\Vert_{L^2(\Gamma_*)}^2\,ds\sim\Vert\overline\Delta f\Vert_{L^2(\Omega_*)}^2+\Vert f\Vert_{L^2(\Omega_*)}^2. 
\end{aligned}\]

\begin{lemma}\label{lem:Hodgexi}
Denote $\xi_{ij}=(\partial_iX_j)\circ X^{-1}$ and denote $\divg_\xi f=\xi^{ij}\partial_jf_i,\curl_\xi f_{ij}=\xi^{ik}\partial_kf_j-\xi^{jk}\partial_kf_i$ for a vector field $f$ on $\Omega_*$. Then we have
\[
\Vert f\Vert_{H^2(\Omega_*)}\lesssim_M\Vert\chi\overline\nabla^2f\Vert_{L^2(\Omega_*)}+\Vert\divg_\xi f\Vert_{H^1(\Omega_*)}+\Vert\curl_\xi f\Vert_{H^1(\Omega_*)}+\Vert f\Vert_{L^2(\Omega_*)}. 
\]
\end{lemma}
\proofheading{Proof}Interior Hodge estimate implies $\Vert(1-\chi)f\Vert_{H^2(\Omega_*)}\lesssim_M\Vert\divg_\xi f\Vert_{H^1(\Omega_*)}+\Vert\curl_\xi f\Vert_{H^1(\Omega_*)}+\Vert f\Vert_{L^2(\Omega_*)}$. As in \eqref{eq:ptHodge} we have $|\nabla^2(\chi f)|\lesssim_M|\overline\nabla^2(\chi f)|+|\nabla\divg_\xi(\chi f)|+|\nabla\curl_\xi(\chi f)|$. The conclusion follows. \proofqed
We first establish the regularity bound for the tangentially regularized velocity:
\begin{equation}\label{eq:ubarepsreg}
\Vert\bar u_\epsilon\Vert_{H^{k+r}(\Omega)}\lesssim_M\epsilon^{-r},\quad r\in[0,2]. 
\end{equation}
Write $\bar u_\epsilon=(\chi\bar w_\epsilon)\circ X+(1-\chi(X))u_\epsilon$. By interior Hodge estimate we have
\[
\Vert(1-\chi(X))u_\epsilon\Vert_{H^{k+r}(\Omega)}\lesssim\Vert\divg u_\epsilon\Vert_{H^{k-1+r}(\Omega)}+\Vert\curl u_\epsilon\Vert_{H^{k-1+r}(\Omega)}+\Vert u_\epsilon\Vert_{H^r(\Omega)}\lesssim_M\epsilon^{-r},\ r\in[0,2]. 
\]
Thus, \eqref{eq:ubarepsreg} is reduced to
\begin{equation}\label{eq:wbarepsreg}
\Vert\bar w_\epsilon\Vert_{H^{k+r}(\Omega_*)}\lesssim_M\epsilon^{-r},\quad r\in[0,2]. 
\end{equation}
By $\bar w_\epsilon-\epsilon^2\overline\Delta\bar w_\epsilon=\tilde\chi w_\epsilon$, \eqref{eq:Deltabarcomm} and Lemma~\ref{lem:Deltabarellip}, we obtain
\[
\sum_{j=0}^2\sum_{|\alpha|\leq k}\epsilon^j\Vert\overline\nabla^j\partial^\alpha\bar w_\epsilon\Vert_{L^2(\Omega_*)}\lesssim_M1. 
\]
To estimate $\divg_\xi\bar w_\epsilon$ and $\curl_\xi\bar w_\epsilon$, start with
\[\begin{aligned}
\divg_\xi\bar w_\epsilon-\epsilon^2\overline\Delta\divg_\xi\bar w_\epsilon&=\tilde\chi\divg_\xi w_\epsilon+\xi^{ij}\partial_j\tilde\chi w_{\epsilon i}+\epsilon^2\xi^{ij}[\partial_j,\overline\Delta]\bar w_\epsilon\\ &\quad-\epsilon^2\left(\overline\Delta\xi^{ij}\partial_j\bar w_{\epsilon i}-\overline\nabla\xi^{ij}\cdot\overline\nabla\partial_j\bar w_{\epsilon i}\right), 
\end{aligned}\]
perform $\partial^{k+1}$ and commute derivatives via \eqref{eq:Deltabarcomm}. This yields
\[
\partial^{k+1}\divg_\xi\bar w_\epsilon-\epsilon^2\overline\Delta \partial^{k+1}\divg_\xi\bar w_\epsilon=\partial^{k+1}(\tilde\chi\divg_\xi w_\epsilon)+R_{k+1}+\epsilon^2(\overline\nabla\cdot R_{k+2}+R_{k+2}'),
\]
where $R_{k+1},R_{k+2},R_{k+2}'$ satisfy
\[
\Vert R_{k+1}\Vert_{L^2(\Omega_*)}\lesssim_M\epsilon^{-\frac{3}{2}}+\Vert\bar w_\epsilon\Vert_{H^{k+1}(\Omega)},\quad \Vert(R_{k+2},R_{k+2}')\Vert_{L^2(\Omega_*)}\lesssim_M\epsilon^{-\frac{5}{2}}+\Vert\bar w_\epsilon\Vert_{H^{k+2}(\Omega)}. 
\] 
By Lemma~\ref{lem:Deltabarellip} and interpolation we derive
\[
\Vert\divg_\xi\bar w_\epsilon\Vert_{H^{k+1}(\Omega_*)}\lesssim_M\epsilon^{-2}+\epsilon\Vert\bar w_\epsilon\Vert_{H^{k+2}(\Omega)}. 
\]
Similarly, we derive 
\[
\Vert\curl_\xi\bar w_\epsilon\Vert_{H^{k+1}(\Omega_*)}\lesssim_M\epsilon^{-2}+\epsilon\Vert\bar w_\epsilon\Vert_{H^{k+2}(\Omega)}. 
\]
Now using Lemma~\ref{lem:Hodgexi} we obtain \eqref{eq:wbarepsreg} and conclude \eqref{eq:ubarepsreg}. As a corollary, from $(\bar u_\epsilon-u_\epsilon)\circ X^{-1}=\epsilon^2\chi\overline\Delta(\tilde\chi\bar w_\epsilon)$ we see
\begin{equation}\label{eq:ubarepserr}
\Vert\bar u_\epsilon-u_\epsilon\Vert_{H^{k-r}(\Omega)}\lesssim_M\epsilon^r,\quad r\in[0,2].
\end{equation}

We next show

\begin{lemma}
We have
\begin{equation}\label{eq:z1barepstr}
\Vert\divg\bar u_\epsilon|_\Gamma\Vert_{H^{k-1-r}(\Gamma)}\lesssim_M\epsilon^r,\ r\in[-\frac{3}{2},2],\quad\Vert\divg\bar u_\epsilon|_\Gamma\Vert_{H^{k-\frac{3}{2}}(\Gamma)}\lesssim_MP\epsilon+\mathcal S^\frac{1}{2}\epsilon^\frac{3}{2}.
\end{equation}
Moreover, if $\bar z_{j,\epsilon}$ is related to $(h_\epsilon,\bar u_\epsilon)$, we have
\begin{equation}\label{eq:zjbarepstr}
\Vert\bar z_{j,\epsilon}\Vert_{H^{k-r-j}(\Gamma)}\lesssim_MP\epsilon^{\frac{1}{2}+r},\quad-\frac{3}{2}\leq r\leq\min\{k-j,\frac{3}{2}\},\ j=2,...,k+1.
\end{equation}
\end{lemma}
\proofheading{Proof}We first prove \eqref{eq:z1barepstr}. By spectral resolution of $\Delta_{\Gamma_*}$ and Moser's estimate, we have
\[\begin{aligned}
\Vert(I-\epsilon^2\overline\Delta)^{-1}\divg_\xi w_\epsilon|_{\Gamma_*}\Vert_{H^{k-1-r}(\Gamma_*)}&\lesssim\Vert\divg_\xi w_\epsilon|_{\Gamma_*}\Vert_{H^{k-1-r}(\Gamma_*)}\\
&\lesssim_M\Vert\divg u_\epsilon|_\Gamma\Vert_{H^{k-1-r}(\Gamma_*)}\lesssim_MP\epsilon^{r+\frac{1}{2}},\quad r\in[-\frac{3}{2},\frac{3}{2}]. 
\end{aligned}\]
By
\[\begin{aligned}
\left.(I-\epsilon^2\overline\Delta)^{-1}\divg_\xi w_\epsilon-\divg_\xi\bar w_\epsilon\,\right|_{\Gamma_*}=-\epsilon^2(I-\epsilon^2\Delta_{\Gamma_*})^{-1}[\xi^{ij}\partial_j,\overline\Delta]\bar w_{\epsilon i}|_{\Gamma_*},
\end{aligned}\]
%and the regularizing effect of $(I-\epsilon^2\Delta_{\Gamma_*})^{-1}$, 
it suffices to show 
\[
\Vert[\xi^{ij}\partial_j,\overline\Delta]\bar w_{\epsilon i}|_{\Gamma_*}\Vert_{H^{k-3}(\Gamma_*)}\lesssim_M1
\]
and
\[
\Vert[\xi^{ij}\partial_j,\overline\Delta]\bar w_{\epsilon i}|_{\Gamma_*}\Vert_{H^{k-\frac{3}{2}}(\Gamma_*)}\lesssim_M\epsilon^{-1}+\epsilon^{-\frac{1}{2}}\Vert\Gamma\Vert_{H^{k+\frac{3}{2}}}. 
\]
But these easily follow from \eqref{eq:wbarepsreg}. 

We next prove \eqref{eq:zjbarepstr}. The case $j=2j_0$ follows easily from \eqref{eq:ubarepserr}, \eqref{eq:ubarepsreg} since the disturbance occurs at lower-order terms. We now consider $j=2j_0+1$, $j_0=1,...,k_0$. As above, we have
\[
\Vert(I-\epsilon^2\overline\Delta)^{-1}z_{j,\epsilon}(X^{-1})\Vert_{H^{k-r-j}(\Gamma_*)}\lesssim_MP\epsilon^{r+\frac{1}{2}},\quad-\frac{3}{2}\leq r\leq\min\{k-j,\frac{3}{2}\}.
\]
Write $\bar z_{j,\epsilon}(X^{-1})=\Delta_\xi^{j_0}\divg_\xi\bar w_\epsilon+\bar N_{j,\epsilon}(X^{-1})$. By spectral resolution of $\Delta_{\Gamma_*}$, also by \eqref{eq:ubarepserr}, \eqref{eq:ubarepsreg}, we have
\[
\Vert N_{j,\epsilon}(X^{-1})-(I-\epsilon^2\overline\Delta)^{-1}N_{j,\epsilon}(X^{-1})\Vert_{H^{k-r-j}(\Gamma_*)}\lesssim_MP\epsilon^{r+\frac{1}{2}},
\]
\[
\Vert N_{j,\epsilon}(X^{-1})-\bar N_{j,\epsilon}(X^{-1})\Vert_{H^{k-r-j}(\Gamma_*)}\lesssim_MP\epsilon^{r+\frac{1}{2}},\quad-\frac{3}{2}\leq r\leq\min\{k-j,\frac{3}{2}\}.
\]
We turn to control $\Delta_\xi^{j_0}\divg_\xi w_\epsilon-\Delta_\xi^{j_0}\divg_\xi\bar w_\epsilon$. Like the proof of \eqref{eq:z1barepstr}, we only need to show
\begin{equation}\label{eq:commtr}
\Vert(I-\epsilon^2\overline\Delta)^{-1}[\Delta_\xi^{j_0}\divg_\xi,\overline\Delta]\bar w_\epsilon\Vert_{H^{k-r-j}(\Gamma_*)}\lesssim_MP\epsilon^{r-\frac{3}{2}},\quad-\frac{3}{2}\leq r\leq\min\{k-j,\frac{3}{2}\}.
\end{equation}
We first identify the structure of the commutator $[\Delta_\xi^{j_0}\divg_\xi,\overline\Delta]\bar w_\epsilon$. It has the form
\[
[\Delta_\xi^{j_0}\divg_\xi,\overline\Delta]\bar w_\epsilon=\overline\nabla\cdot F_{j+1}+F_{j+2},
\]
where $F_{j+1},F_{j+2}$ are multilinear expressions in $(\xi,\bar w_\epsilon)$ respectively involving $j+1,j+2$ derivatives. Moreover, in $F_{j+1}$ exactly $j$ derivatives fall on $\bar w_\epsilon$, while in $F_{j+2}$ at most $j$ derivatives fall on $\xi$ and at most $j$ derivatives fall on $\bar w_\epsilon$. To see this, we expand $[\Delta_\xi^{j_0}\divg_\xi,\overline\Delta]\bar w_\epsilon$ as a multilinear expression involving $j+2$ derivatives. Terms that cannot be attributed to $F_{j+2}$ must look like
\[
T_1=\Delta_\xi^{j_0}\overline\Delta\xi^{lm}\,\partial_m\bar w_{\epsilon l}\quad\text{or}\quad T_2=f(\xi,\partial\xi)\partial^{2j_0}\overline\nabla\partial\bar w_\epsilon. 
\]
From the harmonicity of $\xi^{lm}$: 
\[
\Delta_\xi\xi^{lm}=\Delta\partial_lX_m\circ X^{-1}=0
\]
we see $T_1=[\Delta_\xi^{j_0},\overline\Delta]\xi^{lm}\,\partial_m\bar w_{\epsilon l}$ is in fact of type $F_{j+2}$. In $T_2$, we commute $\overline\nabla$ with everything on the left, then we derive
\[
T_2=\overline\nabla\left(f(\xi,\partial\xi)\partial^{2j_0}\partial\bar w_\epsilon\right)+F_{j+2}=\overline\nabla\cdot F_{j+1}+F_{j+2}, 
\]
as claimed. Now using that
\[
\Vert(I-\epsilon^2\overline\Delta)^{-1}\overline\nabla\cdot F\Vert_{H^{s+1}(\Gamma_*)}\lesssim\epsilon^{-1}\Vert F\Vert_{H^s(\Gamma_*)},\quad s\geq0
\]
and using \eqref{eq:wbarepsreg}, we obtain \eqref{eq:commtr} and complete the proof of the lemma. \proofqed

We finally show

\begin{lemma}\label{lem:Emonobareps}
\[
\mathcal E(h_\epsilon,\bar u_\epsilon,\Omega)\leq\mathcal E(h_\epsilon,u_\epsilon,\Omega)+C(P)\epsilon. 
\]
\end{lemma}
\proofheading{Proof}The quantity $\Vert a_\epsilon^\frac{1}{2}\Delta^{k_0-1}\kappa\Vert_{L^2(\Gamma)}^2$ is unchanged. Using \eqref{eq:ubarepserr} for handling nonlinear terms, we only need to show 
\[
\Vert P_k(\partial)\bar u_\epsilon\Vert_{L^2(\Omega)}^2\leq\Vert P_k(\partial)u_\epsilon\Vert_{L^2(\Omega)}^2+C(P)\epsilon
\]
for $P_k(\partial)=\nabla\Delta^{k_0-1}\divg,\nabla^{k-1}\curl,L_\zeta^{k_0}$ and show 
\[
\Vert P_{k+2}(\partial)\bar u_\epsilon\Vert_{L^2(\Omega)}^2\leq\Vert P_{k+2}(\partial)u_\epsilon\Vert_{L^2(\Omega)}^2+C(P)\epsilon^{-3}
\]
for $P_{k+2}(\partial)=\nabla^{k+1}\curl$. With $Q_m(\partial)\,f(X^{-1})=P_m(\partial)f\circ X^{-1}$, these reduce to
\begin{equation}\label{eq:QDmono}
\Vert Q_m(\partial)\,\bar u_\epsilon(X^{-1})\Vert_{L^2(\Omega_*,J_Xdx)}^2\leq\Vert Q_m(\partial)w_\epsilon\Vert_{L^2(\Omega_*,J_Xdx)}^2+C(P)\epsilon^{1-2(m-k)},\quad m=k,k+2,
\end{equation}
where $J_X=\det(\nabla X^{-1})$. Note that $Q_m(\partial)=\sum_{1\leq|\alpha|\leq m} q_{m,\alpha}\partial^\alpha$ where $q_{m,\alpha}$ is a multilinear expression in $\xi,\zeta$ involving $m-|\alpha|$ derivatives in total. 

We first consider $m=k$. Thanks to \eqref{eq:ubarepserr} it is harmless to drop lower-order terms and assume $Q_k(\partial)=\sum_{|\alpha|=k}q_{k,\alpha}\partial^\alpha$. The polar identity reads as
\[
\Vert Q_k(\partial)\,\bar u_\epsilon(X^{-1})\Vert_{L^2(\Omega_*,J_Xdx)}^2=\Vert Q_k(\partial)w_\epsilon\Vert_{L^2(\Omega_*,J_Xdx)}^2+2I_k+\Vert Q_k(\partial)\left(\bar u_\epsilon(X^{-1})-w_\epsilon\right)\Vert_{L^2(\Omega_*,J_Xdx)}^2
\]
with
\[\begin{aligned}
I_k&=\int_{\Omega_*}Q_k(\partial)\left(\bar u_\epsilon(X^{-1})-w_\epsilon\right)\,Q_k(\partial)\,\bar u_\epsilon(X^{-1})\,J_Xdx\\
&=\epsilon^2\int_{\Omega_*}Q_k(\partial)[\chi\overline\Delta(\tilde\chi\bar w_\epsilon)]\,Q_k(\partial)\,\bar u_\epsilon(X^{-1})\,J_Xdx. 
\end{aligned}\]
Since $\tilde\chi=1$ on $\supp\chi$, it may be omitted from the integrand. Using
$\bar u_\epsilon(X^{-1})=\bar w_\epsilon-\epsilon^2(1-\chi)\overline\Delta\bar w_\epsilon$, which holds on $\supp\chi$, we can rewrite $I_k$ as
\[\begin{aligned}
I_k&=\epsilon^2\int Q_k(\partial)(\chi\overline\Delta\bar w_\epsilon)\,Q_k(\partial)\bar w_\epsilon\,J_Xdx\\ &\quad-\epsilon^4\int Q_k(\partial)(\chi\overline\Delta\bar w_\epsilon)\,Q_k(\partial)\left((1-\chi)\overline\Delta\bar w_\epsilon\right)\,J_Xdx\\ &:=I_k^1+I_k^2. 
\end{aligned}\]
Commuting $Q_k(\partial)$ with $\chi\overline\Delta$ in $I_k^1$ and integrating by parts, we derive
\[
I_k^1=-\epsilon^2\int|\overline\nabla Q_k(\partial)\bar w_\epsilon|^2\chi J_Xdx+O(\epsilon).
\]
Commuting $Q_k(\partial)$ with $\chi,1-\chi$ in $I_k^2$, we derive
\[
I_k^2=-\epsilon^4\int\left(Q_k(\partial)\overline\Delta\bar w_\epsilon\right)^2\chi(1-\chi)J_Xdx+O(\epsilon). 
\]
Now \eqref{eq:QDmono} for $m=k$ easily follows. 

We next consider $m=k+2$. As above, the key is to analyze
\[\begin{aligned}
I_{k+2}&=\int_{\Omega_*}Q_{k+2}(\partial)\left(\bar u_\epsilon(X^{-1})-w_\epsilon\right)\,Q_{k+2}(\partial)\,\bar u_\epsilon(X^{-1})\,J_Xdx\\
&=\int Q_{k+2}(\partial)(\epsilon^2\chi\overline\Delta\bar w_\epsilon)\,Q_{k+2}(\partial)\bar w_\epsilon\,J_Xdx\\
&\quad-\int Q_{k+2}(\partial)[\chi(\bar w_\epsilon-w_\epsilon)]\,Q_{k+2}(\partial)[(1-\chi)(\bar w_\epsilon-w_\epsilon)]\,J_Xdx:=I_{k+2}^1+I_{k+2}^2. 
\end{aligned}\]
The lower-order terms are controlled as above. Retain the principal part $Q_{k+2}(\partial)=\sum_{|\alpha|=k+2}q_{k+2,\alpha}\partial^\alpha$. Since $Q_{k+2}(\partial)w_\epsilon$ only differs from $\nabla^{k+1}\curl u_\epsilon\circ X^{-1}$ by at most $k+1$ derivatives on $w_\epsilon$, we still have
\[
\Vert Q_{k+2}(\partial)w_\epsilon\Vert_{L^2(\Omega_*)}\lesssim_MP\epsilon^{-2}. 
\]
This together with \eqref{eq:uepsreg}, \eqref{eq:wbarepsreg} implies
\[
I_{k+2}^2=-\int [Q_{k+2}(\partial)(\bar w_\epsilon-w_\epsilon)]^2\chi(1-\chi)J_Xdx+O(\epsilon^{-3}). 
\]
To control $I_{k+2}^1$, we commute $\epsilon^2\chi\overline\Delta$ with $Q_{k+2}(\partial)$. The commutator is
\[
\epsilon^2[Q_{k+2}(\partial),\chi\overline\Delta]\bar w_\epsilon=[Q_{k+2}(\partial),\chi](\bar w_\epsilon-w_\epsilon)+\epsilon^2\chi[Q_{k+2}(\partial),\overline\Delta]\bar w_\epsilon. 
\]
The first term, measured in $L^2(\Omega_*)$, is $O(\epsilon^{-1})$. The second term, according to \eqref{eq:Deltabarcomm}, has the form
\[
\epsilon^2\chi[Q_{k+2}(\partial),\overline\Delta]\bar w_\epsilon=\epsilon^2\left(\overline\nabla\cdot\tilde P_{k+2}(\partial)\bar w_\epsilon+\tilde Q_{k+2}(\partial)\bar w_\epsilon\right),
\]
where $\tilde P_{k+2}(\partial),\tilde Q_{k+2}(\partial)$ are differential operators of degree $k+2$ with $C^1$ coefficients. Thus, by \eqref{eq:wbarepsreg} and integration by parts, we have
\[
\int[Q_{k+2}(\partial),\epsilon^2\chi\overline\Delta]\bar w_\epsilon\,Q_{k+2}(\partial)\bar w_\epsilon\,J_Xdx\lesssim_MP\epsilon^{-3}+\sum_{j=0}^1\Vert\overline\nabla^jQ_{k+2}(\partial)\bar w_\epsilon\Vert_{L^2(\Omega_*,\chi J_Xdx)}. 
\]
On the other hand, Lemma~\ref{lem:Deltabarellip} together with $\overline\nabla\chi=0$ gives 
\[\begin{aligned}
\epsilon^2\int\chi\overline\Delta Q_{k+2}(\partial)\bar w_\epsilon\,Q_{k+2}(\partial)\bar w_\epsilon\,J_Xdx&=-\epsilon^2\int|\overline\nabla Q_{k+2}(\partial)\bar w_\epsilon|^2\chi J_Xdx\\
&\quad-\epsilon^2\int\overline\nabla\left(\frac{J_X}{J_\gamma}\right)\cdot\overline\nabla Q_{k+2}(\partial)\bar w_\epsilon\,Q_{k+2}(\partial)\bar w_\epsilon\,\chi J_\gamma dx. 
\end{aligned}\]
Combining these and using Young's inequality we obtain \eqref{eq:QDmono} for $m=k+2$. This completes the proof of the lemma. \proofqed

\subsection{Recovery of compatibility}\label{sub:compat_recovery}

\begin{proposition}\label{prop:Ej}
For an integer $j\geq0$ and a function $f$ on $\Gamma_*$, define $v=E_jf$ via
\[
v(\gamma_s)=\chi\,\frac{s^j}{j!}\exp\left(-s^{2K}(1-\Delta_{\Gamma_*})^K\right)f,
\]
where $\gamma_s=\gamma(s,\cdot)$, $\chi$ is a cut-off function supported near $\Gamma_*$, and $K$ is a sufficiently large integer. We have
\[
\nabla_\nu^jv=f,\quad\nabla_\nu^iv=0,\ i<j
\]
on $\Gamma_*$, where $\nu=\partial_s\gamma\circ\gamma^{-1}$. Moreover, we have
\[
\Vert E_jf\Vert_{H^{r+j+\frac{1}{2}}(\Omega_*)}\lesssim\Vert f\Vert_{H^r(\Gamma_*)},\quad r>0. 
\]
\end{proposition}
The extension estimate is proved in Appendix~\ref{app:extension_lifting}. We now use it to correct the compatibility traces while retaining the energy bounds. The identity $h_\epsilon|_\Gamma=0$ already holds. To cancel the trace of $\divg\bar u_\epsilon$, set
\[
u^{(1)}=\bar u_\epsilon-\nabla\Delta_\Omega^{-1}\mathcal H_\Gamma\divg\bar u_\epsilon. 
\]
For consistent notation, set $h^{(1)}=h_\epsilon$. Estimate \eqref{eq:z1barepstr} gives
\begin{equation}\label{eq:u(1)regerr}
\Vert u^{(1)}\Vert_{H^{k+r}(\Omega)}\lesssim_M\epsilon^{-r},\quad\Vert u^{(1)}-\bar u_\epsilon\Vert_{H^{k-r}(\Omega)}\lesssim_M\epsilon^r,\quad r\in[0,2]. 
\end{equation}
Moreover, since the leading term in $\bar z_{j,\epsilon}$ is unchanged, we have
\begin{equation}\label{eq:zj(1)tr}
\Vert z_j^{(1)}|_\Gamma\Vert_{H^{k-j-r}(\Gamma)}\lesssim_MP\epsilon^{r+\frac{1}{2}},\quad-\frac{3}{2}\leq r\leq\min\{k-j,\frac{3}{2}\},\ j=2,...,k+1. 
\end{equation}
Moreover, by \eqref{eq:Emonoeps}, Lemma~\ref{lem:Emonobareps} and Young's inequality, we have
\[
\mathcal E(h^{(1)},u^{(1)},\Omega)+C_1\epsilon^2\mathcal D\leq\mathcal E(h_0,u_0,\Omega_0)+C(P)\epsilon. 
\]

The higher-order corrections use the linearization of the compatibility expressions. For a state $(h,u)$, let
\[
B_j[h,u;\partial](g,w)=\left.\frac{\partial}{\partial t}\right|_{t=0}z_j(h+tg,u+tw)=\begin{cases}\Delta^ig+\cdots&j=2i\\\Delta^i\divg w+\cdots&j=2i+1\end{cases}
\]
be the linearization of $z_j(h,u)$. Set $g(X^{-1})=G$, $w(X^{-1})=f\nu$, then
\[
B_j[h,u;\partial](g,w)\circ X^{-1}=B_j\left[h(X^{-1}),u(X^{-1});\xi\partial\right](G,f\nu). 
\]
We isolate the highest normal derivatives in these expressions. For $j=2i$,
\[
B_{2i}\left[h(X^{-1}),u(X^{-1});\xi\partial\right](G,f\nu)=\Delta_\xi^iG+\cdots,
\]
where $\cdots$ involves at most $2i-1$ normal derivatives in $(G,f)$. From $\Delta_\xi=\xi_{il}\xi_{im}\partial_l\partial_m+\xi_{il}\partial_l\xi_{im}\partial_m$ and the positivity of $(\xi_{il}\xi_{im})_{lm}$ we see
\[
B_{2i}\left[h(X^{-1}),u(X^{-1});\xi\partial\right](G,f\nu)=c_{2i}\nabla_\nu^{2i}G+\cdots,
\]
where $c_{2i}$ is a polynomial in $\xi$ nonvanishing everywhere on $\Gamma_*$. Similarly, for $j=2i+1$ we have
\[\begin{aligned}
B_{2i+1}\left[h(X^{-1}),u(X^{-1});\xi\partial\right](G,f\nu)&=\Delta_\xi^i\divg_\xi(f\nu)+\cdots\\
&=\Delta_\xi^i(\xi_{lm}\nu_l\nu_m\,\nabla_\nu f+\cdots)+\cdots=c_{2i+1}\nabla_\nu^{2i+1}f+\cdots,
\end{aligned}\]
where $c_{2i+1}$ has the same properties as $c_{2i}$, and $\cdots$ involves at most $2i$ normal derivatives in $(G,f)$. 

With these preparations, we can iteratively construct $(h^{(j)},u^{(j)})=(h^{(j-1)}+g^{(j)},u^{(j-1)}+w^{(j)})$ satisfying
\[
B_l\left[h^{(j-1)},u^{(j-1)};\partial\right](g^{(j)},w^{(j)})=-z_l^{(j-1)},\quad l=0,1,...,j,\quad j=2,3,...,k-2. 
\] 
Choose the corrections as follows:
\[
g^{(2i)}(X^{-1})=E_{2i}\left(-c_{2i}^{-1}z_{2i}^{(2i-1)}\right),\quad w^{(2i)}=0,\quad i=1,...,k_0-1,
\]
\[
g^{(2i+1)}=0,\quad w^{(2i+1)}(X^{-1})=E_{2i+1}\left(-c_{2i+1}^{-1}z_{2i+1}^{(2i)}\right)\nu,\quad i=1,...,k_0-2. 
\]
The remainder
\[
z_l^{(j)}-z_l^{(j-1)}-B_l\left[h^{(j-1)},u^{(j-1)};\partial\right](g^{(j)},w^{(j)})
\]
is quadratic in $(g^{(j)},w^{(j)})$, $l=2,...,j$. Starting from \eqref{eq:zj(1)tr} and using Proposition~\ref{prop:Ej}, by induction we obtain
\[
\Vert(g^{(j)},w^{(j)})\Vert_{H^{k-r}(\Omega)}\lesssim_MP\epsilon^{1+r},\quad r\in[-2,1],
\]
$(h^{(j)},\divg u^{(j)})|_\Gamma=(0,0)$ and 
\[
\Vert z_l^{(j)}|_\Gamma\Vert_{H^{k-\frac{1}{2}-j-r}(\Gamma)}\lesssim_MP\epsilon^{2+r},\ r\in[-2,1],\ l=2,...,j. 
\]

For $2\leq j\leq k-2$, the traces $z_j^{(k-2)}|_\Gamma$ are now smaller than the errors allowed in \eqref{eq:zj0tr0}. It remains to correct $z_{k-1}^{(k-2)}$ and $z_k^{(k-2)}$. Use the operators $I_m,J_m$ from Section~\ref{sub:regop}, now defined on $\Omega$; the analogue of Lemma~\ref{lem:ImJm} holds with constants depending on $M$. Write $I_\epsilon,J_\epsilon$ for these operators at $m=-\log_2\epsilon$, and define
\[
h=h^{(k-2)}-J_\epsilon^{k_0}I_\epsilon\left(z_k^{(k-2)}|_\Gamma\right),\quad u=u^{(k-2)}-\nabla J_\epsilon^{k_0}I_\epsilon\left(z_{k-1}^{(k-2)}|_\Gamma\right). 
\]
Then we have
\[
\Delta^ih=\Delta^{i-1}\divg u=0\text{ on }\Gamma,\quad i\leq k_0-1
\]
and
\[
\Delta^{k_0}h|_\Gamma=\Delta^{k_0}h^{(k-2)}-z_k^{(k-2)},\quad \Delta^{k_0-1}\divg u|_\Gamma=\Delta^{k_0-1}\divg u^{(k-2)}-z_{k-1}^{(k-2)}. 
\]
Moreover, as mentioned above, by (an analogue of) Lemma~\ref{lem:ImJm}, we have
\[
\Vert(h-h^{(k-2)},u-u^{(k-2)})\Vert_{H^{k-r}(\Omega)}\lesssim_MP\epsilon^{1+r},\quad r\in[-2,1]. 
\]
Combining these estimates, we have
\begin{equation}\label{eq:ulterr}
\Vert(h-h^{(1)},u-u^{(1)})\Vert_{H^{k-r}(\Omega)}\lesssim_MP\epsilon^{1+r},\quad r\in[-2,1],
\end{equation}
\begin{equation}\label{eq:zjtr}
(h,\divg u)|_\Gamma=(0,0),\quad\Vert z_j|_\Gamma\Vert_{H^{k-j-r}(\Gamma)}\lesssim_MP\epsilon^{r+\frac{3}{2}},\ -\frac{3}{2}\leq r\leq\min\{k-j,\frac{3}{2}\},\ j=2,...,k,
\end{equation}
and
\begin{equation}\label{eq:zjextra}
\Vert z_{k+1}|_\Gamma\Vert_{H^{\frac{1}{2}-r}(\Gamma)}\lesssim_MP\epsilon^{-1+r},\quad r\in[0,\frac{1}{2}]. 
\end{equation}
Two consequences of \eqref{eq:ulterr} are
\begin{equation}\label{eq:overallreg}
\Vert(h,u,\Omega)\Vert_{\mathbf H^{k+r}}\lesssim_M\epsilon^{-r},\quad r\in[0,2],
\end{equation}
\begin{equation}\label{eq:toterrinreg}
\Vert(h-h_0,u-u_0)\Vert_{H^{k-r}(\Omega)}+\Vert\eta-\eta_0\Vert_{H^{k-r}(\Gamma_*)}\lesssim_M\epsilon^r,\quad r\in[0,2].
\end{equation}

We finish by recording a lower bound for the dissipation $\mathcal D$.

\begin{lemma}\label{lem:Edissip}
We have
\[
\mathcal D':=\int|\nabla\Delta^{k_0}h|^2+(\Delta^{k_0}\divg u)^2\,dx\leq\mathcal D+C(P)\epsilon^{-1}.
\]
As a consequence, we have
\begin{equation}\label{eq:Edissip}
\mathcal E(h,u,\Omega)+C_1\epsilon^2\mathcal D'\leq\mathcal E(h_0,u_0,\Omega_0)+C(P)\epsilon. 
\end{equation}
\end{lemma}
\proofheading{Proof}By \eqref{eq:ulterr} and $\Delta^{k_0}\divg u^{(1)}=\Delta^{k_0}\divg\bar u_\epsilon$ we have
\[
\mathcal D'\leq\int|\nabla\Delta^{k_0}h^{(1)}|^2+(\Delta^{k_0}\divg u^{(1)})^2\,dx+C(P)\epsilon^{-1}=\int|\nabla\Delta^{k_0}h_\epsilon|^2+(\Delta^{k_0}\divg\bar u_\epsilon)^2\,dx+C(P)\epsilon^{-1}. 
\]
It remains to prove $\Vert\Delta^{k_0}\divg\bar u_\epsilon\Vert_{L^2(\Omega)}^2\leq\Vert\Delta^{k_0}\divg u_\epsilon\Vert_{L^2(\Omega)}^2+C(P)\epsilon^{-1}$. The argument is the same commutator and integration-by-parts estimate used in Lemma~\ref{lem:Emonobareps}. \proofqed

\subsection{Transport and the one-step energy bound}\label{sub:transport_step}

Starting from the regularized, corrected state, we advance the domain and fluid variables by one Euler step. Define
\begin{equation}\label{eq:x1def}
\Omega_1=x_1(\Omega),\quad x_1=x+\epsilon u
\end{equation}
and
\begin{equation}\label{eq:expsch}
\left\{\begin{aligned}
&h(x_1)=h-\epsilon\divg u,\\
&u(x_1)=u-\epsilon\nabla h.
\end{aligned}\right.
\end{equation}
With $A=(\nabla x_1)^{-*}$ we have $\partial_if\circ x_1=A_{ij}\partial_jf(x_1)$ and $n_1(x_1)=An/|An|$. Note the formula
\begin{equation}\label{eq:Aexpr}
A_{ij}=\delta_{ij}-\epsilon A_{ik}\partial_ku_j=\delta_{ij}-\epsilon\partial_iu_j+\epsilon^2A_{il}\partial_lu_k\partial_ku_j. 
\end{equation}

\begin{lemma}\label{lem:alg}
(i) With $\partial_i^A=A_{ij}\partial_j$ we have
\begin{equation}\label{eq:voralg}
\curl u_{1,ij}\circ x_1=\curl u_{ij}-\epsilon(\partial_i^Au_l\partial_lu_j-\partial_j^Au_l\partial_lu_i)+\epsilon^2(\partial_i^Au_l\partial_l\partial_jh-\partial_j^Au_l\partial_l\partial_ih). 
\end{equation}
(ii) For $j=1,...,k$ we have
\begin{equation}\label{eq:wavealg}
\left\{\begin{aligned}
&z_{j,1}(x_1)=z_j-\epsilon z_{j+1}+\epsilon^2f_{j+1},\\
&w_{j,1}(x_1)=w_j-\epsilon w_{j+1}+\epsilon^2g_{j+1},
\end{aligned}\right.
\end{equation}
where $f_{j+1},g_{j+1}$ are multilinear expressions in $h,u,A$ involving $j+2$ derivatives among which at most $j+1$ derivatives fall on $h,u$ and at most $j$ derivatives fall on $A$. \\
(iii) For $j=1,...,k_0$ we have
\begin{equation}\label{eq:kappaalg}
\Delta_{\Gamma_1}^{j-1}\kappa_1\circ x_1=\Delta_\Gamma^{j-1}\kappa-\epsilon\Delta_\Gamma^ju\cdot n+\epsilon p_{2j}+\epsilon^2q_{2j+\frac{1}{2}}. 
\end{equation}
Here $p_{2j},q_{2j+\frac{1}{2}}$ are respectively quasilinear expressions in $S=(\nabla u,A,\nabla^\top n)$, $T=(n,\nabla u,A)$ involving $2j-2$, $2j-1$ tangential derivatives on $S,T$ whose coefficients are smooth functions of $T$. 
\end{lemma}
\proofheading{Proof}(i) follows from direct calculation. (ii) is a special case of Lemma~\ref{lem:alg2}. \\
(iii) From \eqref{eq:Aexpr} we note that
\begin{equation}\label{eq:Anmodsq}
|An|^2=\sum_{i=1}^3(A_{ij}n_j)^2=1-2\epsilon\nabla_nu\cdot n+\epsilon^2q_\frac{3}{2}.
\end{equation}
Here $q_{j+1/2}$ denotes a generic expression of the type specified in the lemma. In particular, $q_{3/2}$ is a smooth function of $n,\nabla u,A$. Write
\[
n_1(x_1)=\frac{1-|An|^2}{|An|(1+|An|)}An+An. 
\]
Expand the right-hand side using \eqref{eq:Aexpr}, \eqref{eq:Anmodsq}. Also from \eqref{eq:Anmodsq} we note that $\frac{1}{|An|(1+|An|)}=\frac{1}{2}+\epsilon q_\frac{3}{2}$. It follows that
\begin{equation}\label{eq:n1(x1)}
n_1(x_1)=n-\epsilon\nabla u\cdot n+\epsilon(\nabla_nu\cdot n)n+\epsilon^2q_\frac{3}{2}=n-\epsilon\nabla^\top u\cdot n+\epsilon^2q_\frac{3}{2}. 
\end{equation}
Next, perform $\divg_\Gamma$ on both sides. For vector-valued functions $v,v_1$ on $\Gamma,\Gamma_1$ we have $\divg_\Gamma v=\zeta^{ij}\partial_jv_i$, $\divg_{\Gamma_1}v_1=\zeta_1^{ij}\partial_jv_{1i}$. For a function $f_1$ on $\Gamma_1$ and $f=f_1(x_1)$, by \eqref{eq:Aexpr} we have
\begin{equation}\label{eq:Dtopcomm}
\left(\zeta_1^{ij}\partial_jf_1\right)\circ x_1=\zeta^{ij}\partial_jf+\epsilon q_\frac{3}{2}\cdot\nabla^\top f. 
\end{equation}
These together with $\kappa_1=\divg_{\Gamma_1}n_1$, $\kappa=\divg_\Gamma n$ imply
\[
\kappa_1(x_1)=\kappa-\epsilon\Delta_\Gamma u\cdot n+\epsilon p_2+\epsilon^2q_\frac{5}{2},
\]
which is just \eqref{eq:kappaalg} for $j=1$. Finally, iteratively performing $\Delta_\Gamma=\zeta^{il}\partial_l(\zeta^{im}\partial_m)$ and using \eqref{eq:Dtopcomm}, we obtain \eqref{eq:kappaalg} for $j>1$. \proofqed

\begin{lemma}\label{lem:alg2}
Let $f_p(h,u)$ denote a multilinear expression in $h,u$ involving $p$ derivatives. Define $D_tf_p(h,u)$ via \eqref{eq:CE}, and define $\delta_tf_p(h,u)=\epsilon^{-1}(f_p(h_1,u_1)(x_1)-f_p(h,u))$. Note the difference between $f_p(h_1,u_1)(x_1)$ and $f_p(h_1(x_1),u_1(x_1))$. Then we have:\\
(i) $\delta_tf_p(h,u)$ is a multilinear expression in $h,u,A$ involving $p+1$ derivatives among which at most $p-1$ derivatives fall on $A$. \\
(ii) If $f_{p+q}(h,u)=f_p(h,u)f_q(h,u)$, then
\[
\delta_tf_{p+q}(h,u)=\delta_tf_p(h,u)\ f_q(h,u)+f_p(h,u)\ \delta_tf_q(h,u)+\epsilon\delta_tf_p(h,u)\ \delta_tf_q(h,u). 
\]
(iii) If $f_{p+1}(h,u)=\partial_if_p(h,u)$, then 
\[
\delta_tf_{p+1}(h,u)=\partial_i\delta_tf_p(h,u)-\partial_iu_j\partial_jf_p(h,u)-\epsilon A_{ik}\partial_ku_j\partial_j\delta_tf_p(h,u)+\epsilon A_{il}\partial_lu_k\partial_ku_jf_p(h,u).
\]
(iv) $\delta_tf_p(h,u)-D_tf_p(h,u)$ is $\epsilon$ times a multilinear expression in $h,u,A$ involving $p+2$ derivatives among which at most $p+1$ derivatives fall on $h,u$ and at most $p$ derivatives fall on $A$. 
\end{lemma}
\proofheading{Proof}(i) Suppose $f_p(h,u)=P_p(D)(h,u)$, where $P_p(D)$ is a differential operator of degree $p$. Then $f_p(h_1,u_1)(x_1)=P_p(AD)(h_1(x_1),u_1(x_1))$ and the conclusion is clear from \eqref{eq:expsch}, \eqref{eq:Aexpr}. The general case follows from (ii) and linear combination. \\ 
(ii) is the Leibniz rule for first-order difference quotient. \\
(iii) By \eqref{eq:Aexpr} we have
\[\begin{aligned}
f_{p+1}(h_1,u_1)&=\partial_i^Af_p(h_1,u_1)\\
&=A_{ij}\partial_j(f_p(h,u)+\epsilon\delta_tf_p(h,u))\\
&=\partial_if_p(h,u)-\epsilon\partial_iu_j\partial_jf_p(h,u)+\epsilon^2A_{il}\partial_lu_k\partial_ku_jf_p(h,u)\\
&\quad+\epsilon\partial_i\delta_tf_p(h,u)-\epsilon^2A_{ik}\partial_ku_j\partial_j\delta_tf_p(h,u). 
\end{aligned}\]
By a rearrangement of this identity the desired estimate follows. \\
(iv) Any $f_p(h,u)$ can be generated by $f_0(h,u)=(h,u)$ via linear combination and the operations analyzed in (ii), (iii). Then the result follows from an induction. \proofqed
\proofheading{Proof of Theorem~\ref{thm:1step}}Estimate \eqref{eq:toterrinreg} gives \eqref{eq:1stepappr}. The regularity bound \eqref{eq:overallreg} and the parameter relation \eqref{eq:numrel} imply (A3) for $(h_1,u_1,\Omega_1)$. By \eqref{eq:wavealg}, the regularity estimates above, and \eqref{eq:zjtr}, we have $h_1|_{\Gamma_1}=0$ and
\[\begin{aligned}
\Vert z_{j,1}|_{\Gamma_1}\Vert_{H^{k-j-r}(\Gamma_1)}&\lesssim_M\Vert z_{j,1}(x_1)|_\Gamma\Vert_{H^{k-j-r}(\Gamma)}\\
&\lesssim_MP\epsilon^{r+\frac{3}{2}}+\epsilon^{r+\frac{1}{2}},\quad\frac{1}{2}\leq r\leq\min\{k-j,\frac{3}{2}\},\quad j=1,...,k-1,
\end{aligned}\]
which combined with \eqref{eq:numrel} justifies (A2) for $(h_1,u_1)$. 

To prove the energy bound \eqref{eq:Emono}, consider the vorticity, acoustic, and tangential contributions in turn. First, \eqref{eq:voralg} gives
\[
\mathcal E^v(u_1)\leq\mathcal E^v(u)+C(M)\epsilon. 
\]
Next, we rewrite \eqref{eq:wavealg} as
\[
\left\{\begin{aligned}
&z_{j,1}(x_1)=z_j-\epsilon\divg w_j+\epsilon^2\tilde f_{j+1},\\
&w_{j,1}(x_1)=w_j-\epsilon\nabla z_j+\epsilon^2\tilde g_{j+1},
\end{aligned}\right.
\]
where $\tilde f_{j+1},\tilde g_{j+1}$ have the same form as $f_{j+1},g_{j+1}$. Integrate the square of both sides. By regularity estimate of $h,u$, \eqref{eq:zjtr}, \eqref{eq:zjextra} and divergence formula we have
\[\begin{aligned}
E^w(h_1,u_1)&\leq\sum_{j=0}^k\Vert(z_{j,1}(x_1),w_{j,1}(x_1))\Vert_{L^2(\Omega)}^2+C(M)\epsilon\\
&\leq E^w(h,u)+\epsilon^2\Vert(\divg w_k,\nabla z_k)\Vert_{L^2(\Omega)}^2+C(P)\epsilon. 
\end{aligned}\]
Since the leading terms in $\divg w_k,\nabla z_k$ are $\Delta^{k_0}\divg u,\nabla\Delta^{k_0}h$, we have
\[
\Vert(\divg w_k,\nabla z_k)\Vert_{L^2(\Omega)}^2\leq\mathcal D'+C(M)\epsilon^{-1}. 
\]

Finally, we control $E^b(h_1,u_1,\Omega_1)$. We show that
\begin{equation}\label{eq:zetadiff}
\Vert\zeta_1(x_1)-\zeta\Vert_{H^{k-1}(\Omega)}\lesssim_M\epsilon. 
\end{equation}
By $\Delta\zeta=0$, $\Delta_A\zeta_1(x_1)=\Delta\zeta_1\circ x_1=0$ we have
\[
\zeta_1(x_1)-\zeta=\mathcal H_\Gamma(\zeta_1(x_1)-\zeta)+\Delta_\Omega^{-1}(\Delta-\Delta_A)\zeta_1(x_1). 
\]
From the analysis in Lemma~\ref{lem:alg} (iii) we see $\Vert\zeta_1(x_1)-\zeta\Vert_{H^{k-\frac{3}{2}}(\Gamma)}\lesssim_M\epsilon$, and from \eqref{eq:Aexpr} we see $\Vert(\Delta-\Delta_A)\zeta_1(x_1)\Vert_{H^{k-\frac{5}{2}}(\Omega)}\lesssim_M\epsilon$. Thus, by the elliptic estimate we justify \eqref{eq:zetadiff}. With \eqref{eq:zetadiff}, \eqref{eq:Aexpr} we obtain 
\[
\Vert(L_{\zeta_1}^jh_1\circ x_1-L_\zeta^jh_1(x_1),\ L_{\zeta_1}^ju_1\circ x_1-L_\zeta^ju_1(x_1))\Vert_{L^2(\Omega)}\lesssim_M\epsilon,\quad j=1,...,k_0. 
\]
As a consequence, we have
\[
\left\{\begin{aligned}
&L_{\zeta_1}^jh_1\circ x_1=L_\zeta^jh-\epsilon\divg L_\zeta^ju+O_{L^2}(\epsilon),\\
&L_{\zeta_1}^ju_1\circ x_1=L_\zeta^ju-\epsilon\nabla L_\zeta^jh+O_{L^2}(\epsilon),\quad j=1,...,k_0. 
\end{aligned}\right.
\]
Integrate the square of both sides and use the divergence theorem. This yields
\[\begin{aligned}
\sum_{j=1}^{k_0}\Vert(L_{\zeta_1}^jh_1,L_{\zeta_1}^ju_1)\Vert_{L^2(\Omega_1)}^2\leq\sum_{j=1}^{k_0}\Vert(L_\zeta^jh,L_\zeta^ju)\Vert_{L^2(\Omega)}^2-2\epsilon\int_\Gamma L_\zeta^{k_0}h\,L_\zeta^{k_0}u\cdot n\,dS\\
+\epsilon^2\Vert(\divg L_\zeta^{k_0}u,\nabla L_\zeta^{k_0}h)\Vert_{L^2(\Omega)}^2+C(M)\epsilon. 
\end{aligned}\]
The square integral in the second line can be controlled by
\[
\Vert(\divg L_\zeta^{k_0}u,\nabla L_\zeta^{k_0}h)\Vert_{L^2(\Omega)}^2\lesssim_M\mathcal D'+\epsilon^{-1}. 
\]
This can be seen by commuting derivatives, noting the decomposition
\[
\divg u=\Delta^{-k_0}\Delta^{k_0}\divg u+\sum_{k=1}^{k_0-1}\Delta^{-j}\mathcal H\Delta^j\divg u,\quad h=\Delta^{-k_0}\Delta^{k_0}h+\sum_{k=1}^{k_0-1}\Delta^{-j}\mathcal H\Delta^jh,
\]
the elliptic estimate, and \eqref{eq:zjtr}. The boundary integral in the first line cancels with the increment of the curvature energy. To see this, square \eqref{eq:kappaalg} and integrate with respect to $a\,dS$. The structure of $p_{2j},q_{2j+1/2}$ in Lemma~\ref{lem:alg} and the bound $|a_1(x_1)-a|\lesssim_M\epsilon$ give
\[
\sum_{j=1}^{k_0}\Vert a_1^\frac{1}{2}\Delta_{\Gamma_1}^{j-1}\kappa_1\Vert_{L^2(\Gamma_1)}^2\leq\sum_{j=1}^{k_0}\Vert a^\frac{1}{2}\Delta_\Gamma^{j-1}\kappa\Vert_{L^2(\Gamma)}^2-2\epsilon\int_\Gamma a\Delta_\Gamma^{k_0-1}\kappa\,\Delta_\Gamma^{k_0}u\cdot n\,dS+C(M)\epsilon. 
\]
As shown in Section~\ref{sub:propa}, we have
\[
\Vert L_\zeta^{k_0}h\Vert_{L^2(\Gamma)}+\Vert L_\zeta^{k_0}h+a\Delta_\Gamma^{k_0-1}\kappa\Vert_{H^\frac{1}{2}(\Gamma)}\lesssim_M1,
\] 
\[
\Vert L_\zeta^{k_0}u\Vert_{H^{-\frac{1}{2}}(\Gamma)}+\Vert L_\zeta^{k_0}u-\Delta_\Gamma^{k_0}u\Vert_{L^2(\Gamma)}\lesssim_M1. 
\]
Now combining everything and selecting $C_1=C_1(M)$ sufficiently large we complete the proof of the theorem. \proofqed

\subsection{Convergence of the iteration scheme}\label{sub:cvgitersch}

We now pass from the discrete scheme to an exact solution.

\begin{theorem}[Existence of regular solutions]\label{thm:regsol}
Suppose $(h_0,u_0,\Omega_0)\in\mathbf H^k$ with $\Gamma_0\in\Lambda_*$ satisfies the Taylor condition and the compatibility conditions: 
\[
z_{j,0}|_{\Gamma_0}=0,\quad j=0,1,...,k-1. 
\]
Then there is a $T>0$, depending on $\Vert(h_0,u_0,\Omega_0)\Vert_{\mathbf H^k}$, and a solution $(h,u,\Omega)$ to \eqref{eq:CE} on $[0,T]$ which is uniformly bounded in $\mathbf H^k$. 
\end{theorem}

Define $(h_0^m,u_0^m,\Omega_0^m)=\Psi_m(h_0,u_0,\Omega_0)$. By Proposition~\ref{prop:Psim} there is an $A_0>0$, depending only on $\Vert\Gamma_0\Vert_{C^{1,\delta}}+\Vert(h_0,u_0)\Vert_{C^{1/2+\delta}(\Omega_0)}$, satisfying \eqref{eq:ptass0} for $(h_0^m,u_0^m,\Omega_0^m)$, and there are $M>0$, $N\sim_{A_0}M$, depending only on $\Vert(h_0,u_0,\Omega_0)\Vert_{\mathbf H^k}$, satisfying (A1)-(A3) for $(h_0^m,u_0^m,\Omega_0^m)$ with $\epsilon=2^{-m}$ in Theorem~\ref{thm:1step}. The initial regularization preserves the Taylor sign condition. Now let $(h^m(j\epsilon),u^m(j\epsilon),\Omega^m_{j\epsilon})$ be the $j$th iteration of the scheme constructed in Theorem~\ref{thm:1step}. By \eqref{eq:overallreg} and \eqref{eq:Emono}, induction shows that $\Gamma^m_{j\epsilon}$ remains in a uniform (but possibly fattened) collar neighborhood, that the Taylor coefficient has a uniform lower bound, and that
\begin{equation}\label{eq:unibdofdts}
\Vert(h^m(j\epsilon),u^m(j\epsilon),\Omega^m_{j\epsilon})\Vert_{\mathbf H^k}\leq C(A_0)M,\quad j\epsilon\leq\frac{1}{C(M)}:=T. 
\end{equation}

We extract a subsequential limit of the discrete approximations. Sobolev embedding gives
\begin{equation}\label{eq:huetaC2x}
\Vert(h^m,u^m)\Vert_{C^2(\Omega^m)}\lesssim_A1,\quad \Vert\eta_{\Gamma^m}\Vert_{C^2(\Gamma_*)}\lesssim_M1. 
\end{equation}
From the third relation in \eqref{eq:1stepappr} we see
\begin{equation}\label{eq:etamLiptC1x}
\Vert\eta_{\Gamma^m_t}-\eta_{\Gamma^m_s}\Vert_{C^1(\Gamma_*)}\lesssim_M|t-s|,\quad t,s\in\mathbf T_m:=2^{-m}\mathbb N\cap[0,T]. 
\end{equation}
A consequence of that is, for every $x\in\Omega^m_{(j+l)\epsilon}$ and $y\in\Omega^m_{j\epsilon}$ there exist $z_i\in\Omega^m_{(j+i-1)\epsilon}\cap\Omega^m_{(j+i)\epsilon}$, $i=1,...,l$ with $|y-z_1|+|z_1-z_2|+\cdots+|z_l-x|\lesssim_M|t-s|+|x-y|$. Then \eqref{eq:1stepappr} and \eqref{eq:huetaC2x} imply
\begin{equation}\label{eq:hDhmLiptx}
\sum_{r=0}^1|\nabla^rh^m(t,x)-\nabla^rh^m(s,y)|\lesssim_M|x-y|+|t-s|,\quad t,s\in\mathbf T_m,\ x\in\Omega^m_t,\ y\in \Omega^m_s. 
\end{equation}
The estimate of $u^m$ is similar. By \eqref{eq:etamLiptC1x}, \eqref{eq:hDhmLiptx} and its counterpart for $u^m$, an iteration of \eqref{eq:1stepappr} implies
\begin{equation}\label{eq:multistepappr}
\left\{\begin{aligned}
&h^m(t)-h^m(s)+(t-s)(u^m(s)\cdot h^m(s)+\divg u^m(s))=O(|t-s|^2)&\text{in }\Omega^m_{[s,t]}\\
&u^m(t)-u^m(s)+(t-s)(u^m(s)\cdot u^m(s)+\nabla h^m(s))=O(|t-s|^2)&\text{in }\Omega^m_{[s,t]}\\
&\Omega^m_t=(\text{Id}+(t-s)u^m(s))\Omega^m_s+O(|t-s|^2)
\end{aligned}\right.
\end{equation}
for $t,s\in\mathbf T_m$. Here $\Omega^m_{[s,t]}$ is the intersection of $\Omega^m_\tau$ over intermediate $\tau$ (including $s,t$). 

The set $\mathbf T=\bigcup_m\mathbf T_m$ is countable and dense in $[0,T]$. Compactness of the embedding $C^2(\Gamma_*)\subset\subset C^1(\Gamma_*)$ and a diagonal argument give a subsequence of $\eta_{\Gamma^m}$ converging in $C^1(\Gamma_*)$ at every $t\in\mathbf T$ to a limit $\eta_t$. Passing to the limit in \eqref{eq:etamLiptC1x} yields
\begin{equation}\label{eq:etaLiptC1x}
\Vert\eta_t-\eta_s\Vert_{C^1(\Gamma_*)}\lesssim_M|t-s|,\quad t,s\in\mathbf T.
\end{equation}
Hence $t\mapsto\eta_t$ extends to a Lipschitz map from $[0,T]$ to $C^1(\Gamma_*)$, with the same Lipschitz constant. The convergence is uniform on the discrete time sets: $\sup_{t\in\mathbf T_m}\Vert\eta_{\Gamma_t^m}-\eta_t\Vert_{C^1(\Gamma_*)}\to0$. Let $\Gamma_t$ be the surface parameterized by $\eta_t$. The uniform collar bounds and lower semicontinuity of the $C^{1,\alpha}$ norm keep the limiting surfaces in the fixed collar.

Thanks to the convergence of domain, the sequence $h^m(t)$ can be eventually restricted on $\Omega'\subset\subset\Omega_t$, $t\in\mathbf T$. By compact embedding $C^{0,1}(\Omega')\subset\subset C^0(\Omega')$ and diagonal argument, there is a subsequence of $h^m$ that at every $t\in\mathbf T$ converges uniformly on every compact subset of $\Omega_t$ to a limit $h(t)$. Passing to the limit in \eqref{eq:hDhmLiptx} we see $h$ can be extended to a Lipschitz function on the whole spacetime domain $\bigcup_{\tau\in[0,T]}\{\tau\}\times\Omega_\tau$. Using \eqref{eq:etaLiptC1x}, \eqref{eq:hDhmLiptx} and the Lipschitz bound on $h$ we can strengthen the convergence of $h^m$ as: $\sup_{t\in\mathbf T_m}\Vert h^m(t)-h(t)\Vert_{C^0(\Omega^m_t\cap\Omega_t)}\to0$. The same conclusions hold for $u^m$, $\nabla h^m$, $\nabla u^m$. Moreover, an application of the fundamental theorem of calculus indicates that the limit of $\nabla(h^m,u^m)$ is just the spatial gradient of $(h,u)$. 

By \eqref{eq:etamLiptC1x}, \eqref{eq:etaLiptC1x} and the uniform convergence of $\eta_{\Gamma^m_t}$, an interior point in $\Omega_s$ that is away from $\Gamma_s$ by a multiple of $|t-s|$ must lie in the intersection of $\Omega^m_{[s,t]}$, $m\gg1$. Thus, dividing the first two equations in \eqref{eq:multistepappr} by $t-s$, letting $m\to\infty$ and then $t\to s$ we obtain
\[
\left\{\begin{aligned}
&\partial_th+u\cdot\nabla h+\divg u=0\quad\text{in }\Omega_t,\\
&\partial_tu+u\cdot\nabla u+\nabla h=0\quad\text{in }\Omega_t.
\end{aligned}\right.
\]
The dynamic boundary condition \eqref{eq:vacbc} follows from that for the approximate state, the convergence of $h^m,\Gamma^m$, and the uniform Lipschitz bound on $h,h^m$. Likewise, the Taylor condition follows from the uniform lower bound on Taylor coefficient for the approximate state, the convergence of $\nabla h^m,n_{\Gamma^m}$, and the uniform Lipschitz bound on $\nabla h,\nabla h^m$. Letting $m\to\infty$ in the third equation in \eqref{eq:multistepappr} we deduce 
\[
\Gamma_t=(\text{Id}+(t-s)u(s))\Gamma_s+O(|t-s|^2). 
\]
It follows that for every $\phi\in C^1(\mathbb R^{1+3})$ vanishing on $\bigcup\{t\}\times\Gamma_t$ and every $y\in\Gamma_s$ we have $\phi(t,y+(t-s)u(s,y))=O(|t-s|^2)$, and hence $(\partial_t+u\cdot\nabla)\phi=0$ on $\bigcup\{t\}\times\Gamma_t$. This justifies the kinematic boundary condition \eqref{eq:kinbc}. Therefore, $(h,u,\Omega_t)$ is a classical solution to the free boundary Euler system \eqref{eq:kinbc}--\eqref{eq:vacbc} preserving the Taylor sign \eqref{eq:TS}. Finally, by Lemma~\ref{lem:lsc} we justify the uniform $\mathbf H^k$ bound for $(h,u,\Omega_t)$. This completes the proof of Theorem~\ref{thm:regsol}. 

\begin{lemma}\label{lem:lsc}
Suppose $\Gamma^m\in H^k$ and $f^m\in H^k(\Omega^m)$ are uniformly bounded. Suppose $\Gamma\in\Lambda_*$, $f\in C^{0,1}(\Omega)$ and 
\[
\lim_{m\to\infty}\left(\Vert\eta_{\Gamma^m}-\eta_\Gamma\Vert_{C^0(\Gamma_*)}+\Vert f^m-f\Vert_{C^0(\Omega^m\cap\Omega)}\right)=0. 
\]
Then we have $\Gamma\in H^k$, $f\in H^k(\Omega)$ and
\[
\Vert\Gamma\Vert_{H^k}\leq\liminf_{m\to\infty}\Vert\Gamma^m\Vert_{H^k},\quad\Vert f\Vert_{H^k(\Omega)}\lesssim\liminf_{m\to\infty}\Vert f^m\Vert_{H^k(\Omega^m)}. 
\]
\end{lemma}
\par\smallskip\noindent\textbf{Remark.}\ The second bound includes a uniform extension constant. This weaker form of lower semicontinuity suffices here.
\proofheading{Proof}The first statement, equivalent to the lower semicontinuity of $H^k(\Gamma_*)$ norm, is standard. To see this, recall that $\Vert\eta_{\Gamma^m}\Vert_{H^k(\Gamma_*)}^2=\sum_i\Vert(\alpha_{*i}\eta_{\Gamma^m})\circ H_i(\cdot,0)\Vert_{H^k(\mathbb R^2)}^2$. The uniform convergence of $\eta_{\Gamma^m}$ implies the $L^2$ convergence of $(\alpha_{*i}\eta_{\Gamma^m})\circ H_i(\cdot,0)$, which in turn implies $\Vert(\alpha_{*i}\eta_\Gamma)\circ H_i(\cdot,0)\Vert_{H^k(\mathbb R^2)}^2\leq\liminf\Vert(\alpha_{*i}\eta_{\Gamma^m})\circ H_i(\cdot,0)\Vert_{H^k(\mathbb R^2)}^2$. Summing up we thus justify the claim. 

To prove the second statement, we select for each $m$ an $F^m\in H^k(\mathbb R^3)$ extending $f^m$ with equivalent $H^k$ norm. Select a subsequence of $\Vert f^m\Vert_{H^k(\Omega^m)}$ that realizes the limit inferior, select a $\chi\in C_c^\infty(\mathbb R^3)$ with $\chi|_\Omega=1$, and further select a subsequence of $\chi F^m$ that converges in $L^2(\mathbb R^3)$ to a limit $F$. It follows that $F\in H^k(\mathbb R^3)$ with
\[
\Vert F\Vert_{H^k(\mathbb R^3)}\lesssim\liminf_{m\to\infty}\Vert f^m\Vert_{H^k(\Omega^m)}. 
\]
For every $g\in C_c(\Omega)$ it holds
\[
\int_\Omega Fg\,dx=\lim_{m\to\infty}\int_\Omega\chi F^mg\,dx=\lim_{m\to\infty}\int_\Omega f^mg\,dx=\int_\Omega fg\,dx. 
\]
Then we must have $F|_\Omega=f$. This completes the proof. \proofqed

\section{Rough solutions and continuous dependence}\label{sec:rough}

We prove Theorem~\ref{thm:LWP} by combining the regular solution theory with the distance and frequency-envelope estimates. Choose $\Gamma_*$ by smoothing $\Gamma_0$ at a scale determined by its $C^{1,\delta}$ geometry, and fix a corresponding collar.

For $m>m_0$, let $(h_0^m,u_0^m,\Omega_0^m)=\Psi_m(h_0,u_0,\Omega_0)$. Section~\ref{sec:regsol} gives a regular solution $(h^m,u^m,\Omega^m)$ with these initial data. We first obtain a common existence interval and uniform hybrid Sobolev bounds, then prove strong convergence and continuous dependence.

\subsection{Uniform estimates}

Set $M=\Vert(h_0,u_0,\Omega_0)\Vert_{\mathbf H^s_{1/2}}$ and let $(c_m)$ be the frequency envelope defined by \eqref{eq:feHs1/2}. Set $A^m(t)=\Vert\Gamma^m_t\Vert_{C^{1,\delta}}+\Vert u^m\Vert_{C^{\frac{1}{2}+\delta}(\Omega^m_t)}+\Vert h^m\Vert_{C^{1,\delta}(\Omega^m_t)}$ and let $A$ be that of $(h_0,u_0,\Omega_0)$ plus the term $\Vert\divg u_0\Vert_{C^\delta(\Omega_0)}$. For $m_0(M)\gg1$, for arbitrarily large $m_1$, and for $M_0\gg A_0=3A$ we set the bootstrap assumption: 
\begin{equation}\label{eq:btsass}
A^m(t)\leq A_0,\quad\Vert(h^m,u^m,\Omega^m_t)\Vert_{\mathbf H^s_{1/2}}\leq M_0,\quad m_0<m<m_1. 
\end{equation}
We close the bootstrap on a time interval $[0,T_0]$ independent of $m$ and $m_1$, with
\begin{equation}\label{eq:T0cond}
T_0\ll P(M_0)^{-1},
\end{equation}
where $P(M_0)$ is a polynomial of $M_0$. The dynamic control parameter in Theorems \ref{thm:apest}, \ref{thm:apinHk1} can be controlled by a polynomial of $\Vert(h^m,u^m,\Omega^m_t)\Vert_{\mathbf H^s_{1/2}}$. Thus, under the bootstrap assumption, the Theorems together with Proposition~\ref{prop:Psim2} imply 
\begin{equation}\label{eq:hfbd}
\Vert(h^m,u^m,\Omega^m)\Vert_{\mathbf H^k}\lesssim_A2^{m(k-s)}c_m,
\end{equation}
\begin{equation}\label{eq:hfbd1}
\Vert(h^m,u^m,\Omega^m)\Vert_{\mathbf H^k_1}\lesssim_A2^{m(k+\frac{1}{2}-s)}c_m. 
\end{equation}
On the other hand, Theorem~\ref{thm:L2dist} implies 
\begin{equation}\label{eq:lfdiff}
\Vert(h^m-h^{m+1},u^m-u^{m+1})\Vert_{L^2(\Omega^m\cap\Omega^{m+1})}+\Vert\eta^m-\eta^{m+1}\Vert_{L^2(\Gamma_*)}\lesssim_A2^{-ms}c_m. 
\end{equation}
Another consequence of Proposition~\ref{prop:Psim2} is
\[
\Vert(\divg(u^m-u^{m+1})|_{t=0},\ \nabla(h^m-h^{m+1})|_{t=0})\Vert_{L^2(\Omega_0^m\cap\Omega_0^{m+1})}\lesssim_A2^{-m(s-\frac{1}{2})}c_m. 
\]
Using the decompositions
\[
\nabla\divg(u^m-u^{m+1})=\nabla\Phi_m\divg(u^m-u^{m+1})+\nabla(I-\Phi_m)\divg u^m-\nabla(I-\Phi_m)\divg u^{m+1},
\]
\[
\Delta(h^m-h^{m+1})=\divg\Phi_m\nabla(h^m-h^{m+1})+\divg(I-\Phi_m)\nabla h^m-\divg(I-\Phi_m)\nabla h^{m+1},
\]
by Lemma~\ref{lem:Phim} and \eqref{eq:hfbd1}, we obtain
\[
\Vert(\Delta(h^m-h^{m+1})|_{t=0},\ \nabla\divg(u^m-u^{m+1})|_{t=0})\Vert_{L^2(\Omega_0^m\cap\Omega_0^{m+1})}\lesssim_A2^{-m(s-\frac{3}{2})}c_m. 
\]
Similarly, by \eqref{eq:hfbd} and \eqref{eq:lfdiff} we have
\[
\Vert(h^m-h^{m+1},u^m-u^{m+1})\Vert_{W^{1,2}(\Omega^m\cap\Omega^{m+1})}\lesssim_A2^{-m(s-1)}c_m.
\]
These together with $m>m_0(M)\gg1$ imply
\begin{equation}\label{eq:D2ini}
D_2(\mathbf v^m,\mathbf v^{m+1})|_{t=0}\lesssim_A2^{-m(2s-3)}c_m^2+2^{-m(2s-2)}M^2c_m^2\lesssim_A2^{-m(2s-3)}c_m^2. 
\end{equation}
Invoke Corollary~\ref{cor:PH2}. The dynamic control parameter therein can be dominated by a polynomial of $M_0$ due to \eqref{eq:btsass}, \eqref{eq:hfbd}, interpolation and Sobolev embedding. This implies
\begin{equation}\label{eq:A2vt}
\Vert(\nabla\divg(u^m-u^{m+1}),\ \nabla^2(h^m-h^{m+1}))\Vert_{L^2(\Omega^m\cap\Omega^{m+1})}\lesssim_A2^{-m(s-\frac{3}{2})}c_m. 
\end{equation}

We now derive the uniform hybrid Sobolev bound. For the interface, it suffices to prove
\begin{equation}\label{eq:etaHs}
\sum_j\Vert P_{*j}\eta^m\Vert_{H^s(\Gamma_*)}^2\lesssim_AM^2. 
\end{equation}
For $j\geq m$ we have
\[
\Vert P_{*j}\eta^m\Vert_{H^s(\Gamma_*)}\lesssim2^{-j(k-s)}\Vert\eta^m\Vert_{H^k(\Gamma_*)}\lesssim_A2^{-(j-m)(k-s)}c_m. 
\]
For $m_0<j<m$, we write
\[
\eta^m=\eta^j+\sum_{l=j}^{m-1}(\eta^{l+1}-\eta^l)
\]
and derive
\[\begin{aligned}
\Vert P_{*j}\eta^m\Vert_{H^s(\Gamma_*)}\lesssim2^{-j(k-s)}\Vert\eta^j\Vert_{H^k(\Gamma_*)}+\sum_{l=j}^{m-1}2^{js}\Vert\eta^{l+1}-\eta^l\Vert_{L^2(\Gamma_*)}\lesssim_Ac_j+\sum_{l=j}^{m-1}2^{-(l-j)s}c_l.
\end{aligned}\]
Summing in $j$ and using the convolution inequality, we obtain \eqref{eq:etaHs}. 

For $\mathbf v^m=(h^m,u^m)$, the same frequency argument must account for the varying domains. When $j\geq m$,
\[
\Vert P_j\mathbf v^m\Vert_{H^s(\Omega^m)}\lesssim_A2^{-j(k-s)}\Vert\mathbf v^m\Vert_{H^k(\Omega^m)}\lesssim_A2^{-(k-s)(j-m)}c_m
\]
thanks to \eqref{eq:hfbd}. For $m_0<j<m$, the decomposition 
\begin{equation}\label{eq:telesum}
\mathbf v^m=\Phi_j\mathbf v^j-\sum_{l=j}^{m-1}(\Phi_l\mathbf v^l-\Phi_{l+1}\mathbf v^{l+1})+(I-\Phi_m)\mathbf v^m
\end{equation}
makes sense on $\Omega^m$ due to the domain enlarging effect of $\Phi_l$. The term in the parenthesis is bounded in $L^2(\Omega^m)$ by 
\begin{align}
&\Vert(\Phi_l-\Phi_{l+1})\mathbf v^l\Vert_{L^2(\Omega^m)}+\Vert\Phi_{l+1}(\mathbf v^l-\mathbf v^{l+1})\Vert_{L^2(\Omega^m)}\nonumber\\
&\lesssim_A2^{-lk}\Vert\mathbf v^l\Vert_{H^k(\Omega^l)}+\Vert\mathbf v^l-\mathbf v^{l+1}\Vert_{L^2(\Omega^l\cap\Omega^{l+1})}\lesssim_A2^{-sl}c_l, \nonumber
\end{align}
where we have used Lemma~\ref{lem:Phim} and \eqref{eq:lfdiff}. It follows that
\[
\sum_{l=j}^{m-1}\Vert P_j(\Phi_l\mathbf v^l-\Phi_{l+1}\mathbf v^{l+1})\Vert_{H^s(\Omega^m)}\lesssim_A\sum_{l=j}^{m-1}2^{-s(l-j)}c_l. 
\]
For the other two terms in \eqref{eq:telesum} we have
\[
\Vert P_j\Phi_j\mathbf v^j\Vert_{H^s(\Omega^m)}\lesssim_A2^{-j(k-s)}\Vert\mathbf v^j\Vert_{H^k(\Omega^j)}\lesssim_Ac_j,
\]\[
\Vert\Vert P_j(I-\Phi_m)\mathbf v^m\Vert_{H^s(\Omega^m)}\Vert_{\ell^2_j}\lesssim_A\Vert(I-\Phi_m)\mathbf v^m\Vert_{H^s(\Omega^m)}\lesssim_A2^{-m(k-s)}\Vert\mathbf v^m\Vert_{H^k(\Omega^m)}\lesssim_Ac_m. 
\]
Now by the convolution inequality we conclude $\Vert\mathbf v^m\Vert_{H^s(\Omega^m)}\lesssim_AM$. 

To control $(h^m,\divg u^m)$ in $H^{s+\frac{1}{2}}\times H^{s-\frac{1}{2}}(\Omega^m)$, or equivalently, $(\nabla\divg u^m,\nabla^2h^m)$ in $H^{s-\frac{3}{2}}(\Omega^m)$, we adopt a similar argument as above based on \eqref{eq:hfbd1} and \eqref{eq:A2vt}. If $s\in\mathbb N$, we still need to control $\Vert z_s^m\Vert_{H^{1/2}_{00}(\Omega^m)}$. Using the analogue of \eqref{eq:telesum} for $(\divg u^m,\nabla h^m)$, using Lemma~\ref{lem:Phim} and using that $c_l$ is slowly varying, we derive $\Vert(\divg u^m,\nabla h^m)\Vert_{H^s(\Omega^m)}\lesssim_A2^\frac{m}{2}c_m$. Using the uniform bound on $\Vert\mathbf v^m\Vert_{H^s(\Omega^m)}$, \eqref{eq:trilest} and Sobolev embedding, and using that $z_s^m|_{\Gamma^m}=0$, we derive
\begin{equation}\label{eq:zsH01}
\Vert z_s^m\mathbf1_{\Omega^m}\Vert_{H^1(\mathbb R^3)}=\Vert z_s^m\Vert_{H_0^1(\Omega^m)}\lesssim_A2^\frac{m}{2}c_m+C(M). 
\end{equation}
Stein's extension theorem applies to the Lipschitz intersection $\tilde\Omega_m=\Omega^m\cap\Omega^{m+1}$, as explained in Appendix~\ref{app:estimates}. Thus Sobolev interpolation and the equivalence $W^{j,2}(\tilde\Omega_m)=H^j(\tilde\Omega_m)$ hold uniformly. Consequently, $\Vert\mathbf v^m-\mathbf v^{m+1}\Vert_{H^{s-1}(\tilde\Omega_m)}\lesssim_A2^{-m}c_m$ and $\Vert(\divg(u^m-u^{m+1}),\nabla(h^m-h^{m+1}))\Vert_{H^{s-1}(\tilde\Omega_m)}\lesssim_A2^{-m/2}c_m$. It follows that
\begin{equation}\label{eq:zsL2int}
\Vert z_s^m-z_s^{m+1}\Vert_{L^2(\tilde\Omega_m)}\lesssim_A2^{-\frac{m}{2}}c_m+2^{-m}C(M). 
\end{equation}
By \eqref{eq:zsH01} and Hardy's inequality we have
\[
\left(\int_{\Omega^m}\frac{|z_s^m(t,x)|^2}{d(x,\Gamma^m)^2}\ dx\right)^\frac{1}{2}+\left(\int_{\Omega^{m+1}}\frac{|z_s^{m+1}(t,x)|^2}{d(x,\Gamma^{m+1})^2}\ dx\right)^\frac{1}{2}\lesssim_A2^\frac{m}{2}c_m+C(M).
\]
On the other hand, now that $s\geq4$, by interpolation and Sobolev embedding we have $\Vert\eta^m-\eta^{m+1}\Vert_{L^\infty(\Gamma_*)}\lesssim_A2^{-2m}M$. Combining these and \eqref{eq:zsL2int} we get
\begin{equation}\label{eq:zsL2glo}
\Vert z_s^m\mathbf1_{\Omega^m}-z_s^{m+1}\mathbf1_{\Omega^{m+1}}\Vert_{L^2(\mathbb R^3)}\lesssim_A\left(2^{-\frac{m}{2}}c_m+2^{-m}C(M)\right)\left(1+2^{-2m}M\right).  
\end{equation}
By \eqref{eq:zsH01}, \eqref{eq:zsL2glo}, adopting the argument for $\Vert\eta^m\Vert_{H^s(\Gamma_*)}$, and noticing that $m_0(M)\gg1$, we obtain
\[
\Vert z_s^m\Vert_{H^{1/2}_{00}(\Omega^m)}\lesssim\Vert z_s^m\mathbf1_{\Omega^m}\Vert_{H^\frac{1}{2}(\mathbb R^3)}\lesssim_AM+1. 
\]
To sum up, we have
\[
\Vert(h^m,u^m,\Omega^m)\Vert_{\mathbf H^s_{1/2}}\lesssim_AM+1. 
\]
This improves the assumption on $\Vert(h^m,u^m,\Omega_t^m)\Vert_{\mathbf H^s_{1/2}}$ provided that $M_0$ is sufficiently large. It remains to improve the assumption on $A^m(t)$. Let $\phi^m_t:\Omega^m_0\to\Omega^m_t$ be the flow map: $\partial_t\phi^m_t=u^m(t,\phi^m_t)$. By Sobolev embedding we have $\Vert u^m\Vert_{C^{1,\delta}(\Omega^m_t)}\lesssim_AM_0$ and hence
\[
\Vert\phi_t^m-\text{Id}\Vert_{C^{1,\delta}(\Omega_0^m)}\lesssim_AM_0t. 
\]
Since $\Gamma^m_t$ moves with $\phi^m_t$, it follows that
\[
\Vert\Gamma^m_t\Vert_{C^{1,\delta}}\leq\Vert\Gamma^m_0\Vert_{C^{1,\delta}}+C(A)M_0t. 
\]
Similarly, we have
\[\begin{aligned}
\Vert u^m\Vert_{C^{\frac{1}{2}+\delta}(\Omega^m_t)}&\leq\Vert u^m(t,\phi^m_t)\Vert_{C^{\frac{1}{2}+\delta}(\Omega^m_0)}+C(A)M_0t\\
&\leq\Vert u^m\Vert_{C^{\frac{1}{2}+\delta}(\Omega^m_0)}+\int_0^t\Vert D_tu^m\circ(\tau,\phi^m_\tau)\Vert_{C^{\frac{1}{2}+\delta}(\Omega^m_0)}\,d\tau+C(A)M_0t\\
&\leq\Vert u^m\Vert_{C^{\frac{1}{2}+\delta}(\Omega^m_0)}+C(A)M_0t, 
\end{aligned}\]
\[\begin{aligned}
\Vert h^m\Vert_{C^{1,\delta}(\Omega^m_t)}&\leq\Vert h^m(t,\phi^m_t)\Vert_{C^{1,\delta}(\Omega^m_0)}+C(A)M_0t\\
&\leq\Vert h^m\Vert_{C^{1,\delta}(\Omega^m_0)}+\int_0^t\Vert D_th^m\circ(\tau,\phi^m_\tau)\Vert_{C^{1,\delta}(\Omega^m_0)}\,d\tau+C(A)M_0t\\
&\leq\Vert h^m\Vert_{C^{1,\delta}(\Omega^m_0)}+C(A)M_0t. 
\end{aligned}\]
Now by \eqref{eq:T0cond} we indeed improve the assumption \eqref{eq:btsass}. 

\subsection{Strong convergence to the limit solution}

The uniform estimates and compactness give convergence to a limit $(h,u,\Omega)$ in a weaker topology. By \eqref{eq:hfbd}--\eqref{eq:hfbd1}, each approximating solution is continuous in $\mathbf H^s_{1/2}$. We now strengthen the convergence to $C\mathbf H^s_{1/2}$.

We first prove convergence of $\Gamma^m$. For $j\in\mathbb N$ write
\[
P_{*j}(\eta^m-\eta)=\sum_{l=m}^jP_{*j}(\eta^l-\eta^{l+1})+\sum_{l>j}P_{*j}(\eta^l-\eta^{l+1}). 
\]
By Bernstein's estimate, the high-frequency bound \eqref{eq:hfbd} and the distance estimate \eqref{eq:lfdiff}, we have
\[
\Vert P_{*j}(\eta^m-\eta)\Vert_{H^s(\Gamma_*)}\lesssim_A\sum_{l=m}^{j-1}2^{-(k-s)(j-l)}c_l+\sum_{l=j}^\infty2^{-s(l-j)}c_l\lesssim_A\sum_l2^{-\min\{k-s,s\}|j-l|}\mathbf1_{l>m}c_l. 
\]
By the convolution inequality we obtain
\begin{equation}\label{eq:etacvg}
\Vert\eta^m-\eta\Vert_{H^s(\Gamma_*)}\lesssim_A\Vert c_{>m}\Vert_{\ell^2}, 
\end{equation}
which justifies the convergence of $\Gamma^m$. 

We next prove convergence of $\mathbf v^m=(h^m,u^m)$. For $p>m>m_0$ and for $j\in\mathbb N$ the decomposition
\[
\Phi_m\mathbf v^m-\Phi_p\mathbf v^p=\sum_{l=m}^{j-1}(\Phi_l\mathbf v^l-\Phi_{l+1}\mathbf v^{l+1})+\sum_{l=j}^{p-1}(\Phi_l\mathbf v^l-\Phi_{l+1}\mathbf v^{l+1})
\]
makes sense on $\Omega^p\cap\Omega$. Performing Littlewood--Paley decomposition $I=\sum P_j$ on $\Omega^p,\Omega$ and arguing as above, we obtain
\begin{equation}\label{eq:truncdiff}
\Vert\Phi_m\mathbf v^m-\Phi_p\mathbf v^p\Vert_{H^s(\Omega^p)}+\Vert\Phi_m\mathbf v^m-\mathbf v\Vert_{H^s(\Omega)}\lesssim_A\Vert c_{>m}\Vert_{\ell^2}. 
\end{equation}
%Recall the extension operator depending continuously on the domain. 
By Proposition~\ref{prop:ctext} we have
\begin{equation}\label{eq:target}
\lim_{p\to\infty}\Vert E_{\Omega_t^p}\Phi_m\mathbf v^m(t)-E_{\Omega_t}\Phi_m\mathbf v^m(t)\Vert_{H^s(\mathbb R^3)}=0,\quad\forall t.
\end{equation}
To obtain uniform convergence in time, we use the following elementary operator lemma.

\begin{lemma}\label{lem:funcanal}
Let $X$ be a Banach space and $Y$ a normed space. Let $E^j(t),E(t):X\to Y$ be bounded linear operators on a compact time interval. Suppose $E(t)$ is strongly continuous and
\[
E^j(t)f\to E(t_0)f\,\quad j\to\infty,\ t\to t_0
\]
for every $f\in X$. Suppose $v^j(t)$ is a sequence of continuous curves in $X$ that converge uniformly to $v(t)$. Then we have
\[
E^j(t)v^j(t)\to E(t)v(t)\text{ uniformly in }t. 
\]
\end{lemma}
\proofheading{Proof}Fix sequences $j_\ell\to\infty$ and $t_\ell\to t_0$. The hypothesis gives strong convergence $E^{j_\ell}(t_\ell)\to E(t_0)$. The uniform boundedness principle therefore bounds their operator norms. Since $v^{j_\ell}(t_\ell)\to v(t_0)$, the corresponding images converge to $E(t_0)v(t_0)$. Compactness of the time interval yields the asserted uniform convergence. \proofqed

Apply the lemma with $X=H^s(\tilde\Omega)$, $Y=H^s(\mathbb R^3)$, $E^j(t)f=E_{\Omega_t^j}(f|_{\Omega_t^j})$, $E(t)f=E_{\Omega_t}(f|_{\Omega_t})$, and $v^j(t)=v(t)=\Phi_m\mathbf v^m(t)$. Here $m$ is fixed and $\tilde\Omega$ is a neighborhood of $\overline{\Omega_{t_0}}$. Proposition~\ref{prop:ctext} and convergence of the boundaries verify the hypotheses on $[t_0,t_0+c2^{-m}]$. A finite cover gives uniformity in \eqref{eq:target}. Now write
\[\begin{aligned}
E_{\Omega^p}\mathbf v^p-E_\Omega\mathbf v=E_{\Omega^p}(I-\Phi_p)\mathbf v^p-E_{\Omega^p}(\Phi_m\mathbf v^m-\Phi_p\mathbf v^p)+E_\Omega(\Phi_m\mathbf v^m-\mathbf v)\\
+E_{\Omega^p}\Phi_m\mathbf v^m-E_\Omega\Phi_m\mathbf v^m. 
\end{aligned}\]
The second line vanishes in $H^s(\mathbb R^3)$ as $p\to\infty$ first, followed by $m\to\infty$. For the first line, we use \eqref{eq:truncdiff}, Lemma~\ref{lem:Phim} and the boundedness of $E_{\Omega^p},E_\Omega$, yielding that
\[
\Vert E_{\Omega^p}(I-\Phi_p)\mathbf v^p\Vert_{H^s(\mathbb R^3)}+\Vert E_{\Omega^p}(\Phi_m\mathbf v^m-\Phi_p\mathbf v^p)\Vert_{H^s(\mathbb R^3)}+\Vert E_\Omega(\Phi_m\mathbf v^m-\mathbf v)\Vert_{H^s(\mathbb R^3)}\lesssim_M\Vert c_{>m}\Vert_{\ell^2}. 
\]
This justifies the convergence of $\mathbf v^m$. 

The verification of the convergence of $(h^m,\divg u^m)$ in $H^{s+\frac{1}{2}}\times H^{s-\frac{1}{2}}$ is similar. From \eqref{eq:hfbd1}, \eqref{eq:A2vt} we derive: 
\[\begin{aligned}
&\Vert\Phi_m(\nabla\divg u^m,\nabla^2h^m)-\Phi_p(\nabla\divg u^p,\nabla^2h^p)\Vert_{H^{s-\frac{3}{2}}(\Omega^p)}\lesssim_A\Vert c_{>m}\Vert_{\ell^2},\\
&\Vert\Phi_m(\nabla\divg u^m,\nabla^2h^m)-(\nabla\divg u,\nabla^2h)\Vert_{H^{s-\frac{3}{2}}(\Omega)}\lesssim_A\Vert c_{>m}\Vert_{\ell^2}. 
\end{aligned}\]
By the estimate of $[\nabla,\Phi_m]$, it follows that
\begin{equation}\label{eq:truncdiffcomp}
\begin{aligned}
&\Vert\Phi_mh^m-\Phi_ph^p\Vert_{H^{s+\frac{1}{2}}(\Omega^p)}+\Vert\Phi_mh^m-h\Vert_{H^{s+\frac{1}{2}}(\Omega)}\lesssim_A\Vert c_{>m}\Vert_{\ell^2},\\
&\Vert\Phi_m\divg u^m-\Phi_p\divg u^p\Vert_{H^{s-\frac{1}{2}}(\Omega^p)}+\Vert\Phi_m\divg u^m-\divg u\Vert_{H^{s-\frac{1}{2}}(\Omega)}\lesssim_A\Vert c_{>m}\Vert_{\ell^2}.
\end{aligned}
\end{equation}
Then the argument above indicates that $E_{\Omega^m}(h^m,\divg u^m)\to E_\Omega(h,\divg u)$ in $C([0,T_0],H^{s+\frac{1}{2}}\times H^{s-\frac{1}{2}}(\mathbb R^3))$. If $s\in\mathbb N$, the convergence of $z_s^m\mathbf1_{\Omega^m}$ in $H^\frac{1}{2}(\mathbb R^3)$ follows from \eqref{eq:zsH01}, \eqref{eq:zsL2glo} and the argument for $\eta^m$. In fact, we have
\begin{equation}\label{eq:zsH1/2cvg}
\Vert z_s^m\mathbf1_{\Omega^m}-z_s\mathbf1_\Omega\Vert_{H^\frac{1}{2}(\mathbb R^3)}\lesssim_A\Vert c_{>m}\Vert_{\ell^2}+2^{-\frac{m}{2}}C(M). 
\end{equation}
This concludes the convergence of $(h^m,u^m,\Omega^m)$ to the limit solution in $C\mathbf H^s_{1/2}$. 

\subsection{Continuous dependence}

Suppose $(h_l,u_l,\Omega_l)$ is a sequence of solutions whose initial data converge to $(h_0,u_0,\Omega_0)$ in $\mathbf H^s_{1/2}$. Let $(h_l^m,u_l^m,\Omega_l^m)$ be the solution subject to $(h_l^m,u_l^m,\Omega_l^m)|_{t=0}=\Psi_m(h_l,u_l,\Omega_l)|_{t=0}$. For each fixed $m$ we claim that $(h_l^m,u_l^m,\Omega_l^m)\to(h^m,u^m,\Omega^m)$ in $C\mathbf H^s_{1/2}$. Indeed, combining \eqref{eq:Psimconti} and Theorem~\ref{thm:L2dist} we see
\[
\lim_{l\to\infty}\sup_t\left(\Vert\eta_l^m(t)-\eta^m(t)\Vert_{L^2(\Gamma_*)}+\Vert\mathbf v_l^m(t)-\mathbf v^m(t)\Vert_{L^2(\Omega_l^m\cap\Omega^m)}\right)=0. 
\]
From the kinematic boundary condition \eqref{eq:kinbc} we see $\Vert\partial_t\eta_l^m\Vert_{H^{k-1/2}(\Gamma_*)}$ is uniformly bounded. Using interpolation, we thus verify the two-index convergence $\Gamma_l^m(t)\to\Gamma^m(t_0)$ in $H^{k'}$ for every $k'<k$. Next, write
\[\begin{aligned}
E_{\Omega_l^m}\mathbf v_l^m-E_{\Omega^m}\mathbf v^m=E_{\Omega_l^m}(I-\Phi_p)\mathbf v_l^m+E_{\Omega_l^m}\Phi_p(\mathbf v_l^m-\mathbf v^m)-E_{\Omega^m}(I-\Phi_p)\mathbf v^m\\
+E_{\Omega_l^m}\Phi_p\mathbf v^m-E_{\Omega^m}\Phi_p\mathbf v^m. 
\end{aligned}\]
By Lemma~\ref{lem:Phim} the right-hand side in the first line vanishes in $C([0,T],H^{k'}(\mathbb R^3))$ as $l\to\infty$ followed by $p\to\infty$. By Lemma~\ref{lem:funcanal} the second line vanishes in $C([0,T],H^{k'}(\mathbb R^3))$ as $l\to\infty$. This justifies the convergence of $\mathbf v_l^m$. If $s\in\mathbb N$, this together with Hardy's inequality implies the $L^2$-convergence of $z_{s,l}^m\mathbf1_{\Omega_l^m}$. Interpolating with the uniform $H_0^1$-bound, this implies the $H^\frac{1}{2}$-convergence of $z_{s,l}^m\mathbf1_{\Omega_l^m}$ and thus concludes the proof of the claim. 

It remains to prove convergence of the regularized solutions uniformly over the sequence of initial states:
\begin{equation}\label{eq:unicvg}
(h_l^m,u_l^m,\Omega_l^m)\to(h_l,u_l,\Omega_l)\text{ in }C\mathbf H^s_{1/2},\text{ uniformly in }l. 
\end{equation}

Let $(c_j^l)$ be the frequency envelope of $(h_l,u_l,\Omega_l)|_{t=0}$ in $\mathbf H^s_{1/2}$. Corollary~\ref{cor:unil2} implies
\[
\sigma_m:=\sup_l\Vert c_{>m}^l\Vert_{\ell^2}\to0\text{ as }m\to\infty. 
\]
We have the following analogues of \eqref{eq:hfbd}, \eqref{eq:hfbd1}, \eqref{eq:etacvg}, \eqref{eq:truncdiff}, \eqref{eq:truncdiffcomp}: 
\begin{equation}\label{eq:hfbd2}
\Vert(h_l^m,u_l^m,\Omega_l^m)\Vert_{\mathbf H^k}\lesssim_A2^{(k-s)m}c_m^l,
\end{equation}
\begin{equation}\label{eq:hfbd3}
\Vert(h_l^m,u_l^m,\Omega_l^m)\Vert_{\mathbf H^k_1}\lesssim_A2^{(k+\frac{1}{2}-s)m}c_m^l,
\end{equation}
\begin{equation}\label{eq:etacvg1}
\Vert\eta_l^m-\eta_l\Vert_{H^s(\Gamma_*)}\lesssim_A\sigma_m, 
\end{equation}
\begin{equation}\label{eq:truncdiff1}
\Vert\Phi_m\mathbf v_l^m-\Phi_p\mathbf v_l^p\Vert_{H^s(\Omega_l^p)}+\Vert\Phi_m\mathbf v_l^m-\mathbf v_l\Vert_{H^s(\Omega_l)}\lesssim_A\sigma_m. 
\end{equation}
\begin{equation}\label{eq:truncdiffcomp1}
\begin{aligned}
&\Vert\Phi_mh_l^m-\Phi_ph_l^p\Vert_{H^{s+\frac{1}{2}}(\Omega_l^p)}+\Vert\Phi_mh_l^m-h_l\Vert_{H^{s+\frac{1}{2}}(\Omega_l)}\lesssim_A\sigma_m,\\
&\Vert\Phi_m\divg u_l^m-\Phi_p\divg u_l^p\Vert_{H^{s-\frac{1}{2}}(\Omega_l^p)}+\Vert\Phi_m\divg u_l^m-\divg u_l\Vert_{H^{s-\frac{1}{2}}(\Omega_l)}\lesssim_A\sigma_m.
\end{aligned}
\end{equation}
If $s\in\mathbb N$, as in \eqref{eq:zsH1/2cvg} we have
\begin{equation}\label{eq:zsunicvg}
\Vert z_{s,l}^m\mathbf1_{\Omega_l^m}-z_{s,l}\mathbf1_{\Omega_l}\Vert_{H^\frac{1}{2}(\mathbb R^3)}\lesssim_A\sigma_m+2^{-\frac{m}{2}}C(M). 
\end{equation}
\eqref{eq:etacvg1} and \eqref{eq:zsunicvg} already imply the uniform-in-$l$ convergence of $\Gamma_l^m$ and $\mathbf1_{s\in\mathbb N}z_{s,l}^m\mathbf1_{\Omega_l^m}$. For $\mathbf v^m$, by \eqref{eq:hfbd2}, \eqref{eq:truncdiff1} and the argument in the last subsection we have
\begin{equation}\label{eq:extvlp}
E_{\Omega^p_l}\mathbf v_l^p-E_{\Omega_l}\mathbf v_l=E_{\Omega^p_l}\Phi_m\mathbf v_l^m-E_{\Omega_l}\Phi_m\mathbf v_l^m+O_{H^s(\mathbb R^3)}(\sigma_m). 
\end{equation} 
For the first two terms we invoke Lemma~\ref{lem:funcanal} with a double index $\mathbf j=(j,l)$. Set $X=H^s(\tilde\Omega),Y=H^s(\mathbb R^3)$ where $\tilde\Omega$ is some neighborhood of $\Omega(t_0)$. Set $E^{\mathbf j}(t)=E_{\Omega_l^j(t)},E(t)=E_{\Omega(t)}$, $v^{\mathbf j}(t)=\Phi_m\mathbf v_l^m(t)$ and $v(t)=\Phi_m\mathbf v^m(t)$, Since $m$ is fixed, we have $v^{\mathbf j}\to v$ uniformly. Further by the convergence of $\Gamma_l^j(t)$ and by Proposition~\ref{prop:ctext}, the assumptions of the lemma are satisfied on $[t_0,t_0+c2^{-m}]$, and it follows that $E_{\Omega^p_l}\Phi_m\mathbf v_l^m\to E_\Omega\Phi_m\mathbf v^m$ in $C([0,T_0],H^s(\mathbb R^3))$ as $p,l\to\infty$. Similarly, setting $E^{\mathbf j}(t)=E_{\Omega_l(t)}$ we see $E_{\Omega_l}\Phi_m\mathbf v_l^m\to E_\Omega\Phi_m\mathbf v^m$ in $C([0,T_0],H^s(\mathbb R^3))$ as $l\to\infty$. Combining these conclusions and the convergence result for fixed $l$ we conclude that the difference of the first two terms in the right-hand side of \eqref{eq:extvlp} vanishes uniformly in $t,l$ as $p\to\infty$. Now the uniform-in-$l$ convergence of $\mathbf v_l^p$ is clear. 

The above argument combined with \eqref{eq:hfbd3}, \eqref{eq:truncdiffcomp1} also indicates the uniform-in-$l$ convergence of $E_{\Omega^p_l}(h^p_l,\divg u_l^p)$ in $H^{s+\frac{1}{2}}\times H^{s-\frac{1}{2}}(\mathbb R^3))$. This concludes \eqref{eq:unicvg} and completes the proof of Theorem~\ref{thm:LWP}.

\section*{Acknowledgement}
Luo's research is supported by a grant from the Research Grants Council of the
Hong Kong Special Administrative Region, China (Project No. 11310023). Zhang is grateful to Professor Huihui Zeng for her discussions and guidance. 
Zhang is also supported by a grant of NSFC (Grant number: 12671261) through Professor Huihui Zeng. 

\appendix

\section{Sobolev spaces and elliptic estimates}\label{app:estimates}

This appendix collects the extension, product, trace, and elliptic estimates used above. Throughout, assume $\Gamma\in\Lambda_*$ and $\Vert\Gamma\Vert_{C^{1,\delta}}\leq A$. All implicit constants may depend on the fixed collar.

\subsection{Extension operators and multilinear estimates}

If $\phi:\mathbb R^{d-1}\to\mathbb R$ has finite Lipschitz constant $M$, define an open set $U_\phi=\{(x',x_d)\mid x_d>\phi(x')\}$.  A
classical result of Stein, Theorem 5', p. 181 in \cite{Stein}, asserts that there exists a linear operator $\mathcal E=\mathcal E_\phi$ mapping functions on $U_\phi$ to functions on $\mathbb R^d$ with the property that $\mathcal E f|_{U_\phi}=f$ and $\mathcal E:W^{k,p}(U_\phi)\to W^{k,p}(\mathbb R^d)$ is bounded for all $k\geq0$, $p\in[1,\infty]$. Moreover, the operator norm depends only on $d,k$ and $M$. 

A partition of unity gives an extension operator $\mathcal E_\Omega$ with the same mapping properties on a Lipschitz domain. For boundaries in $\Lambda_*$, the relevant operator norms are uniform.

Recall that $H^s(\Omega)$ is defined as the quotient space $H^s(\mathbb R^d)/\{f\text{ vanishes on }\Omega\}$. Also recall that $(H^{s_0}(\Omega),H^{s_1}(\Omega))_{\theta,2}=H^s(\Omega)$ if $s=\theta s_1+(1-\theta)s_0$. Thus, the mapping property of $\mathcal E$ can be extended to fractional regularity spaces. 

\begin{proposition}\label{prop:SteinExt}
Stein's extension operator $\mathcal E_\Omega$ enjoys the following properties: 
\[
\Vert\mathcal E_\Omega f\Vert_{C^\alpha(\mathbb R^d)}\lesssim_A\Vert f\Vert_{C^\alpha(\Omega)},\quad \alpha\in[0,1+\delta],
\]
\[
\Vert\mathcal E_\Omega f\Vert_{W^{k,p}(\mathbb R^d)}\lesssim_A\Vert f\Vert_{W^{k,p}(\Omega)},\quad k\geq0,\ p\in[1,\infty]. 
\]
\end{proposition}
\proofheading{Proof}See Proposition 5.1 in \cite{IE}. \proofqed
\par\smallskip\noindent\textbf{Remark.}\ The estimates actually hold for a larger class of domains, for example, $\tilde\Omega=\Omega_1\cap\Omega_2$ with $\Gamma_i\in\Lambda_*$. Indeed, as explained in Section~\ref{sub:l2distance}, $\tilde\Gamma$ is parameterized by $\min\{\eta_{\Gamma_1},\eta_{\Gamma_2}\}$, which is a Lipschitz function on $\Gamma_*$. The upper bound on Lipschitz constant and the partition of unity are uniform for such domains. As a corollary, we have
\begin{equation}\label{eq:Holdinter}
\Vert f\Vert_{C^\alpha(\tilde\Omega)}\lesssim\Vert f\Vert_{C^{\alpha_0}(\tilde\Omega)}^{1-\theta}\Vert f\Vert_{C^{\alpha_0}(\tilde\Omega)}^\theta,\quad\alpha=(1-\theta)\alpha_0+\theta\alpha_1,
\end{equation}
where the implicit constant depends only on the collar $\Lambda_*$.\par\smallskip

We use Stein's extension to define Littlewood--Paley operators on $\Omega$: by abuse of notation, $P_jv$ denotes the restriction of $P_j\mathcal E_\Omega v$ to $\Omega$. Thus $v=\sum_{j\geq0}P_jv$ on $\Omega$. The operators $P_{<j}$ and $P_{\geq j}$ denote the corresponding low- and high-frequency sums.

\begin{proposition}\label{prop:bilest2}
Let $s>0$, $\alpha_1,\alpha_2,\beta\in[0,1]$. For every $r\geq0$ and every sequence of partitions $g=g_j^1+g_j^2$, we have
\[\begin{aligned}
\Vert\partial_kf\partial_lg\Vert_{H^s(\Omega)}\lesssim_A\Vert f\Vert_{C^{\alpha_1}(\Omega)}\Vert g\Vert_{H^{s+2-\alpha_1}(\Omega)}&+\Vert f\Vert_{H^{s+r+1}(\Omega)}\sup_{j>0}2^{-j(r+\alpha_2-1)}\Vert g_j^1\Vert_{C^{\alpha_2}(\Omega)}\\
&+\Vert f\Vert_{C^{\alpha_1+\delta}(\Omega)}\sup_{j>0}2^{j(s+1-\alpha_1+\beta)}\Vert g_j^2\Vert_{H^{1-\beta}(\Omega)}.
\end{aligned}\]
\end{proposition}
\proofheading{Proof}See Proposition 5.3 in \cite{IE}. Although there is a modification of the index in the last term, the proofs are essentially the same. \proofqed

\begin{proposition}\label{prop:trilest}
Let $s\geq0$, $r_i\geq0,\alpha_i\in[0,r_i],i=1,...,m$. Then if $s>0$, we have
\begin{equation}\label{eq:trilest}
\Vert\nabla^{r_1}f_1\cdots\nabla^{r_m}f_m\Vert_{H^s(\Omega)}\lesssim_A\sum_{i=1}^m\Vert f_i\Vert_{H^{\sigma_i}(\Omega)}\prod_{j\neq i}\Vert f_j\Vert_{C^{\alpha_j}(\Omega)},
\end{equation}
where $\sigma_i=s+r_i+\sum_{j\neq i}(r_j-\alpha_j)$. If $s=0$, the inequality remains true if the non-integral $\alpha_j$ are replaced by $\alpha_j+\delta'$, $\delta'>0$. 
\end{proposition}
\proofheading{Proof}Proposition~\ref{prop:SteinExt} reduces the argument to $\Omega=\mathbb R^d$.
Let us first consider $s>0$. Denote $F=\nabla^{r_1}f_1\cdots\nabla^{r_m}f_m$. Note the identity
\begin{align}
F&=\sum_{j_1,...,j_m}P_{j_1}\nabla^{r_1}f_1\cdots P_{j_m}\nabla^{r_m}f_m\nonumber\\
&=\sum_jP_j\nabla^{r_1}f_1P_{<j}\nabla^{r_2}f_2\cdots P_{<j}\nabla^{r_m}f_m+\cdots+\sum_jP_{\leq j}\nabla^{r_1}f_1\cdots P_{\leq j}\nabla^{r_{m-1}}f_{m-1}P_j\nabla^{r_m}f_m\nonumber\\
&=\sum_j\left(T_j^1+\cdots+T_j^m\right)\label{eq:Tji}. 
\end{align}
The Fourier transform of $T_k^i$ is supported in a ball of radius $2^{k+N/2}$, thus we have
\[
P_jF=P_j\sum_{k>j-N}\left(T_k^1+\cdots+T_k^m\right). 
\]
By Bernstein's estimate, for $i=1,...,m$ we have
\begin{equation}\label{eq:TkiL2}
\Vert T_k^i\Vert_{L^2}\lesssim2^{-k(\sigma_i-r_i)}\prod_{h\neq i}2^{k(r_h-\alpha_h)}\,c_k=2^{-ks}c_k,
\end{equation}
where $\Vert c_k\Vert_{\ell^2_k}$ is dominated by the right-hand side of \eqref{eq:trilest}. It follows that
\[
2^{js}\Vert P_jF\Vert_{L^2}\lesssim\sum_{k>j-N}2^{-(k-j)s}c_k. 
\]
By the convolution inequality we conclude \eqref{eq:trilest} for $s>0$. 

We next consider $s=0$. If $r_1-\alpha_1=0$, we have $\nabla^{r_1}f_1\in L^\infty$ and can simply drop it. If $r_1-\alpha_1$ is a positive integer, we iteratively use the identity 
\[
\nabla^{r_1}f_1\cdots\nabla^{r_m}f_m=\nabla\left(\nabla^{r_1-1}f_1\nabla^{r_2}f_2\cdots\nabla^{r_m}f_m\right)-\nabla^{r_1-1}f_1\nabla\left(\nabla^{r_2}f_2\cdots\nabla^{r_m}f_m\right)
\]
until the number of derivatives on $f_1$ reduces to $\alpha_1$. Note that the first item on the right-hand side can be handled by the case $s>0$. Thus, we may assume all $\alpha_i$ in \eqref{eq:trilest} are added by $\delta'>0$. Go back to \eqref{eq:Tji}--\eqref{eq:TkiL2}. This time we have 
\[
\Vert T_k^i\Vert_{L^2}\lesssim2^{-k(s+\delta'')}c_k,
\]
where $\delta''=(m-1)\delta'$. The resulting geometric decay of $\Vert P_jF\Vert_{L^2}$ gives the required summability. \proofqed

\subsection{Sobolev norms on the boundary}

Compactness of $\Gamma_*$ gives a finite collection of rotated cylinders $R_i(2r_i)$, of radius and length $2r_i$, with the following properties:

(i) The union of $R_i(r_i)$ covers a neighborhood of $\Gamma_*$. 

(ii) There are $f_{*i}:B_i(2r_i)\to(-\sigma r_i,\sigma r_i)$ such that
\[
\Vert Df_{*i}\Vert_{C^0}<\sigma,\quad\Omega_*\cap R_i(2r_i)=\{y_3>f_{*i}(\bar y)\},
\]
where $B_i(2r_i)$ is the base disk of $R_i(2r_i)$ and $(\bar y,y_3)$ denotes the Euclidean coordinates on $R_i$. 

If $\Gamma\in\Lambda_*$ and $\delta_0\ll1$, (i) holds for $\Gamma$ and there are (unique) functions $f_i$ satisfying (ii) for $\Omega$. Moreover, we have
\[
\Vert f_i\Vert_{H^s}\lesssim_A1+\Vert\Gamma\Vert_{H^s},\quad\Vert f_i\Vert_{C^{k,\alpha}}\lesssim_A1+\Vert\Gamma\Vert_{C^{k,\alpha}}. 
\]

Let $\bar\beta:[0,\infty)\to[0,1]$ be a smooth function that equals $1$ on $[0,\frac{5}{4}]$ and vanishes on $[\frac{3}{2},\infty)$. Let $\beta_i(\bar y)=\bar\beta(\frac{|\bar y|}{r_i})$. Define an extension of $\beta_if_i$ by
\[
F_i(\bar y,y_3)=\int_{\mathbb R^2}e^{2\pi i\bar y\cdot\xi'}e^{-(1+|\xi'|^2)y_3^2}\widehat{\beta_if_i}(\xi')\,d\xi', 
\]
where $\hat f(\xi')$ is the Fourier transform of $f(\bar y)$. By a standard analysis we can show $\Vert F_i\Vert_{C^{k,\alpha}(\mathbb R^3)}\lesssim\Vert\beta_if_i\Vert_{C^{k,\alpha}(\mathbb R^2)}$ and $\Vert F_i\Vert_{H^{s+\frac{1}{2}}(\mathbb R^3)}\lesssim\Vert\beta_if_i\Vert_{H^s(\mathbb R^2)}$. If $\sigma$, the uniform bound on $Df_i$, is sufficiently small, the map 
\[
H_i(\bar y,y_3)=(\bar y,\ y_3+F_i(\bar y,y_3))
\] 
is a diffeomorphism on $\mathbb R^3$. Let $G_i$ be the inverse map of $H_i$ and let $g_i$ be the $3$rd component of $G_i$. Then for some uniform constant $\delta_*>0$ we have
\[
\left(B_i(\frac{5}{4}r_i)\times I_i(\frac{5}{4}\delta_*r_i)\right)\cap\Omega=\left(B_i(\frac{5}{4}r_i)\times I_i(\frac{5}{4}\delta_*r_i)\right)\cap\{g_i>0\}. 
\]
Moreover, we have $\Vert H_i-\text{Id}\Vert_{C^{1,\delta}}+\Vert G_i-\text{Id}\Vert_{C^{1,\delta}}\lesssim_A1$ and
\[
\Vert H_i-\text{Id}\Vert_{H^{s+\frac{1}{2}}(\mathbb R^3)}+\Vert G_i-\text{Id}\Vert_{H^{s+\frac{1}{2}}(\mathbb R^3)}\lesssim_A\Vert\Gamma\Vert_{H^s}. 
\]

We now construct a partition of unity for $\Omega$ with bounds uniform in $\Lambda_*$. Let $\bar\alpha:[0,\infty)\to[0,1]$ be a smooth function that equals $1$ on $[0,\frac{9}{8}]$ and vanishes on $[\frac{5}{4},\infty)$. Let $\bar\rho:[0,\infty)\to[\frac{1}{3},\infty)$ be a smooth function such that $\bar\rho=\frac{1}{3}$ on $[0,\frac{1}{3}]$ and $\bar\rho(x)=x$ for $x\geq\frac{2}{3}$. Fix $0<\rho_*\leq\delta_*$ and define
\[
\tilde\alpha_{*i}(\bar y,y_3)=\bar\alpha(\frac{|\bar y|}{r_i})\bar\alpha(\frac{|y_3|}{\rho_*r_i}),\quad\rho=\bar\rho\circ\sum_i\tilde\alpha_{*i}(G_i)
\]
and
\[
\alpha_{*i}=\frac{\tilde\alpha_{*i}(G_i)}{\rho},\quad\alpha_{*0}=(1-\sum_i\alpha_{*i})\mathbf1_\Omega. 
\]

These charts and the partition of unity define Sobolev and H\"older norms on $\Gamma$. For $f,g$ on $\Gamma$, define their $H^s(\Gamma)$ inner product by summing the $H^s(\mathbb R^2)$ inner products of $(\alpha_{*i}f)\circ H_i(\cdot,0)$ and $(\alpha_{*i}g)\circ H_i(\cdot,0)$. For a $C^{k,\alpha}$ surface, define the $C^{k,\alpha}(\Gamma)$ norm by the corresponding sum of local $C^{k,\alpha}(\mathbb R^2)$ norms.
\par\smallskip\noindent\textbf{Remark.}\ By definition we have
\[
\Vert f\Vert_{H^s(\Gamma)}\leq\Vert f\Vert_{H^{s_0}(\Gamma)}^{1-\theta}\Vert f\Vert_{H^{s_1}(\Gamma)}^\theta,\quad s=(1-\theta)s_0+\theta s_1
\]
and
\[
\Vert f\Vert_{C^\alpha(\Gamma)}\leq\Vert f\Vert_{C^{\alpha_0}(\Gamma)}^{1-\theta}\Vert f\Vert_{C^{\alpha_1}(\Gamma)}^\theta,\quad \alpha=(1-\theta)\alpha_0+\theta \alpha_1.
\]

\begin{proposition}\label{prop:Moser2}
Let $G:\mathbb R^d\to\mathbb R^d$ be a diffeomorphism and $\Vert\nabla G\Vert_{L^\infty}+\Vert\nabla G^{-1}\Vert_{L^\infty}\leq B^\sharp$, $\Vert\nabla G\Vert_{C^\delta}+\Vert\nabla G^{-1}\Vert_{C^\delta}\leq B$. For $s\in[0,1]$ we have
\[
\Vert f\circ G\Vert_{H^s(\mathbb R^d)}\sim_{B^\sharp}\Vert f\Vert_{H^s(\mathbb R^d)}. 
\]
For $s,r\geq0$, $\alpha,\beta\in[0,1]$ and partitions $f=f_j^1+f_j^2$ we have
\[\begin{aligned}
\Vert f\circ G\Vert_{H^s(\mathbb R^d)}\lesssim_B\Vert f\Vert_{H^s}+\Vert G-\text{Id}\Vert_{H^{s+r}}\sup_{j\geq0}2^{-j(r+\alpha-1)}\Vert f_j^1\Vert_{C^\alpha}+\sup_{j\geq0}2^{j(s+\beta-1-\frac{\delta}{2})}\Vert f_j^2\Vert_{H^{1-\beta}}. 
\end{aligned}\]
\end{proposition}
\proofheading{Proof}See Proposition 5.5 in \cite{IE}. \proofqed

\begin{proposition}\label{prop:tr}
For $s\in(\frac{1}{2},1+\delta)$ we have
\[
\Vert f|_\Gamma\Vert_{H^{s-\frac{1}{2}}(\Gamma)}\lesssim_A\Vert f\Vert_{H^s(\Omega)}. 
\]
For $s>\frac{1}{2}$, $r\geq0$, $\alpha,\beta\in[0,1]$ and partitions $f=f_j^1+f_j^2$, we have
\[\begin{aligned}
\Vert f|_\Gamma\Vert_{H^{s-\frac{1}{2}}(\Gamma)}\lesssim_A\Vert f\Vert_{H^s(\Omega)}+\Vert\Gamma\Vert_{H^{s+r-\frac{1}{2}}}\sup_{j>0}2^{-j(r+\alpha-1)}\Vert f_j^1\Vert_{C^\alpha(\Omega)}+\sup_{j>0}2^{j(s-1+\beta-\delta)}\Vert f_j^2\Vert_{H^{1-\beta}(\Omega)}.
\end{aligned}\]
\end{proposition}
\proofheading{Proof}See Proposition 5.11 in \cite{IE}. \proofqed

\begin{proposition}\label{prop:mlinb}
Let $\nabla_\Gamma$ be the covariant differentiation on $\Gamma$ and $r_i\geq0,\alpha_i\in\{0,1\},i=1,...,m$. Then we have
\[
\Vert\nabla_\Gamma^{r_1}v_1\cdots\nabla_\Gamma^{r_m}v_m\Vert_{L^2(\Gamma)}\lesssim_A\sum_{i=1}^m\Vert v_i\Vert_{H^{\sigma_i}(\Gamma)}\prod_{j\neq i}\Vert v_j\Vert_{C^{\alpha_j}(\Gamma)}+\Vert\Gamma\Vert_{H^{\sigma_0}}\prod_{i=1}^m\Vert v_i\Vert_{C^{\alpha_i}(\Gamma)},
\]
where $\sigma_i=r_i+\sum_{j\neq i}(r_j-\alpha_j),\ i=1,...,m$ and $\sigma_0=1+\sum_{i=1}^m(r_i-\alpha_i)$. 
\end{proposition}
\proofheading{Proof}We first record the geometry of a graph. For $\Gamma_f=\{x=(\bar y,f(\bar y))\}$,
\[
n=\frac{(-\nabla_{\bar y}f,1)}{\sqrt{1+|\nabla_{\bar y}f|^2}},\quad\partial_i^\top\partial_{y_i}-n_in^j\partial_{y_j},\quad\partial_3^\top=-n_3n^j\partial_{y_j}
\]
where $n$ is unit normal and $(\partial_i^\top,\partial_3^\top)$ is the covariant differentiation on $\Gamma_f$. 

Localize with $\alpha_{*i}$ and use the chart $H_i$. Fix $i$ and set $\tilde v_j=v_j(H_i)$. In these coordinates, $\nabla_\Gamma^{r_1}v_1\cdots\nabla_\Gamma^{r_m}v_m$ is a sum of terms of the form
\[
F(n)\,\nabla_{\bar y}^{t_1}n\cdots\nabla_{\bar y}^{t_l}n\,\nabla_{\bar y}^{s_1}\tilde v_1\cdots\nabla_{\bar y}^{s_m}\tilde v_m
\]
where $F(n)$ is a smooth function, $t_j\geq0$, $s_j\geq1$ and $t_1+\cdots+t_l+s_1+\cdots+s_m=r_1+\cdots+r_m$. By Proposition~\ref{prop:trilest}, which also holds for planar domains, we have
\[
\Vert\nabla_\Gamma^{r_1}v_1\cdots\nabla_\Gamma^{r_m}v_m\Vert_{L^2(U_i)}\lesssim_A\sum_{j=1}^m\Vert\tilde v_j\Vert_{H^{\sigma_j}(B_i)}\prod_{k\neq j}\Vert\tilde v_k\Vert_{C^{\alpha_k}(B_i)}+\Vert\Gamma\Vert_{H^{\sigma_0}}\prod_{j=1}^m\Vert\tilde v_j\Vert_{C^{\alpha_j}(B_i)}, 
\] 
where $U_i=H_i(B_i,0)$. Write
\[
\tilde v_j=\sum_{i'}(\alpha_{*i'}v_j)\circ H_{i'}\circ(G_{i'}\circ H_i).
\]
Since all indices $\alpha_j$ belong to $\{0,1\}$, we have $\Vert\tilde v_j\Vert_{C^{\alpha_j}(B_i)}\lesssim_A\Vert v_j\Vert_{C^{\alpha_j}(\Gamma)}$. For the Sobolev norm of $\tilde v_j$, we use Moser's estimate (Proposition~\ref{prop:Moser2}). This completes the proof. \proofqed

\begin{proposition}\label{prop:ctext}
For each domain $\Omega$ with boundary $\Gamma\in H^s\cap\Lambda_*$ there is an $E_\Omega:H^s(\Omega)\to H^s(\mathbb R^3)$ with the following properties: (i)
\[
\Vert E_\Omega v\Vert_{H^s(\mathbb R^3)}\lesssim_A(1+\Vert\Gamma\Vert_{H^{s-\frac{1}{2}}})\Vert v\Vert_{H^s(\Omega)}.
\]
(ii) For $v\in H^s(\mathbb R^3)$ and sequence of domains $\Omega_m$ with $\Gamma_m\to\Gamma$ in $H^s$ we have
\[
\Vert E_{\Omega_m}v-E_\Omega v\Vert_{H^s(\mathbb R^3)}\to0. 
\]
\end{proposition}
\proofheading{Proof}See Proposition 5.12 in \cite{IE}. \proofqed

\subsection{Elliptic estimates}

\begin{proposition}\label{prop:Schauder}
(i) For $\alpha\in[0,1)$ we have
\begin{equation}\label{eq:calHCalp}
\Vert\mathcal H\psi\Vert_{C^\alpha(\Omega)}\lesssim_A\Vert\psi\Vert_{C^\alpha(\Gamma)}.
\end{equation}
(ii) For $k\geq2$ we have
\begin{equation}\label{eq:DiriCkalp}
\Vert v\Vert_{C^{k,\delta}(\Omega)}\lesssim_A\Vert\Delta v\Vert_{C^{k-2,\delta}(\Omega)}+\Vert v\Vert_{C^{k,\delta}(\Gamma)}+\Vert\Gamma\Vert_{C^{k,\delta}}\Vert v\Vert_{W^{1,\infty}(\Omega)}. 
\end{equation}
\end{proposition}
\proofheading{Proof}(i) is Corollary 5.16 in \cite{IE}. (ii) is a generalization of Proposition 5.14 in \cite{IE} and the proof is similar. \proofqed

\begin{proposition}\label{prop:LapDiri}
(i) Suppose that $v$ solves the Dirichlet problem
\[\left\{\begin{aligned}
\Delta v&=g\quad\text{in }\Omega,\\
v&=\psi\quad\text{on }\Gamma,
\end{aligned}\right.\]
and let $s\geq2$. Then for $r\geq0$, $\alpha,\beta\in[0,1]$ and every sequence of partitions $v=v_j^1+v_j^2$, we have
\begin{equation}\label{eq:LapDiri}
\begin{aligned}
\Vert v\Vert_{H^s(\Omega)}\lesssim_A\Vert g\Vert_{H^{s-2}(\Omega)}+\Vert \psi\Vert_{H^{s-\frac{1}{2}}(\Gamma)}+\Vert\Gamma\Vert_{H^{s+r-\frac{1}{2}}}\sup_{j>0}2^{-j(r+\alpha-1)}\Vert v_j^1\Vert_{C^\alpha(\Omega)}\\
+\sup_{j>0}2^{j(s-1+\beta-\frac{\delta}{2})}\Vert v_j^2\Vert_{H^{1-\beta}(\Omega)}.
\end{aligned}
\end{equation}
(ii) For $k\geq1$, $s\geq0$, $r\geq0$, $\alpha,\beta\in[0,1]$ and partitions $v=v_j^1+v_j^2$ we have
\begin{equation}\label{eq:LDk}
\begin{aligned}
\Vert v\Vert_{H^{s+2k}(\Omega)}\lesssim_M\Vert\Delta^kv\Vert_{H^s(\Omega)}+\sum_{l=1}^k\Vert\Delta^{k-l}v\Vert_{H^{s+2l-\frac{1}{2}}(\Gamma)}+\Vert\Gamma\Vert_{H^{s+2k+r-\frac{1}{2}}}\sup_{j>0}2^{-j(r+\alpha-1)}\Vert v_j^1\Vert_{C^\alpha(\Omega)}\\
+\sup_{j>0}2^{j(s+2k-1+\beta-\frac{\delta}{2})}\Vert v_j^2\Vert_{H^{1-\beta}(\Omega)}.
\end{aligned}
\end{equation}
\end{proposition}
\proofheading{Proof}(i) is Proposition 5.19 in \cite{IE}. \\
(ii) The case $k=1$ is just (i). Suppose the conclusion is true with $k$ replaced by $k-1$. Set $\tilde v=\Delta^{k-1}v=\tilde v_j^1+\tilde v_j^2$, where
\[
\tilde v_j^1=\Delta^{k-1}P_{<j}v_j^1,\quad \tilde v_j^2=\Delta^{k-1}P_{<j}v_j^1+\Delta^{k-1}P_{\geq j}v. 
\]
Bernstein's estimate gives 
\[
2^{-j(r+2k-2+\alpha-1)}\Vert\tilde v_j^1\Vert_{C^\alpha(\Omega)}\lesssim_A2^{-j(r+\alpha-1)}\Vert v_j^1\Vert_{C^\alpha(\Omega)},\quad\forall j,
\]
\[
2^{j(s+1+\beta-\frac{\delta}{2})}\Vert\tilde v_j^2\Vert_{H^{1-\beta}(\Omega)}\lesssim_A2^{j(s+2k-1+\beta-\frac{\delta}{2})}\Vert v_j^2\Vert_{H^{1-\beta}(\Omega)}+\Vert v\Vert_{H^{s+2k-\frac{\delta}{4}}(\Omega)},\quad\forall j. 
\]
Then \eqref{eq:LapDiri} implies
\[\begin{aligned}
\Vert\tilde v\Vert_{H^{s+2}(\Omega)}&\lesssim_A\Vert v\Vert_{H^{s+2k-\frac{\delta}{4}}(\Omega)}+\Vert\Delta^kv\Vert_{H^s(\Omega)}+\Vert\Delta^{k-1}v\Vert_{H^{s+\frac{3}{2}}(\Gamma)}\\ &\quad+\sup_{j>0}2^{j(s+2k-1+\beta-\frac{\delta}{2})}\Vert v_j^2\Vert_{H^{1-\beta}(\Omega)}\\
&\quad+\Vert\Gamma\Vert_{H^{s+2k+r-\frac{1}{2}}}\sup_{j>0}2^{-j(r+\alpha-1)}\Vert v_j^1\Vert_{C^\alpha(\Omega)}.
\end{aligned}\]
Now by \eqref{eq:LDk} with $k,s$ replaced by $k-1,s+2$ and by interpolation, the desired estimate follows. \proofqed

\begin{proposition}\label{prop:kappa}
Let $\kappa$ be the mean curvature of $\Gamma$. For $s\geq2$ we have
\[
\Vert\Gamma\Vert_{H^s}+\Vert n\Vert_{H^{s-1}(\Gamma)}\lesssim_A 1+\Vert \kappa\Vert_{H^{s-2}(\Gamma)}. 
\]
Conversely, we have
\[
\Vert \kappa\Vert_{H^{s-2}(\Gamma)}\lesssim_A\Vert\Gamma\Vert_{H^s}. 
\]
\end{proposition}
\proofheading{Proof}The first statement is Proposition 5.22 in \cite{IE}, and the second statement follows from the explicit formula for mean curvature in local coordinates as well as tame estimates on $\mathbb R^2$. \proofqed

\begin{proposition}\label{prop:LapGam}
Let $\Delta_\Gamma$ be the Laplacian operator on $\Gamma$. (i) For $s\geq2$, $r\geq0$, $\alpha,\beta\in[0,1]$ and any sequence of partitions $v=v_j^1+v_j^2$, we have
\begin{equation}\label{eq:LapGam}
\begin{aligned}
\Vert v\Vert_{H^s(\Gamma)}\lesssim_A\Vert v\Vert_{L^2(\Gamma)}+\Vert\Delta_\Gamma v\Vert_{H^{s-2}(\Gamma)}+\Vert\Gamma\Vert_{H^{s+r}}\sup_{j>0}2^{-j(r+\alpha-1)}\Vert v_j^1\Vert_{C^\alpha(\Gamma)}\\
+\sup_{j>0}2^{j(s-1+\beta-\frac{\delta}{2})}\Vert v_j^2\Vert_{H^{1-\beta}(\Gamma)}.
\end{aligned}
\end{equation}
(ii) For $k\geq1$, $s\geq0$, $r\geq0$, $\alpha,\beta\in[0,1]$ and any sequence of partitions $v=v_j^1+v_j^2$, we have
\begin{equation}\label{eq:LGk}
\begin{aligned}
\Vert v\Vert_{H^{s+2k}(\Gamma)}\lesssim_A\Vert v\Vert_{L^2(\Gamma)}+\Vert\Delta_\Gamma^kv\Vert_{H^s(\Gamma)}+\Vert\Gamma\Vert_{H^{s+2k+r}}\sup_{j>0}2^{-j(r+\alpha-1)}\Vert v_j^1\Vert_{C^\alpha(\Gamma)}\\
+\sup_{j>0}2^{j(s+2k-1+\beta-\frac{\delta}{2})}\Vert v_j^2\Vert_{H^{1-\beta}(\Gamma)}.
\end{aligned}
\end{equation}
\end{proposition}
\proofheading{Proof}(i) Let $u=(\alpha_{*i}v)\circ H_i$ and define $u_j^1,u_j^2$ similarly. We clearly have $\Vert u_j^1\Vert_{C^\alpha(\mathbb R^2)}\lesssim_A\Vert v_j^1\Vert_{C^\alpha(\Gamma)}$ and $\Vert u_j^2\Vert_{H^{1-\beta}(\mathbb R^2)}\lesssim_A\Vert v_j^2\Vert_{H^{1-\beta}(\Gamma)}$. In local coordinates, we have $\Delta_\Gamma=g^{ij}\partial_i\partial_j+h^k\partial_k$, where $(g^{ij})=g^{-1}$ is a smooth function of $\partial H_i$ in $B_i(\frac{5}{4}r_i)$ and $h^k$ is trilinear in $(g^{-1},g^{-1},\partial g)$. Since $u$ is supported in $B_i(\frac{5}{4}r_i)$, $g$ can be extended to the whole of $\mathbf R^2$ such that $g^{ij}=\delta^{ij}$ outside $B_i(2r_i)$ and
\[
\Vert g\Vert_{H^{s+r-1}(\mathbb R^2)}+\Vert g^{-1}\Vert_{H^{s+r-1}(\mathbb R^2)}\lesssim_A\Vert\Gamma\Vert_{H^{s+r}}. 
\]
Separate the constant-coefficient Laplacian from $g^{ij}\partial_i\partial_j$ and apply the Euclidean elliptic estimate:
\begin{equation}\label{eq:flatellipest}
\Vert u\Vert_{H^s(B)}\lesssim\Vert u\Vert_{L^2(B)}+\Vert\Delta_\Gamma u\Vert_{H^{s-2}(B)}+\Vert(g^{-1}-\text{Id})\partial^2u\Vert_{H^{s-2}(B)}+\Vert g^{-2}\partial g\partial u\Vert_{H^{s-2}(B)}. 
\end{equation}
To control the first term on the r.h.s., we use Bony decomposition:
\[\begin{aligned}
\tilde P_j[(g^{-1}-\text{Id})\partial^2u]&=P_{<j}(g^{-1}-\text{Id})P_j\partial^2u+P_j(g^{-1}-\text{Id})P_{<j}\partial^2u\\ &\quad+\sum_{k\geq j}\sum_{|k'-k|\leq1}P_k(g^{-1}-\text{Id})P_{k'}\partial^2u\\
&:=T_j^1+T_j^2+R_j,
\end{aligned}\]
where $P_j$ is the frequency localization operator on $\mathbb R^2$ and $\tilde P_j$ is a fattened version. By adjusting the collar parameters $\delta_0,\sigma$, we may assume that $\Vert g^{-1}-\text{Id}\Vert_{L^\infty}$ is sufficiently small. Hence, we can estimate $T_j^1$, $R_j$ by 
\[
\Vert T_j^1\Vert_{L^2}\leq \frac{1}{100}2^{2j}\Vert P_ju\Vert_{L^2},
\]
\[
\Vert R_j\Vert_{L^2}\lesssim_A\sum_{k\geq j}2^{-k(s-2+\delta)}\Big(2^{ks}\Vert P_ku\Vert_{L^2}\Big), 
\]
and after $\ell^2_j$ summation, 
\[
\Vert2^{j(s-2)}\Vert T_j^1\Vert_{L^2}\Vert_{\ell^2_j}\leq \frac{1}{10}\Vert u\Vert_{H^s},
\]
\[ 
\Vert2^{j(s-2)}\Vert R_j\Vert_{L^2}\Vert_{\ell^2_j}\leq C(A)N2^{2N}\Vert u\Vert_{L^2}+\frac{1}{10}\Vert u\Vert_{H^s}, 
\]
where in deriving the second inequality we split the summation into $j<N$ and $j\geq N$ for sufficiently large $N$. For $T_j^2$ we use that $u=u_j^1+u_j^2$: 
\[
T_j^2=P_j(g^{-1}-\text{Id})P_{<j}\partial^2u_j^1+P_j(g^{-1}-\text{Id})P_{<j}\partial^2u_j^2. 
\]
By Bernstein's estimate we have
\[
\Vert P_j(g^{-1}-\text{Id})P_{<j}\partial^2u_j^1\Vert_{L^2}\lesssim 2^{j(2-\alpha)}\Vert P_j(g^{-1}-\text{Id})\Vert_{L^2}\Vert u_j^1\Vert_{C^\alpha}, 
\]
\[
\Vert P_j(g^{-1}-\text{Id})P_{<j}\partial^2u_j^2\Vert_{L^2}\lesssim 2^{j(1+\beta-\delta)}\Vert u_j^2\Vert_{H^{1-\beta}},
\]
and after $\ell^2_j$ summation, 
\[\begin{aligned}
\Vert2^{j(s-2)}\Vert T_j^2\Vert_{L^2}\Vert_{\ell^2_j}\lesssim_A\Vert g^{-1}\Vert_{H^{s+r-1}}\sup_{j>0}2^{-j(r+\alpha-1)}\Vert u_j^1\Vert_{C^\alpha}+\sup_{j>0}2^{j(s-1+\beta-\frac{\delta}{2})}\Vert u_j^2\Vert_{H^{1-\beta}}.
\end{aligned}\]
This estimates the first coefficient-error term in \eqref{eq:flatellipest}. For the second, expand the three factors by frequency. The interactions split into terms with two lower frequencies and one higher frequency, and terms with one lower frequency and two comparable higher frequencies. At frequency $2^j$, write
\[
P_j(g^{-2}\partial g\partial u)=\widetilde T_j^1+\widetilde T_j^2+\widetilde T_j^3+\widetilde R_j^1+\widetilde R_j^2+\widetilde R_j^3, 
\]
where
\[
\widetilde T_j^1\sim P_{<j-1}(g^{-2})\,P_{<j-1}\partial g\,P_j\partial u,\quad \widetilde T_j^2\sim P_{<j-1}(g^{-2})\,P_j\partial g\,P_{<j-1}\partial u,
\]
\[
\widetilde T_j^3\sim P_j(g^{-2})\,P_{<j-1}\partial g\,P_{<j-1}\partial u,
\]
and
\[
\widetilde R_j^1\sim \sum_{k>j}P_{<j+2}(g^{-2})\,P_k\partial g\,P_k\partial u,\quad \widetilde R_j^2\sim \sum_{k>j}P_k(g^{-2})\,P_{<j+2}\partial g\,P_k\partial u,
\]
\[
\widetilde R_j^3\sim \sum_{k>j}P_k(g^{-2})\,P_k\partial g\,P_{<j+2}\partial u. 
\]
Terms $\widetilde T_j^1,\widetilde R_j^1,\widetilde R_j^2$ are treated like $T_j^1,R_j$: 
\[
\Vert\widetilde T_j^1\Vert_{L^2}\lesssim_A2^{j(2-\delta)}\Vert P_ju\Vert_{L^2},
\]
\[
\Vert\widetilde R_j^1\Vert_{L^2}+\Vert\widetilde R_j^2\Vert_{L^2}\lesssim_A\sum_{k>j}2^{k(2-\delta)}\Vert P_ku\Vert_{L^2}. 
\]
Terms $\widetilde T_j^2,\widetilde T_j^3,\widetilde R_j^3$ are treated like $T_j^2$: 
\[\begin{aligned}
\Vert\widetilde T_j^2\Vert_{L^2}+\Vert\widetilde T_j^3\Vert_{L^2}\lesssim_A&2^{j(2-\alpha)}\Vert u_j^1\Vert_{C^\alpha}(\Vert P_jg\Vert_{L^2}+\Vert P_j(g^{-2})\Vert_{L^2})+2^{j(1-\delta+\beta)}\Vert u_j^2\Vert_{H^{1-\beta}},
\end{aligned}\]
\[
\Vert\widetilde R_j^3\Vert_{L^2}\lesssim_A\sum_{k>j}\Big(2^{k(2-\alpha-\delta)}\Vert u_k^1\Vert_{C^\alpha}+2^{k(1+\beta-\delta)}\Vert u_k^2\Vert_{H^{1-\beta}}\Big). 
\]
Sum in $\ell^2_j$, treating $j\leq N$ and $j>N$ separately, to obtain
\[
\big\Vert2^{j(s-2)}\Vert(\widetilde T_j^1,\widetilde R_j^1,\widetilde R_j^2)\Vert_{L^2}\big\Vert_{\ell^2_j}\leq C(A)N2^{2N}\Vert u\Vert_{L^2}+\frac{1}{10}\Vert u\Vert_{H^s},
\]
\[\begin{aligned}&
\big\Vert2^{j(s-2)}\Vert(\widetilde T_j^2,\widetilde T_j^3,\widetilde R_j^3)\Vert_{L^2}\big\Vert_{\ell^2_j}\\ &\quad\lesssim_A\Vert(g^{-1},g)\Vert_{H^{s+r-1}}\sup_{j>0}2^{-j(r+\alpha-1)}\Vert u_j^1\Vert_{C^\alpha}\\ &\quad+\sup_{j>0}2^{j(s-1+\beta-\frac{\delta}{2})}\Vert u_j^2\Vert_{H^{1-\beta}}. 
\end{aligned}\]
Combining the estimates, we complete the proof of (i). \\
(ii) %follows from an induction on $k$ as we have done for \eqref{eq:LDk}. 
One can adjust the above proof to obtain \eqref{eq:LGk}. \proofqed

\section{Boundary lifting with prescribed normal derivatives}\label{app:extension_lifting}

We prove the boundary lifting estimate from Proposition~\ref{prop:Ej}. Work first on the product manifold $\mathbb R\times\Gamma_*$. For integer $r$, define $H^r(\mathbb R\times\Gamma_*)$ as the completion of test functions under the norm
\[
\Vert F\Vert_{H^r(\mathbb R\times\Gamma_*)}^2=\int\sum_{j=0}^r\Vert\partial_s^jF(s,\cdot)\Vert_{H^{r-j}(\Gamma_*)}^2\,ds
\]
for $r\in\mathbb N$, and is defined via interpolation for fractional $r=j+\alpha$: 
\[
H^{j+\alpha}(\mathbb R\times\Gamma_*)=\left(H^j(\mathbb R\times\Gamma_*),\ H^{j+1}(\mathbb R\times\Gamma_*)\right)_{\alpha,2},\quad\alpha\in(0,1). 
\]
We prove that
\begin{equation}\label{eq:nebd}
\Vert\chi F(\gamma^{-1})\Vert_{H^r(\Omega_*)}\lesssim\Vert F\Vert_{H^r(\mathbb R\times\Gamma_*)}. 
\end{equation}
It suffices to treat integer $r$ and then interpolate. Choose local charts $H_i:B_i\to U_i$ covering $\Gamma_*$, with each $B_i$ a disk, and a subordinate partition of unity $\chi_i$. Set $G_i(s,y)=\chi_i(H_i(y))F(s,H_i(y))$ on $\mathbb R\times B_i$. The equivalence of $W^{r,2}$ and $H^r$ on these cylinders gives
\[
\Vert F\Vert_{H^r(\mathbb R\times\Gamma_*)}\sim\sum_i\Vert G_i\Vert_{H^r(\mathbb R\times B_i)}. 
\]
Next, define $\gamma_i(s,y)=\gamma(s,H_i(y))$, which maps $[0,\delta_0)\times B_i$ diffeomorphically onto a subset of $\Omega_*$. Moser's estimate gives
\[
\Vert\chi G_i(\gamma_i^{-1})\Vert_{H^r(\Omega_*)}\lesssim\Vert G_i\Vert_{H^r(\mathbb R\times B_i)}. 
\]
Finally, by $\chi F(\gamma^{-1})=\sum_i\chi G_i(\gamma_i^{-1})$ we conclude \eqref{eq:nebd}. 

Let $\hat F(\sigma,x)$ denote the Fourier transform of $F$ in the variable $s$: $\hat F(\sigma,x)=\int e^{-2\pi i\sigma s}F(s,x)\,ds$. We claim that
\begin{equation}\label{eq:Hrequiv}
\Vert F\Vert_{H^r(\mathbb R\times\Gamma_*)}^2\sim\int\Vert(1-\Delta_{\Gamma_*}+\sigma^2)^\frac{r}{2}\hat F(\sigma,\cdot)\Vert_{L^2(\Gamma_*)}^2\,d\sigma,\quad r\geq0. 
\end{equation}
For even integer $r$, this follows from the identity $(\partial_s^jG)\hat\ =(2\pi i\sigma)^j\hat G$, the operator interpolation
\[
\Vert(1-\Delta_{\Gamma_*}+\sigma^2)^\frac{r}{2}f\Vert_{L^2(\Gamma_*)}\sim\sum_{j=0}^r\Vert\sigma^j(1-\Delta_{\Gamma_*})^\frac{r-j}{2}f\Vert_{L^2(\Gamma_*)},
\]
the ellipticity
\[
\Vert(1-\Delta_{\Gamma_*})^\frac{r-j}{2}f\Vert_{L^2(\Gamma_*)}\sim\Vert f\Vert_{H^{r-j}(\Gamma_*)},
\]
and the Parseval identity. For real $r$, the claim follows from interpolation. 

We now prove Proposition~\ref{prop:Ej}. Direct differentiation verifies the prescribed boundary jets. With $F(s,\cdot)=\frac{s^j}{j!}\exp\left(-s^{2K}(1-\Delta_{\Gamma_*})^K\right)f$, it remains to show
\begin{equation}\label{eq:hatFbd}
\int\Vert(1-\Delta_{\Gamma_*}+\sigma^2)^{\frac{r+j}{2}+\frac{1}{4}}\hat F(\sigma,\cdot)\Vert_{L^2(\Gamma_*)}^2\,d\sigma\lesssim\Vert(1-\Delta_{\Gamma_*})^\frac{r}{2}f\Vert_{L^2(\Gamma_*)}^2. 
\end{equation}
Apply the spectral theorem to $1-\Delta_{\Gamma_*}$. Let $m$ be a measure on $[1,\infty)$ and $U:L^2(\Gamma_*)\to L^2([1,\infty),m)$ a unitary representation such that
\[
U\sqrt{1-\Delta_{\Gamma_*}}f=\alpha Uf(\alpha). 
\]
Then \eqref{eq:hatFbd} is reduced to 
\begin{equation}\label{eq:tildeFbd}
\int\int_{[1,\infty)}(\alpha^2+\sigma^2)^{r+j+\frac{1}{2}}|\tilde F(\sigma,\alpha)|^2\,dm(\alpha)\,d\sigma\lesssim\int_{[1,\infty)}\alpha^{2r}|Uf(\alpha)|^2\,dm(\alpha),
\end{equation}
where $\tilde F(\sigma,\alpha)=U\hat F(\sigma,\cdot)=M(\sigma,\alpha)Uf(\alpha)$ and
\[
M(\sigma,\alpha)=\int e^{-2\pi i\sigma s}\frac{s^j}{j!}\exp\left(-s^{2K}\alpha^{2K}\right)\,ds. 
\]
By a change of variables $t=\alpha s$ we see
\begin{equation}\label{eq:M1}
|M(\sigma,\alpha)|\lesssim\alpha^{-1-j}. 
\end{equation}
Integrating by parts $N$ times: 
\[
M(\sigma,\alpha)=\frac{1}{(2\pi i\sigma)^N}\int e^{-2\pi i\sigma s}\partial_s^N\left(\frac{s^j}{j!}\exp\left(-s^{2K}\alpha^{2K}\right)\right)\,ds,
\]
and using again the variable substitution $t=\alpha s$, we see
\[
|M(\sigma,\alpha)|\lesssim\alpha^{N-1-j}|\sigma|^{-N}. 
\]
A geometric mean of this and \eqref{eq:M1} is $|M(\sigma,\alpha)|\lesssim\alpha^{r+1}|\sigma|^{-j-r-2}$. Putting these together we obtain
\begin{equation}\label{eq:M2}
|M(\sigma,\alpha)|\lesssim\alpha^{r+1}(\alpha+|\sigma|)^{-j-r-2}. 
\end{equation}
Consider the integral
\[
\int|M(\sigma,\alpha)|^2(\alpha^2+\sigma^2)^{r+j+\frac{1}{2}}\,d\sigma=\left(\int_{|\cdot|\leq\alpha}+\int_{|\cdot|>\alpha}\right)|M(\sigma,\alpha)|^2(\alpha^2+\sigma^2)^{r+j+\frac{1}{2}}\,d\sigma:=I^0+I^1. 
\]
From \eqref{eq:M1}, \eqref{eq:M2} it follows
\[
|I^0|\lesssim\alpha^{2r},\quad|I^1|\lesssim\alpha^{2r}. 
\]
Substituting these bounds into \eqref{eq:tildeFbd} completes the proof. \proofqed

\end{document}